\documentclass[11pt,twoside]{amsart} 
\usepackage{amssymb}
\usepackage{amscd}
\usepackage[matrix,arrow]{xy}

    \newtheorem{thm}{Theorem}                     [section]
    \newtheorem{thm*}{Theorem}
    \newtheorem{prop}[thm]{Proposition}
    \newtheorem{lemma}[thm]{Lemma}
    \newtheorem{cor}[thm]{Corollary}

    \newtheorem{lemma*}{Lemma}    

    \newtheorem{defn}[thm]{Definition}                 
    
    \newtheorem{rems*}{Remark}   

\newcommand{\ndef}{\newcommand*}
\def\rndef{\renewcommand}

\ndef{\myaddress}[1]{\begin{center} \it\small #1 \end{center}}

\ndef{\clA}{{\mathcal A}} \ndef{\rmA}{{\mathrm A}} \ndef{\mbA}{{\mathbb A}} \ndef{\bfA}{{\mathbf A}} \ndef{\euA}{{\EuScript A}} \ndef{\frA}{{\mathfrak A}}
\ndef{\clB}{{\mathcal B}} \ndef{\rmB}{{\mathrm B}} \ndef{\mbB}{{\mathbb B}} \ndef{\bfB}{{\mathbf B}} \ndef{\euB}{{\EuScript B}} \ndef{\frB}{{\mathfrak B}}
\ndef{\clC}{{\mathcal C}} \ndef{\rmC}{{\mathrm C}} \ndef{\mbC}{{\mathbb C}} \ndef{\bfC}{{\mathbf C}} \ndef{\euC}{{\EuScript C}} \ndef{\frC}{{\mathfrak C}}
\ndef{\clD}{{\mathcal D}} \ndef{\rmD}{{\mathrm D}} \ndef{\mbD}{{\mathbb D}} \ndef{\bfD}{{\mathbf D}} \ndef{\euD}{{\EuScript D}} \ndef{\frD}{{\mathfrak D}}
\ndef{\clE}{{\mathcal E}} \ndef{\rmE}{{\mathrm E}} \ndef{\mbE}{{\mathbb E}} \ndef{\bfE}{{\mathbf E}} \ndef{\euE}{{\EuScript E}} \ndef{\frE}{{\mathfrak E}}
\ndef{\clF}{{\mathcal F}} \ndef{\rmF}{{\mathrm F}} \ndef{\mbF}{{\mathbb F}} \ndef{\bfF}{{\mathbf F}} \ndef{\euF}{{\EuScript F}} \ndef{\frF}{{\mathfrak F}}
\ndef{\clG}{{\mathcal G}} \ndef{\rmG}{{\mathrm G}} \ndef{\mbG}{{\mathbb G}} \ndef{\bfG}{{\mathbf G}} \ndef{\euG}{{\EuScript G}} \ndef{\frG}{{\mathfrak G}}
\ndef{\clH}{{\mathcal H}} \ndef{\rmH}{{\mathrm H}} \ndef{\mbH}{{\mathbb H}} \ndef{\bfH}{{\mathbf H}} \ndef{\euH}{{\EuScript H}} \ndef{\frH}{{\mathfrak H}}
\ndef{\clI}{{\mathcal I}} \ndef{\rmI}{{\mathrm I}} \ndef{\mbI}{{\mathbb I}} \ndef{\bfI}{{\mathbf I}} \ndef{\euI}{{\EuScript I}} \ndef{\frI}{{\mathfrak I}}
\ndef{\clJ}{{\mathcal J}} \ndef{\rmJ}{{\mathrm J}} \ndef{\mbJ}{{\mathbb J}} \ndef{\bfJ}{{\mathbf J}} \ndef{\euJ}{{\EuScript J}} \ndef{\frJ}{{\mathfrak J}}
\ndef{\clK}{{\mathcal K}} \ndef{\rmK}{{\mathrm K}} \ndef{\mbK}{{\mathbb K}} \ndef{\bfK}{{\mathbf K}} \ndef{\euK}{{\EuScript K}} \ndef{\frK}{{\mathfrak K}}
\ndef{\clL}{{\mathcal L}} \ndef{\rmL}{{\mathrm L}} \ndef{\mbL}{{\mathbb L}} \ndef{\bfL}{{\mathbf L}} \ndef{\euL}{{\EuScript L}} \ndef{\frL}{{\mathfrak L}}
\ndef{\clM}{{\mathcal M}} \ndef{\rmM}{{\mathrm M}} \ndef{\mbM}{{\mathbb M}} \ndef{\bfM}{{\mathbf M}} \ndef{\euM}{{\EuScript M}} \ndef{\frM}{{\mathfrak M}}
\ndef{\clN}{{\mathcal N}} \ndef{\rmN}{{\mathrm N}} \ndef{\mbN}{{\mathbb N}} \ndef{\bfN}{{\mathbf N}} \ndef{\euN}{{\EuScript N}} \ndef{\frN}{{\mathfrak N}}
\ndef{\clO}{{\mathcal O}} \ndef{\rmO}{{\mathrm O}} \ndef{\mbO}{{\mathbb O}} \ndef{\bfO}{{\mathbf O}} \ndef{\euO}{{\EuScript O}} \ndef{\frO}{{\mathfrak O}}
\ndef{\clP}{{\mathcal P}} \ndef{\rmP}{{\mathrm P}} \ndef{\mbP}{{\mathbb P}} \ndef{\bfP}{{\mathbf P}} \ndef{\euP}{{\EuScript P}} \ndef{\frP}{{\mathfrak P}}
\ndef{\clQ}{{\mathcal Q}} \ndef{\rmQ}{{\mathrm Q}} \ndef{\mbQ}{{\mathbb Q}} \ndef{\bfQ}{{\mathbf Q}} \ndef{\euQ}{{\EuScript Q}} \ndef{\frQ}{{\mathfrak Q}}
\ndef{\clR}{{\mathcal R}} \ndef{\rmR}{{\mathrm R}} \ndef{\mbR}{{\mathbb R}} \ndef{\bfR}{{\mathbf R}} \ndef{\euR}{{\EuScript R}} \ndef{\frR}{{\mathfrak R}}
\ndef{\clS}{{\mathcal S}} \ndef{\rmS}{{\mathrm S}} \ndef{\mbS}{{\mathbb S}} \ndef{\bfS}{{\mathbf S}} \ndef{\euS}{{\EuScript S}} \ndef{\frS}{{\mathfrak S}}
\ndef{\clT}{{\mathcal T}} \ndef{\rmT}{{\mathrm T}} \ndef{\mbT}{{\mathbb T}} \ndef{\bfT}{{\mathbf T}} \ndef{\euT}{{\EuScript T}} \ndef{\frT}{{\mathfrak T}}
\ndef{\clU}{{\mathcal U}} \ndef{\rmU}{{\mathrm U}} \ndef{\mbU}{{\mathbb U}} \ndef{\bfU}{{\mathbf U}} \ndef{\euU}{{\EuScript U}} \ndef{\frU}{{\mathfrak U}}
\ndef{\clV}{{\mathcal V}} \ndef{\rmV}{{\mathrm V}} \ndef{\mbV}{{\mathbb V}} \ndef{\bfV}{{\mathbf V}} \ndef{\euV}{{\EuScript V}} \ndef{\frV}{{\mathfrak V}}
\ndef{\clW}{{\mathcal W}} \ndef{\rmW}{{\mathrm W}} \ndef{\mbW}{{\mathbb W}} \ndef{\bfW}{{\mathbf W}} \ndef{\euW}{{\EuScript W}} \ndef{\frW}{{\mathfrak W}}
\ndef{\clX}{{\mathcal X}} \ndef{\rmX}{{\mathrm X}} \ndef{\mbX}{{\mathbb X}} \ndef{\bfX}{{\mathbf X}} \ndef{\euX}{{\EuScript X}} \ndef{\frX}{{\mathfrak X}}
\ndef{\clY}{{\mathcal Y}} \ndef{\rmY}{{\mathrm Y}} \ndef{\mbY}{{\mathbb Y}} \ndef{\bfY}{{\mathbf Y}} \ndef{\euY}{{\EuScript Y}} \ndef{\frY}{{\mathfrak Y}}
\ndef{\clZ}{{\mathcal Z}} \ndef{\rmZ}{{\mathrm Z}} \ndef{\mbZ}{{\mathbb Z}} \ndef{\bfZ}{{\mathbf Z}} \ndef{\euZ}{{\EuScript Z}} \ndef{\frZ}{{\mathfrak Z}}

\ndef{\tA}{{\widetilde A}} \ndef{\tcA}{{\widetilde\clA}} \ndef{\ttcA}{\widetilde{\tcA}} \ndef{\sfA}{{\textsf A}} \ndef{\ttA}{\widetilde{\tA}} \ndef{\dzA}{{A^\sharp}}
\ndef{\tB}{{\widetilde B}} \ndef{\tcB}{{\widetilde\clB}} \ndef{\ttcB}{\widetilde{\tcB}} \ndef{\sfB}{{\textsf B}} \ndef{\ttB}{\widetilde{\tB}} \ndef{\dzB}{{B^\sharp}}
\ndef{\tC}{{\widetilde C}} \ndef{\tcC}{{\widetilde\clC}} \ndef{\ttcC}{\widetilde{\tcC}} \ndef{\sfC}{{\textsf C}} \ndef{\ttC}{\widetilde{\tC}} \ndef{\dzC}{{C^\sharp}}
\ndef{\tD}{{\widetilde D}} \ndef{\tcD}{{\widetilde\clD}} \ndef{\ttcD}{\widetilde{\tcD}} \ndef{\sfD}{{\textsf D}} \ndef{\ttD}{\widetilde{\tD}} \ndef{\dzD}{{D^\sharp}}
\ndef{\tE}{{\widetilde E}} \ndef{\tcE}{{\widetilde\clE}} \ndef{\ttcE}{\widetilde{\tcE}} \ndef{\sfE}{{\textsf E}} \ndef{\ttE}{\widetilde{\tE}} \ndef{\dzE}{{E^\sharp}}
\ndef{\tF}{{\widetilde F}} \ndef{\tcF}{{\widetilde\clF}} \ndef{\ttcF}{\widetilde{\tcF}} \ndef{\sfF}{{\textsf F}} \ndef{\ttF}{\widetilde{\tF}} \ndef{\dzF}{{F^\sharp}}
\ndef{\tG}{{\widetilde G}} \ndef{\tcG}{{\widetilde\clG}} \ndef{\ttcG}{\widetilde{\tcG}} \ndef{\sfG}{{\textsf G}} \ndef{\ttG}{\widetilde{\tG}} \ndef{\dzG}{{G^\sharp}}
\ndef{\tH}{{\widetilde H}} \ndef{\tcH}{{\widetilde\clH}} \ndef{\ttcH}{\widetilde{\tcH}} \ndef{\sfH}{{\textsf H}} \ndef{\ttH}{\widetilde{\tH}} \ndef{\dzH}{{H^\sharp}}
\ndef{\tI}{{\widetilde I}} \ndef{\tcI}{{\widetilde\clI}} \ndef{\ttcI}{\widetilde{\tcI}} \ndef{\sfI}{{\textsf I}} \ndef{\ttI}{\widetilde{\tI}} \ndef{\dzI}{{I^\sharp}}
\ndef{\tJ}{{\widetilde J}} \ndef{\tcJ}{{\widetilde\clJ}} \ndef{\ttcJ}{\widetilde{\tcJ}} \ndef{\sfJ}{{\textsf J}} \ndef{\ttJ}{\widetilde{\tJ}} \ndef{\dzJ}{{J^\sharp}}
\ndef{\tK}{{\widetilde K}} \ndef{\tcK}{{\widetilde\clK}} \ndef{\ttcK}{\widetilde{\tcK}} \ndef{\sfK}{{\textsf K}} \ndef{\ttK}{\widetilde{\tK}} \ndef{\dzK}{{K^\sharp}}
\ndef{\tL}{{\widetilde L}} \ndef{\tcL}{{\widetilde\clL}} \ndef{\ttcL}{\widetilde{\tcL}} \ndef{\sfL}{{\textsf L}} \ndef{\ttL}{\widetilde{\tL}} \ndef{\dzL}{{L^\sharp}}
\ndef{\tM}{{\widetilde M}} \ndef{\tcM}{{\widetilde\clM}} \ndef{\ttcM}{\widetilde{\tcM}} \ndef{\sfM}{{\textsf M}} \ndef{\ttM}{\widetilde{\tM}} \ndef{\dzM}{{M^\sharp}}
\ndef{\tN}{{\widetilde N}} \ndef{\tcN}{{\widetilde\clN}} \ndef{\ttcN}{\widetilde{\tcN}} \ndef{\sfN}{{\textsf N}} \ndef{\ttN}{\widetilde{\tN}} \ndef{\dzN}{{N^\sharp}}
\ndef{\tO}{{\widetilde O}} \ndef{\tcO}{{\widetilde\clO}} \ndef{\ttcO}{\widetilde{\tcO}} \ndef{\sfO}{{\textsf O}} \ndef{\ttO}{\widetilde{\tO}} \ndef{\dzO}{{O^\sharp}}
\ndef{\tP}{{\widetilde P}} \ndef{\tcP}{{\widetilde\clP}} \ndef{\ttcP}{\widetilde{\tcP}} \ndef{\sfP}{{\textsf P}} \ndef{\ttP}{\widetilde{\tP}} \ndef{\dzP}{{P^\sharp}}
\ndef{\tQ}{{\widetilde Q}} \ndef{\tcQ}{{\widetilde\clQ}} \ndef{\ttcQ}{\widetilde{\tcQ}} \ndef{\sfQ}{{\textsf Q}} \ndef{\ttQ}{\widetilde{\tQ}} \ndef{\dzQ}{{Q^\sharp}}
\ndef{\tR}{{\widetilde R}} \ndef{\tcR}{{\widetilde\clR}} \ndef{\ttcR}{\widetilde{\tcR}} \ndef{\sfR}{{\textsf R}} \ndef{\ttR}{\widetilde{\tR}} \ndef{\dzR}{{R^\sharp}}
\ndef{\tS}{{\widetilde S}} \ndef{\tcS}{{\widetilde\clS}} \ndef{\ttcS}{\widetilde{\tcS}} \ndef{\sfS}{{\textsf S}} \ndef{\ttS}{\widetilde{\tS}} \ndef{\dzS}{{S^\sharp}}
\ndef{\tT}{{\widetilde T}} \ndef{\tcT}{{\widetilde\clT}} \ndef{\ttcT}{\widetilde{\tcT}} \ndef{\sfT}{{\textsf T}} \ndef{\ttT}{\widetilde{\tT}} \ndef{\dzT}{{T^\sharp}}
\ndef{\tU}{{\widetilde U}} \ndef{\tcU}{{\widetilde\clU}} \ndef{\ttcU}{\widetilde{\tcU}} \ndef{\sfU}{{\textsf U}} \ndef{\ttU}{\widetilde{\tU}} \ndef{\dzU}{{U^\sharp}}
\ndef{\tV}{{\widetilde V}} \ndef{\tcV}{{\widetilde\clV}} \ndef{\ttcV}{\widetilde{\tcV}} \ndef{\sfV}{{\textsf V}} \ndef{\ttV}{\widetilde{\tV}} \ndef{\dzV}{{V^\sharp}}
\ndef{\tW}{{\widetilde W}} \ndef{\tcW}{{\widetilde\clW}} \ndef{\ttcW}{\widetilde{\tcW}} \ndef{\sfW}{{\textsf W}} \ndef{\ttW}{\widetilde{\tW}} \ndef{\dzW}{{W^\sharp}}
\ndef{\tX}{{\widetilde X}} \ndef{\tcX}{{\widetilde\clX}} \ndef{\ttcX}{\widetilde{\tcX}} \ndef{\sfX}{{\textsf X}} \ndef{\ttX}{\widetilde{\tX}} \ndef{\dzX}{{X^\sharp}}
\ndef{\tY}{{\widetilde Y}} \ndef{\tcY}{{\widetilde\clY}} \ndef{\ttcY}{\widetilde{\tcY}} \ndef{\sfY}{{\textsf Y}} \ndef{\ttY}{\widetilde{\tY}} \ndef{\dzY}{{Y^\sharp}}
\ndef{\tZ}{{\widetilde Z}} \ndef{\tcZ}{{\widetilde\clZ}} \ndef{\ttcZ}{\widetilde{\tcZ}} \ndef{\sfZ}{{\textsf Z}} \ndef{\ttZ}{\widetilde{\tZ}} \ndef{\dzZ}{{Z^\sharp}}

\ndef{\bfc}{{\bf c}}

  \ndef{\eps}{\varepsilon}

\let\geq\geqslant
\let\leq\leqslant

\ndef{\lims}[1]{\lim\limits_{#1}}
\ndef{\sums}[1]{\sum\limits_{#1}}
\ndef{\ints}[1]{\int\limits_{#1}}
\ndef{\sups}[1]{\sup\limits_{#1}}
\ndef{\liminfty}[1]{\lims{#1\to\infty}}
\ndef{\suminf}[1]{\sums{#1=1}^\infty}

\ndef{\limo}[1]{\omega\mbox{-}\!\!\!\lims{#1\to\infty}}          
\ndef{\limL}[1]{\rmL\mbox{-}\!\!\!\lims{#1\to\infty}}            
\ndef{\limLOne}[1]{\clL_1\mbox{-}\!\!\!\lims{#1}}
\ndef{\tildelimo}[1]{\tilde\omega\mbox{-}\!\!\!\lims{#1\to\infty}}

\ndef{\Aut}{\operatorname{Aut}}      
\ndef{\Ch}{\operatorname{ch}}        
\ndef{\End}{\operatorname{End}}      
\ndef{\Hom}{\operatorname{Hom}}      
\ndef{\Ker}{\operatorname{Ker}}      
\ndef{\Log}{\operatorname{Log}}      
\ndef{\OP}{\operatorname{OP}}        
\ndef{\Op}{\operatorname{Op}}        
\ndef{\Symb}{\operatorname{Symb}}    
\ndef{\Tr}{\operatorname{Tr}}        
\ndef{\Wres}{\operatorname{Wres}}    
\ndef{\cl}{\operatorname{cl}}        
\ndef{\com}{\operatorname{com}}
\ndef{\const}{\operatorname{const}}  
\ndef{\conv}{\operatorname{conv}}    
\rndef{\det}{\operatorname{det}}     

\ndef{\detFK}[1]{\Delta\brs{#1}} 
\ndef{\detFKrel}[2]{\Delta_{#2}\brs{#1}} 

\ndef{\diag}{\operatorname{diag}}    
\ndef{\dist}{\operatorname{dist}}    
\ndef{\dom}{\operatorname{dom}}      
\ndef{\ec}{\operatorname{ec}}        
\ndef{\id}{1}                        
\ndef{\ind}{\operatorname{ind}}      
\ndef{\mydeg}{\operatorname{deg}}    
\ndef{\op}{\operatorname{op}}
\ndef{\rank}{\operatorname{rank}}
\ndef{\res}{\operatorname{res}}      
\ndef{\rng}{\operatorname{ran}}      
\ndef{\sflow}{\operatorname{sf}}     
\ndef{\isf}{\operatorname{isf}}      
\ndef{\sign}{\operatorname{sign}}    
\ndef{\sgn}{\operatorname{sgn}}      
\ndef{\sing}{\operatorname{sing}}    
\ndef{\supp}{\operatorname{supp}}    
\ndef{\tr}{\operatorname{tr}}        
\ndef{\var}{\operatorname{var}}      
\ndef{\vol}{\operatorname{vol}}      
\ndef{\wn}{\operatorname{wn}}        
\ndef{\wres}{\operatorname{wres}}    
\rndef{\Im}{\operatorname{Im}}       
\rndef{\Re}{\operatorname{Re}}       

\ndef{\prng}[1]{\mathrm R_{#1}} 
\ndef{\pker}[1]{\mathrm N_{#1}} 
\ndef{\rprng}[2]{\mathrm R_{#1}^{#2}}           
\ndef{\rpker}[2]{\mathrm N_{#1}^{#2}}           
\ndef{\rsupp}[1]{\supp_r(#1)}
\ndef{\lsupp}[1]{\supp_l(#1)}
\ndef{\rslv}[1]{R_z(#1)}      
\ndef{\HH}{H}                 
\ndef{\tHH}{\tilde \HH}       
\ndef{\VV}{V}                 
\ndef{\Rz}{R_z}               
\ndef{\tRz}{\tR_z}            
\ndef{\psif}[1]{#1^{[1]}} 
\ndef{\CPlus}[1]{W_{#1}(\mbR)}
\ndef{\bndl}{\xi}                         
\ndef{\bndlA}{\eta}                       
\ndef{\GlueMap}{\varphi}                  
\ndef{\ChartMap}{h}                       

\ndef{\hilb}{\clH}                     
   \ndef{\hilbasargument}{(\hilb)} 
\ndef{\LpH}[1]{\clL^{#1}\hilbasargument}       
\ndef{\clBH}{\clB\hilbasargument}              
\ndef{\ubBH}{\clB_1\hilbasargument}            
\ndef{\clCH}{\clC\hilbasargument}              
\ndef{\clKH}{\clK\hilbasargument}              
\ndef{\clFH}{\clF\hilbasargument}              
\ndef{\clUH}{\clU\hilbasargument}              
\ndef{\clCFH}{{\clC\clF}\hilbasargument}       
\ndef{\saBH}{\clB_{sa}\hilbasargument}         
\ndef{\saCH}{\clC_{sa}\hilbasargument}         
\ndef{\saFH}{\clF_{sa}\hilbasargument}         
\ndef{\saKH}{\clK_{sa}\hilbasargument}         
\ndef{\saCFH}{\clC\clF_{sa}\hilbasargument}    
\ndef{\clUFH}{\clU\clF\hilbasargument}         
\ndef{\Uinj}{\clU_{inj}\hilbasargument}        
\ndef{\UFinj}{\clU\clF_{inj}\hilbasargument}   

\ndef{\spproj}[2]{E^{#1}_{#2}}                      
\ndef{\spprojb}[2]{E^{#2}_{#1}}                     

\ndef{\LpN}[1]{\clL^{#1}(\clN,\tau)}     
\ndef{\saLpN}[1]{\clL^{#1}_{sa}(\clN,\tau)} 
\ndef{\rLpN}[1]{L^{#1}(\clN,\tau)}       
\ndef{\clAND}{(\clA,\clN,D)}             
\ndef{\clBA}{{\clB(\clA)}}
\ndef{\saKN}{{\clK_{sa}(\clN,\tau)}}          
\ndef{\clKN}{{\clK(\clN,\tau)}}          
\ndef{\clKtN}{{\clK(\tilde\clN,\tau)}}   
\ndef{\clFN}{{\clF(\clN,\tau)}}          
\ndef{\saFN}{{\clF_{sa}(\clN,\tau)}}     
\ndef{\clPN}{\clP(\clN)}                 
\ndef{\clQN}{\clQ(\clN,\tau)}            
\ndef{\infPN}{{\clP_\tau^\infty(\clN)}}  
\ndef{\clOF}[2]{\clF_{#1\mbox{-}#2}(\clN,\tau)}         
\ndef{\oind}[2]{{\rm \tau\mbox{-}ind}_{#1\mbox{-}#2}}   
\ndef{\tind}{\tau\mbox{-}\ind}                  
\ndef{\DInd}{\ind_{\clD,\tau}}           
\ndef{\BF}{Breuer-Fredholm}              
\ndef{\skewfred}[2]{$(#1\cdot #2)$ $\tau$\tire Fredholm}   
\ndef{\affl}{\eta}                       
\ndef{\vNa}{von Neumann algebra}         
\ndef{\nsf}{faithful normal semifinite } 
\ndef{\taubrs}[1]{\tau\brackets{#1}}     
\ndef{\sqbrs}[1]{\left[#1\right]}        

\ndef{\domd}{\bigcap\limits_{n\ge 0} \dom\;\delta^n}         
\ndef{\DiffOP}{{\rm \clD}}
\ndef{\ADA}{\clA \cup [\clD,\clA]}
\ndef{\DixIdeal}[1]{\LpH{#1,\infty}}               
\ndef{\dixideal}{\ell^{1,\infty}}                  
\ndef{\WDixIdeal}{\LpH{1,\mathrm w}}               
\ndef{\DixIdealPos}[1]{\DixIdeal{#1}_+}            
\ndef{\DixIdealN}[1]{\LpN{#1,\infty}}              
\ndef{\DixIdealNPar}[2]{\clL^{#1,\infty}_{#2}(\clN,\tau)}    
\ndef{\DixIdealNPos}[1]{\LpN{#1,\infty}_+}                   
\ndef{\TrD}{\Tr_\omega}                                      
\ndef{\tauD}{{\tau_\omega}}                                  
\ndef{\ILog}{\frac 1{\log(1+t)}}
\ndef{\ILogN}{\frac 1{\log(1+N)}}
\ndef{\DixNorm}[1]{\norm{#1}_{(1,\infty)}}                   
\ndef{\DixInt}[1]{\ints 0^t \mu_s(#1)\,ds}
\ndef{\DixIntL}[1]{\ints 0^{\lambda_{1/t}(#1)}\mu_s(#1)\,ds}
    \ndef{\SmallIdeal}{{\clL^{1, \mathrm w}}}
    \ndef{\SmallIdealMeas}{{\clL^{1, \mathrm w}_m}}
    \ndef{\DixIntII}[1]{\ints 0^t \mu_s(#1)\,ds}
    \ndef{\DixIntf}[1]{f_t(#1)}
    \ndef{\DixIntg}[1]{g_t(#1)}

\ndef{\lpi}{\clL^{1,\pi}(\clN,\tau)}
\ndef{\IIinfty}{$\mathrm{II}_\infty$\ }

\ndef{\fourier}[1]{\clF(#1)}          
\ndef{\HaarMeasBohrs}{\nu}            
\ndef{\BrownsMeas}{\mu}               
\ndef{\BohrCont}[1]{\tilde{#1}}       
\ndef{\APMean}{{M}}                   
\ndef{\CDSS}{{\clA_B}}                
\ndef{\matr}{{\rm Mat}}               
\ndef{\seque}[1]{\ensuremath{\{#1_j\}_{j=1}^\infty}}    
\ndef{\sequen}[2]{\ensuremath{\{#1_#2\}_{#2=1}^\infty}}    
\ndef{\Seque}[1]{\ensuremath{\left(#1_0,#1_1,#1_2,\dots\right)}}    
\ndef{\Cesaro}{H}                           
\ndef{\CesaroRPlus}{M}                      
\ndef{\Dilation}{D}                         
\ndef{\Shift}{T}                            

\ndef{\norm}[1]{\left\Vert#1\right\Vert}    
\ndef{\TrNorm}[1]{\norm{#1}_1}              
\ndef{\HSNorm}[1]{\norm{#1}_2}              
\ndef{\InftyNorm}[1]{\norm{#1}_\infty}      
\ndef{\normQN}[1]{\norm{#1}_{\clQN}}        
\ndef{\clLnorm}[1]{\norm{#1}_{1,\infty}}    

\ndef{\ccurve}{\gamma}                      

\ndef{\abs}[1]{\left\lvert#1\right\rvert}   
\ndef{\set}[1]{\left\{#1\right\}}           
\ndef{\brackets}[1]{\left(#1\right)}        
\ndef{\brs}[1]{\brackets{#1}}               
\ndef{\Brs}[1]{\big(#1\big)}                
\ndef{\BRS}[1]{\Big(#1\Big)}                
\ndef{\scal}[2]{\left\la #1,#2\right\ra}               
\ndef{\precprec}{\prec\!\!\!\prec}
\ndef{\qeq}{\stackrel?=}
\ndef{\spectrum}[1]{\sigma_{#1}} 
\ndef{\numrange}[1]{\mathrm{W}(#1)}                         
\rndef{\emptyset}{\varnothing}                              
\ndef{\csupp}{c}                           
\ndef{\closure}[1]{\overline{#1}}
\ndef{\linspan}[1]{\mathrm{span}\ {#1}}
\ndef{\bddborel}[1]{B(#1)}                 
\ndef{\charfunc}{\chi}
\rndef{\ln}{\log}
\ndef{\FrDer}{\euD}                        
\ndef{\LieDer}[1]{\pounds_{#1}\,}          
\ndef{\dds}{\left.\frac d{ds} \right|_{s = 0}}
\ndef{\ortcmp}[1]{#1^{\scriptscriptstyle \perp}}            
\ndef{\Laplace}{\Delta}                    

\ndef{\matrPQ}[3]
{
    \left(
      \begin{array}{cc}
        #1_{11} & #1_{12} \\
        #1_{21} & #1_{22}
      \end{array}
    \right)_{[#2,#3]}
}

\newcounter{margcomcount}
\ndef{\margcom}[1]{\marginpar{\bf \small #1} \addtocounter{margcomcount}{1}}

\newcounter{margproof}
\ndef{\margproof}{\marginpar{\bf \small Proof} \addtocounter{margproof}{1}}

\newcounter{margdetails}
\ndef{\margdetails}{\marginpar{\bf \small has been changed} \addtocounter{margdetails}{1}}

\newcounter{margproofb}
\ndef{\margproofb}{\marginpar{\bf \small Proof (B)} \addtocounter{margproof}{1}}

\newcounter{margdetailsb}
\ndef{\margdetailsb}{\marginpar{\bf \small Details (B)} \addtocounter{margdetailsb}{1}}

\ndef{\mytimes}{\!\times\!}
\ndef{\sss}[1]{\subsubsection{}\label{#1}}
\rndef{\phi}{\varphi}
\ndef{\OpenUnitDisk}{D}
\ndef{\RHS}{RHS}                            
\ndef{\LHS}{LHS} 
\ndef{\ttt}{\Leftrightarrow}
\ndef{\then}{\Rightarrow}
\ndef{\tto}{\longrightarrow}
\ndef{\nno}{\nonumber\\}
\ndef{\newn}[1]{\index{#1} \emph{#1}}       
\ndef{\la}{\langle}
\ndef{\ra}{\rangle}
\ndef{\dbar}{{\;\bar{\phantom{o}} \!\!\!\! d}}
\ndef{\stl}[1]{\stackrel{\vbox to 0pt{\vss\hbox{$\scriptstyle #1$}}}}
\ndef{\mathcomment}[1]{{\scriptstyle\text{(#1)}}\qquad}        
\ndef{\details}[1]{\smallskip\begin{center} {\bf Here:} #1\end{center}\medskip}
\ndef{\indexcom}[1]{ --- #1}
\ndef{\longsim}{\ \sim \ }              
\ndef{\tire}{-}              
\ndef{\intinfinf}{\int_{-\infty}^\infty}
\ndef{\refnsftrace}{\cite[V.\,2.\,1]{TakI}} 
\ndef{\refaffloper}{\cite[IV.\,5, Exercise 3]{TakI}} 
\ndef{\refsemifinvNa}{\cite[V.\,1.\,21]{TakI}} 
\ndef{\reftaumeasurable}{\cite[Definition 1.2]{FK86PJM}} 
\ndef{\reftautraceclassaffl}{\cite[V.2, p.\,320]{TakI}} 
\ndef{\refinvoperideal}{\cite[Appendix A.2]{CP2}} 
\ndef{\reftautracenorm}{\cite[V.2, p.\,320]{TakI}} 
\ndef{\reftaucompact}{\cite{}} 
\ndef{\reftauFredholm}{\cite[Appendix B]{PR94JFA}} 

     \ndef{\npartial}{\slash\!\!\!\partial}
     \ndef{\Heis}{\operatorname{Heis}}
     \ndef{\Solv}{\operatorname{Solv}}
     \ndef{\Spin}{\operatorname{Spin}}
     \ndef{\SO}{\operatorname{SO}}
     \ndef{\Index}{\operatorname{index}}

             \ndef{\coker}{{\mbox coker}}
             \ndef{\p}{\partial}
             \ndef{\dd}{|\clD|}
             \ndef{\n}{\parallel}

     \ndef{\gf}[2]{\genfrac{}{}{0pt}{}{#1}{#2}}
     \ndef{\ta}{\widetilde{\alpha}}
     \ndef{\tb}{\widetilde{\beta}}
     \ndef{\txi}{\widetilde{\xi}}
     \ndef{\tk}{\widetilde{K}}
     \ndef{\CGh}{\widetilde{\CG}}
     \ndef{\boe}{{\bf e}}\ndef{\bt}{{\bf t}}
     \ndef{\vth}{\vartheta}
     \ndef{\db}{\overline{\partial}}
     \ndef{\hV}{\hat{V}}
     \ndef{\cag}{{\clA^\Gamma}}
     \ndef{\sind}{\sigma{\rm -ind}}

\let\LatexCite=\cite  

\let\ifnumref\iftrue 

\ndef{\ifuncited}[4]{\expandafter\ifx\csname used#4\endcsname\relax}

\ndef{\ifcited}[4]{\expandafter\ifx\csname used#4\endcsname\relax\else}

  \ndef{\papertitle}[1]{ \emph{#1}, }
  \ndef{\paperauthor}[2]{#2}  
  \ndef{\pbbi}[9]{%
      \ifcited{#1}{#2}{#3}{#5}%
        \ifnumref%
          \bibitem{#5}\paperauthor{#1}{#6},\papertitle{#7}#8.%
        \else%
          \advance #9 by 1%
          \ifnum#9<1%
            \bibitem[#4]{#5}\paperauthor{#1}{#6}, \papertitle{#7}#8.%
          \else%
            \bibitem[#4$\!_{\the#9}\!$]{#5}\paperauthor{#1}{#6},\papertitle{#7}#8.%
          \fi%
        \fi%
      \fi%
  }
  \ndef{\mbbi}[8]{%
     \ifcited{#1}{#2}{#3}{#5}%
        \ifnumref%
          \bibitem{#5}\paperauthor{#1}{#6},\papertitle{#7}#8.%
        \else%
          \bibitem[#4]{#5}\paperauthor{#1}{#6},\papertitle{#7}#8.%
        \fi%
     \fi%
  }

\ndef{\AddCite}[1]{%
   \ifuncited{0}{0}{0}{#1}%
     \expandafter\gdef\csname used#1\endcsname {}%
   \fi%
}

\def\ProcessCite#1,{%
     \ifx\relax#1%
         \let\next=\relax%
     \else%
         \AddCite{#1}%
         \let\next=\ProcessCite%
     \fi%
     \next%
}

\ndef{\AddCites}[1]{\ProcessCite#1,\relax,}

\ndef{\CiteWithoutExtension}[1]{%
   \AddCites{#1}%
   \LatexCite{#1}%
}

\def\CiteWithExtension[#1]#2{%
   \AddCites{#2}%
   \LatexCite[#1]{#2}%
}

\ndef{\CleverCite}{%
    \ifx\NChar[ %
       \let\MyCite=\CiteWithExtension %
    \else %
       \let\MyCite=\CiteWithoutExtension %
    \fi %
    \MyCite%
}

\renewcommand{\cite}{\futurelet\NChar\CleverCite}

      \ndef{\volume}[1]{{\bf #1}}
      \ndef{\VolYearPP}[3]{\ifnum#2=0 (to appear)\else\volume{#1} (#2), #3\fi}
      \ndef{\VolNoYearPP}[4]{\ifnum#3=0 (to appear)\else\volume{#1} #2 (#3), #4\fi}
      \ndef{\libcode}[1]{}

\ndef{\jnActaMath}[3]{Acta Math. \VolYearPP{#1}{#2}{#3}}                       
\ndef{\jnAdvMath}[3]{Adv. in~Math. \VolYearPP{#1}{#2}{#3}}                     
\ndef{\jnAlgAnal}[3]{Algebra i~Analiz \VolYearPP{#1}{#2}{#3}}
\ndef{\jnAmerMathMonth}[3]{Amer. Math. Monthly \VolYearPP{#1}{#2}{#3}}         
\ndef{\jnAnnMath}[4]{Ann. of~Math. \VolNoYearPP{#1}{#2}{#3}{#4}}               
\ndef{\jnAnalMath}[3]{J. Anal. Math. \VolYearPP{#1}{#2}{#3}}                   
\ndef{\jnBullLondMathSoc}[3]{Bull. London Math. Soc. \VolYearPP{#1}{#2}{#3}}   
\ndef{\jnBullAMS}[3]{Bull. Amer. Math. Soc. \VolYearPP{#1}{#2}{#3}}   
\ndef{\jnCanMathBull}[3]{Canad. Math. Bull. \VolYearPP{#1}{#2}{#3}}            
\ndef{\jnCanMath}[4]{Canad. J.~Math. \VolYearPP{#1}{#2}{#3}}             
\ndef{\jnCommMathPhys}[3]{Comm. Math. Phys \VolYearPP{#1}{#2}{#3}}             
\ndef{\jnCommPDE}[3]{Comm. Partial Differential Equations \VolYearPP{#1}{#2}{#3}}             
\ndef{\jnComptRendue}[3]{C.\,R.~Acad. Sci. Paris S\'er. A-B \VolYearPP{#1}{#2}{#3}}      
\ndef{\jnDiffGeom}[3]{J.~Diff. Geom. \VolYearPP{#1}{#2}{#3}}                   
\ndef{\jnErgodicTheory}[3]{Ergodic Theory and Dynamical Systems \VolYearPP{#1}{#2}{#3}} 
\ndef{\jnFuncAnal}[3]{J.~Functional Analysis \VolYearPP{#1}{#2}{#3}}           
\ndef{\jnFunkAnalPril}[4]{Функциональный анализ и его приложения \VolNoYearPP{#1}{#2}{#3}{#4}}  
\ndef{\jnGAFA}[3]{GAFA \VolYearPP{#1}{#2}{#3}}                                 
\ndef{\jnIHES}[3]{IHES Publ. Math. (Paris) \VolYearPP{#1}{#2}{#3}}             
\ndef{\jnIEOT}[3]{Integral Equations Operator Theory   \VolYearPP{#1}{#2}{#3}} 
\ndef{\jnIsrMath}[3]{Israel J.~Math. \VolYearPP{#1}{#2}{#3}}                   
\ndef{\jnKTheory}[3]{K-Theory \VolYearPP{#1}{#2}{#3}}                          
\ndef{\jnLetMathPhys}[3]{Lett. Math. Phys. \VolYearPP{#1}{#2}{#3}}             
\ndef{\jnMathAnn}[3]{Math. Ann. \VolYearPP{#1}{#2}{#3}}                        
\ndef{\jnMathAnalAppl}[3]{J.~Math. Anal. and Appl. \VolYearPP{#1}{#2}{#3}}     
\ndef{\jnMathNachr}[3]{Math. Nachr. \VolYearPP{#1}{#2}{#3}}
\ndef{\jnMathPhys}[3]{J. Math. Phys. \VolYearPP{#1}{#2}{#3}}
\ndef{\jnOperTheory}[3]{J.~Operator Theory \VolYearPP{#1}{#2}{#3}}             
\ndef{\jnPacJMath}[3]{Pacific J.~Math. \VolYearPP{#1}{#2}{#3}}                  
\ndef{\jnPositivity}[3]{Positivity \VolYearPP{#1}{#2}{#3}}
\ndef{\jnProcAmerMS}[3]{Proc. Amer. Math. Soc. \VolYearPP{#1}{#2}{#3}}         
\ndef{\jnProcCambPhilSoc}[3]{Math. Proc. Camb. Phil. Soc. \VolYearPP{#1}{#2}{#3}}
\ndef{\jnReineAngew}[3]{J.~Reine Angew. Math. \VolYearPP{#1}{#2}{#3}}          
\ndef{\jnTokyoMath}[3]{Tokyo J.~Math. \VolYearPP{#1}{#2}{#3}}
\ndef{\jnTopology}[3]{Topology \VolYearPP{#1}{#2}{#3}}
\ndef{\jnTransAmerMathSoc}[3]{Trans. Amer. Math. Soc. \VolYearPP{#1}{#2}{#3}}
\ndef{\jnIzvANSSSR}[3]{Izv. Akad. Nauk SSSR, Ser. Mat. \VolYearPP{#1}{#2}{#3}}
\ndef{\jnIzvVyshUchZav}[3]{Izv. Vyssh. Uch. Zav., Mat. \VolYearPP{#1}{#2}{#3} (Russian)}
\ndef{\jnIzdatLenUniv}[2]{Izdat. Leningrad. Univ., Leningrad, (#1), #2 (Russian)}
\ndef{\jnFieldsInsComm}[3]{Fields Inst. Comm. \VolYearPP{#1}{#2}{#3}}
\ndef{\jnDoklANSSSR}[3]{Dokl. Akad. Nauk SSSR \VolYearPP{#1}{#2}{#3}}
\ndef{\jnMatZametki}[3]{Matem. zametki \VolYearPP{#1}{#2}{#3}}
\ndef{\jnRussMathSurvey}[3]{Russian Math. Surveys \VolYearPP{#1}{#2}{#3}}
\ndef{\jnSibMathJ}[3]{Sib. Math.~J. \VolYearPP{#1}{#2}{#3}}
\ndef{\jnSovMath}[3]{J.~Soviet math. \VolYearPP{#1}{#2}{#3}}
\ndef{\jnTransMoscMathSoc}[3]{Trans. Moscow Math. Soc. \VolYearPP{#1}{#2}{#3}}
\ndef{\jnUMN}[3]{Uspekhi Mat. Nauk \VolYearPP{#1}{#2}{#3}}

\ndef{\bkTransMathMon}[2]{Trans. Math. Monographs, AMS, \volume{#1}, #2}

\ndef{\pbBirkhauser}[1]{Birkh\"auser, Boston, #1}
\ndef{\pbFactorial}[1]{Moscow, Factorial, #1}
\ndef{\pbGauthier}[1]{Gauthier-Villars, Paris, #1}
\ndef{\pbNauka}[1]{Moscow, Nauka, #1 (Russian)}
\ndef{\pbNaukaR}[1]{Москва, Наука, #1}
\ndef{\pbPrinceton}[1]{Princeton University Press, Princeton, New Jersey, #1}
\ndef{\pbPublPerish}[1]{Publish or Perish Inc., Berkeley, #1}
\ndef{\pbSpringer}[1]{Springer-Verlag, #1}

\ndef{\myauthor}[1]{\mbox{#1}}

\ndef{\Ahiezer}{\myauthor{N.\,I.\,Ahiezer }}
\ndef{\Arazy}{\myauthor{J.\,Arazy}}
\ndef{\Astashkin}{\myauthor{S.\,V.\,Astashkin}}
\ndef{\Atiyah}{\myauthor{M.\,Atiyah}}
\ndef{\Avron}{\myauthor{J.\,Avron}}
\ndef{\Azamov}{\myauthor{N.\,A.\,Azamov}}
\ndef{\Banach}{\myauthor{S.\,Banach}}
\ndef{\Benameur}{\myauthor{M-T.\,Benameur}}
\ndef{\Bennett}{\myauthor{C.\,Bennett}}
\ndef{\Berezin}{\myauthor{F.\,A.\,Berezin}}
\ndef{\Berline}{\myauthor{N.\,Berline}}
\ndef{\Birman}{\myauthor{M.\,Sh.\,Birman}}
\ndef{\Blackadar}{\myauthor{B.\,Blackadar}}
\ndef{\Bogolyubov}{\myauthor{N.\,N.\,Bogolyubov}}
\ndef{\Bonsall}{\myauthor{F.\,F.\,Bonsall}}
\ndef{\BoosBavnbek}{\myauthor{B.\,Boo$\beta$-Bavnbek}}
\ndef{\Bott}{\myauthor{R.\,Bott}}
\ndef{\Bratteli}{\myauthor{O.\,Bratteli}}
\ndef{\Bredon}{\myauthor{G.\,E.\,Bredon}}
\ndef{\Breuer}{\myauthor{M.\,Breuer}}
\ndef{\Brown}{\myauthor{L.\,G.\,Brown}}
\ndef{\Bruneau}{\myauthor{V.\,Bruneau}}
\ndef{\Buslaev}{\myauthor{V.\,S.\,Buslaev}}
\ndef{\Carey}{\myauthor{A.\,L.\,Carey}}
\ndef{\CareyRW}{\myauthor{R.\,W.\,Carey}} 
\ndef{\Cartan}{\myauthor{H.\,Cartan}}
\ndef{\Chilin}{\myauthor{V.\,I.\,Chilin}}
\ndef{\Coburn}{\myauthor{L.\,A.\,Coburn}}
\ndef{\Connes}{\myauthor{A.\,Connes}}
\ndef{\Cornfeld}{\myauthor{I.\,P.\,Cornfeld}}
\ndef{\Daletskii}{\myauthor{Yu.\,L.\,Daletski\u\i}}   
\ndef{\Dixmier}{\myauthor{J.\,Dixmier}}
\ndef{\DoddsPG}{\myauthor{P.\,G.\,Dodds}}
\ndef{\DoddsTK}{\myauthor{T.\,K.\,Dodds}}
\ndef{\Douglas}{\myauthor{R.\,G.\,Douglas}}
\ndef{\Dubrovin}{\myauthor{B.\,A.\,Dubrovin}}
\ndef{\Dugundji}{\myauthor{J.\,Dugundji}}
\ndef{\Duncan}{\myauthor{J.\,Duncan}}
\ndef{\Dunford}{\myauthor{N.\,Dunford}}
\ndef{\Dykema}{\myauthor{K.\,J.\,Dykema}}
\ndef{\Edwards}{\myauthor{R.\,E.\,Edwards}}
\ndef{\Eilenberg}{\myauthor{S.\,Eilenberg}}
\ndef{\Fack}{\myauthor{T.\,Fack}} 
\ndef{\Faddeev}{\myauthor{L.\,D.\,Faddeev}}
\ndef{\Farber}{\myauthor{M.\,Farber}}
\ndef{\Farforovskaya}{\myauthor{Yu.\,B.\,Farforovskaya}}
\ndef{\Federer}{\myauthor{H.\,Federer}}
\ndef{\Fedosov}{\myauthor{B.\,V.\,Fedosov}}
\ndef{\Figiel}{\myauthor{T.\,Figiel}} 
\ndef{\Figueroa}{\myauthor{H.\,Figueroa}}
\ndef{\Fillmore}{\myauthor{P.\,A.\,Fillmore}}
\ndef{\Fomenko}{\myauthor{A.\,T.\,Fomenko}} 
\ndef{\Fomin}{\myauthor{S.\,V.\,Fomin}}
\ndef{\Frohlich}{\myauthor{J.\,Fr\"ohlich}}
\ndef{\Fuglede}{\myauthor{B.\,Fuglede}}
\ndef{\Furutani}{\myauthor{K.\,Furutani}}
\ndef{\Gelfand}{\myauthor{I.\,M.\,Gelfand}}
\ndef{\Gesztesy}{\myauthor{F.\,Gesztesy}}     
\ndef{\Getzler}{\myauthor{E.\,Getzler}} 
\ndef{\Gilkey}{\myauthor{P.\,B.\,Gilkey}}
\ndef{\Gitler}{\myauthor{S.\,Gitler}}
\ndef{\Glazman}{\myauthor{I.\,M.\,Glazman}}
\ndef{\Glimm}{\myauthor{J.\,Glimm}}
\ndef{\Gohberg}{\myauthor{I.\,C.\,Gohberg}}
\ndef{\Golze}{\myauthor{F.\,Golze}}
\ndef{\GraciaBondia}{\myauthor{J.\,M.\,Gracia-Bond\'{i}a}}
\ndef{\Greenleaf}{\myauthor{F.\,P.\,Greenleaf}}
\ndef{\Gromov}{\myauthor{M.\,Gromov}}
\ndef{\Gunning}{\myauthor{R.\,C.\,Gunning}}
\ndef{\Haagerup}{\myauthor{U.\,Haagerup}}
\ndef{\Haag}{\myauthor{R.\,Haag}}
\ndef{\Halmos}{\myauthor{Halmos}}
\ndef{\Hardy}{\myauthor{G.\,H.\,Hardy}}
\ndef{\Higson}{\myauthor{N.\,Higson}}  
\ndef{\Hoermander}{\myauthor{L.\,Hoermander}} 
\ndef{\Hoffman}{\myauthor{K.\,Hoffman}} 
\ndef{\Ito}{\myauthor{K.\,Ito}}
\ndef{\Jaffe}{\myauthor{A.\,Jaffe}}
\ndef{\James}{\myauthor{I.\,M.\,James}}
\ndef{\Javrjan}{\myauthor{V.\,A.\,Javrjan}}
\ndef{\Kadison}{\myauthor{R.\,V.\,Kadison}}
\ndef{\Kalton}{\myauthor{N.\,J.\,Kalton}} 
\ndef{\Kato}{\myauthor{T.\,Kato}} 
\ndef{\Kobayashi}{\myauthor{S.\,Kobayashi}}
\ndef{\Koplienko}{\myauthor{L.\,S.\,Koplienko}}
\ndef{\Korotyaev}{\myauthor{E.\,Korotyaev}}
\ndef{\Kosaki}{\myauthor{H.\,Kosaki}}
\ndef{\KreinMG}{\myauthor{M.\,G.\,Kre\u\i n}}
\ndef{\KreinSG}{\myauthor{S.\,G.\,Kre\u\i n}}
\ndef{\Leichtnam}{\myauthor{E.\,Leichtnam}}
\ndef{\Lesch}{\myauthor{M.\,Lesch}}
\ndef{\Lesniewski}{\myauthor{A.\,Lesniewski}}
\ndef{\Levitan}{\myauthor{B.\,M.\,Levitan}}
\ndef{\Lidskii}{\myauthor{V.\,B.\,Lidskii}}
\ndef{\Lifshits}{\myauthor{I.\,M.\,Lifshits}}
\ndef{\Lindenstrauss}{\myauthor{J.\,Lindenstrauss}}
\ndef{\Loday}{\myauthor{J.-L.\,Loday}}
\ndef{\Lord}{\myauthor{S.\,Lord}}      
\ndef{\Lorentz}{\myauthor{G.\,Lorentz}}
\ndef{\Magnus}{\myauthor{W.\,Magnus}}
\ndef{\Makarov}{\myauthor{K.\,A.\,Makarov}}
\ndef{\Mathai}{\myauthor{V.\,Mathai}}         
\ndef{\McKean}{\myauthor{H.\,P.\,McKean}}
\ndef{\Mishchenko}{\myauthor{A.\,S.\,Mishchenko}}
\ndef{\Moore}{\myauthor{C.\,C.\,Moore}}
\ndef{\Moscovici}{\myauthor{H.\,Moscovici}}  
\ndef{\Motovilov}{\myauthor{A.\,K.\,Motovilov}}
\ndef{\Moyer}{\myauthor{R.\,D.\,Moyer}}
\ndef{\Naboko}{\myauthor{S.\,N.\,Naboko}}
\ndef{\Narasimhan}{\myauthor{R.\,Narasimhan}}
\ndef{\Nomizu}{\myauthor{K.\,Nomizu}}
\ndef{\Novikov}{\myauthor{S.\,P.\,Novikov}}
\ndef{\Osterwalder}{\myauthor{K.\,Osterwalder}}
\ndef{\Patodi}{\myauthor{V.\,Patodi}}
\ndef{\Pagter}{\myauthor{B.\,de~Pagter}}  
\ndef{\Pavlov}{\myauthor{B.\,S.\,Pavlov}}
\ndef{\Pedersen}{\myauthor{G.\,K.\,Pedersen}}
\ndef{\Peller}{\myauthor{V.\,V.\,Peller}}
\ndef{\Perera}{\myauthor{V.\,S.\,Perera}}
\ndef{\Petunin}{\myauthor{Ju.\,I.\,Petunin}}
\ndef{\Phillips}{\myauthor{J.\,Phillips}}  
\ndef{\Piazza}{\myauthor{P.\,Piazza}}   
\ndef{\Pincus}{\myauthor{J.\,D.\,Pincus}}   
\ndef{\Poincare}{Poincar\'e}
\ndef{\Postnikov}{\myauthor{M.\,M.\,Postnikov}} 
\ndef{\Prinzis}{\myauthor{R.\,Prinzis}}
\ndef{\Privalov}{\myauthor{I.\,I.\,Privalov}}
\ndef{\Pushnitski}{\myauthor{A.\,B.\,Pushnitski}} 
\ndef{\Raeburn}{\myauthor{I.\,Raeburn}}
\ndef{\Raikov}{\myauthor{G.\,Raikov}}
\ndef{\Reed}{\myauthor{M.\,Reed}}
\ndef{\Rennie}{\myauthor{A.\,Rennie}}
\ndef{\Rickart}{\myauthor{C.\,E.\,Rickart}}
\ndef{\Riesz}{\myauthor{F.\,Riesz}}
\ndef{\Ringrose}{\myauthor{J.\,Ringrose}}
\ndef{\Robinson}{\myauthor{D.\,Robinson}}
\ndef{\Rossi}{\myauthor{H.\,Rossi}}
\ndef{\Rudin}{\myauthor{W.\,Rudin}}
\ndef{\Ruelle}{\myauthor{D.\,Ruelle}}
\ndef{\Ruzhansky}{\myauthor{M.\,Ruzhansky}}
\ndef{\Sakai}{\myauthor{Sh.\,Sakai}}
\ndef{\Sargsjan}{\myauthor{I.\,S.\,Sargsjan}}
\ndef{\Sato}{\myauthor{H.\,Sato}}
\ndef{\Schaeffer}{\myauthor{D.\,G.\,Schaeffer}}
\ndef{\Schluchtermann}{\myauthor{G.\,Schluchtermann}}
\ndef{\Schochet}{\myauthor{C.\,Schochet}}
\ndef{\Schrodinger}{\myauthor{E.\,Schr\"odinger}}
\ndef{\Schrohe}{\myauthor{E.\,Schrohe}}
\ndef{\Schwartz}{\myauthor{J.\,T.\,Schwartz}}
\ndef{\Sedaev}{\myauthor{A.\,A.\,Sedaev}}
\ndef{\Seiler}{\myauthor{R.\,Seiler}}
\ndef{\Semenov}{\myauthor{E.\,M.\,Semenov}}
\ndef{\Shabat}{\myauthor{B.\,V.\,Shabat}}
\ndef{\Shafarevich}{\myauthor{I.\,R.\,Shafarevich}}
\ndef{\Sharpley}{\myauthor{R.\,Sharpley}}
\ndef{\Shilov}{\myauthor{G.\,E.\,Shilov}}
\ndef{\Shirkov}{\myauthor{D.\,V.\,Shirkov}}
\ndef{\Shubin}{\myauthor{M.\,A.\,Shubin}}
\ndef{\Silverman}{\myauthor{H.\,Silverman}}
\ndef{\Simon}{\myauthor{B.\,Simon}}
\ndef{\Sinai}{\myauthor{Ya.\,G.\,Sinai}}
\ndef{\Singer}{\myauthor{I.\,M.\,Singer}}
\ndef{\Solomyak}{\myauthor{M.\,Z.\,Solomyak}}
\ndef{\Soloviev}{\myauthor{Yu.\,P.\,Soloviev}}
\ndef{\Spivak}{\myauthor{M.\,Spivak}}
\ndef{\Stenkin}{\myauthor{V.\,V.\,Sten'kin}}
\ndef{\Stratila}{\myauthor{S.\,Stratila}}
\ndef{\Sucheston}{\myauthor{L.\,Sucheston}}
\ndef{\Sukochev}{\myauthor{F.\,A.\,Sukochev}}
\ndef{\Switzer}{\myauthor{R.\,M.\,Switzer}}
\ndef{\SzNagy}{\myauthor{B.\,Sz.-Nagy}}
\ndef{\Takesaki}{\myauthor{M.\,Takesaki}}
\ndef{\Taylor}{\myauthor{M.\,E.\,Taylor}}
\ndef{\Treves}{\myauthor{F.\,Treves}}
\ndef{\Troitsky}{\myauthor{E.\,V.\,Troitsky}}
\ndef{\Tzafriri}{\myauthor{L.\,Tzafriri}}
\ndef{\Varilly}{\myauthor{J.\,C.\,V\'{a}rilly}}
\ndef{\Vergne}{\myauthor{M.\,Vergne}}
\ndef{\Vladimirov}{\myauthor{V.\,S.\,Vladimirov}}
\ndef{\Voiculescu}{\myauthor{D.\,Voiculescu}}
\ndef{\Weiss}{\myauthor{G.\,Weiss}}
\ndef{\Wells}{\myauthor{R.\,O.\,Wells}}
\ndef{\Williams}{\myauthor{J.\,P.\,Williams}}
\ndef{\Winkler}{\myauthor{S.\,Winkler}}
\ndef{\Witten}{\myauthor{E.\,Witten}}
\ndef{\Wodzicki}{\myauthor{M.\,Wodzicki}}
\ndef{\Wojciechowski}{\myauthor{K.\,P.\,Wojciechowski}}
\ndef{\Yafaev}{\myauthor{D.\,R.\,Yafaev}}
\ndef{\Yosida}{\myauthor{K.\,Yosida}}
\ndef{\Zsido}{\myauthor{L.\,Zsido}}

\ndef{\epspi}{\sqrt{\frac{\eps}{\pi}}}
\rndef{\clAND}{(\clA,\clN,\HH)}
\rndef{\HH}{D}
\rndef{\VV}{V}
\rndef{\tHH}{\HH_1}
\ndef{\rr}{r}
\rndef{\isf}{\Phi}     
\ndef{\TraceComp}{\clT\clC}
\ndef{\FNab}[2]{\clF^{#1,#2}(\clN,\tau)}
\ndef{\clFNpm}{\FNab {-1}{1}}
\ndef{\clFNab}{\FNab a b}
\ndef{\h}{h}                   
\ndef{\af}{\kappa}
\ndef{\BanLDCT}{Lebesgue dominated convergence theorem~\cite[Corollary III.6.16]{DS}}
\rndef{\CPlus}[1]{C^{#1,+}(\mbR)}
\begin{document}
\begin{center} THE SPECTRAL SHIFT FUNCTION AND SPECTRAL FLOW.
\end{center}
\vskip 1 cm
\begin{center}
   \Azamov$^a$, \ \Carey$^{b},$ \ \Sukochev$^{c}$
\end{center}

\myaddress{$^{a,c}$ School of Informatics and Engineering,
   \\ Flinders University of South Australia,
   \\ Bedford Park, 5042, SA, Australia.
   \\ Email: azam0001,sukochev@infoeng.flinders.edu.au}
\myaddress{$^{b}$ Mathematical Sciences Institute,
   \\ Australian National University,
   \\ Canberra, ACT 0200, Australia,
   \\ Email: acarey@maths.anu.edu.au}
\begin{center} \today\end{center}

 {\bf Abstract. }
At the 1974 International Congress, I. M. Singer proposed that eta invariants and hence
spectral flow should be thought of as the integral of a one form.
In the intervening years this idea has lead to many interesting
developments in the study of both eta invariants and
spectral flow. Using ideas of~\cite{Ge93Top} Singer's proposal
was brought to an advanced level in~\cite{CP2} where a very general
formula for spectral flow as the integral of a one form was produced in the
framework of noncommutative geometry.
This formula can be used for computing spectral flow in a general semifinite
von Neumann algebra as described and reviewed in~\cite{BCPRSW}.
In the present paper we take the analytic approach to spectral flow
much further by giving a large family of formulae for spectral flow
between a pair of unbounded self-adjoint operators~$\HH$ and~$\HH+\VV$
with~$\HH$ having compact resolvent belonging to
a general semifinite von Neumann algebra~$\clN$
and the perturbation $\VV\in \clN.$ In noncommutative geometry terms we remove summability hypotheses.
This level of generality is
made possible by introducing a new idea from~\cite{ACDS}.
There it was observed that M.\,G.\,Krein's
spectral shift function (in certain restricted cases with~$\VV$ trace class)
computes spectral flow. The present paper  extends Krein's theory
to the setting of semifinite spectral triples where~$\HH$
has compact resolvent belonging to~$\clN$
and~$\VV$ is any bounded self-adjoint operator in~$\clN.$
We give a definition of the  spectral shift function under these
hypotheses and show that it computes
spectral flow. This is made possible by the understanding
discovered in the present paper of the interplay between
spectral shift function
theory and the analytic theory of spectral flow.
It is this interplay that enables us to
take Singer's idea much further to create a large class of
one forms whose integrals calculate spectral flow. These advances
depend critically on a new approach to the
calculus of functions of non-commuting
operators discovered in~\cite{ACDS} which
generalizes the double operator integral formalism of~\cite{BS66DOI,BS67DOI,BS73DOI}.
One surprising conclusion that follows from our
results is that the Krein spectral shift function is computed,
in certain circumstances, by the Atiyah-Patodi-Singer index theorem
~\cite{APS76}.

 \vskip 1 cm

\section*{Introduction}

\subsection*{Overview}
In~\cite{ACDS} we gave an analytic proof that
the spectral shift function of M.\,G.\,Krein computes
the spectral flow under certain restricted circumstances.
Spectral flow stems from the work
of Atiyah-Patodi-Singer~\cite{APS76} where it is introduced primarily in
a topological sense. Subsequently, starting with a suggestion of I.\,M.\,Singer at the 1974 Vancouver ICM,
the idea that spectral flow could be expressed as the integral of a one form has been
extensively studied in the framework of unbounded Fredholm modules
(or spectral triples) beginning with~\cite{Ge93Top} and continuing in~\cite{CP98CJM,CP2}.
The current state of knowledge using an analytic approach due
to Phillips~\cite{Ph96CMB,Ph97FIC} is described in detail in~\cite{BCPRSW}.

On the other hand the Krein spectral shift function has an extensive history in perturbation
theory which may be partly traced from~\cite{Kr63,BP98IEOT,BYa92AA,SimTrId}.
The spectral shift theory was developed by both physicists and
mathematicians in the context of perturbation theory
of Schrodinger operators. By
the Birman-Krein formula, it is related to
the phase of the scattering matrix in the latter's spectral
representation. Krein's formula for the spectral
shift function, which is the motivation for this paper,
is restricted to the case
where one perturbs an arbitrary unbounded self adjoint operator $\HH_0$ by a trace class operator.
The spectral shift function compares in a sense the relation between
the spectrum of $\HH_0$ and that of  its perturbation.
In this paper we will extend the theory of the spectral shift function to the situation where we perturb
$\HH_0$  by an arbitrary bounded self adjoint operator but adopt the assumption
that $\HH_0$ has compact resolvent motivated by the notion of  spectral triple.
One surprising conclusion that follows is
that the Atiyah-Patodi-Singer (or APS) index theorem
\cite{APS76} and its generalizations are,
in certain instances, computing the spectral shift function.

The generalizations of the APS index theorem we are referring to in the previous paragraph
have to do with situations where von Neumann algebras
other than all the bounded operators on a Hilbert space arise.
 The standard theory of spectral flow deals with unbounded self
adjoint Fredholm operators. Using an approach due to Phillips
{\em op cit} spectral flow may be defined analytically
(but as yet not topologically)
for certain operators affiliated to semifinite von Neumann algebras,
the so-called Breuer-Fredholm theory, which is reviewed in~\cite{BCPRSW}.
In addition the connection between spectral flow in semifinite spectral triples
and generalizations of the APS index theorem is explained in
\cite{BCPRSW} with reference to the extensive previous history of the
matter.

Spectral flow is related to odd degree $K$-theory and the odd local index theorem in
noncommutative geometry~\cite{CM95GAFA,CPRS2}. For
the purposes of using spectral flow to obtain information about $K$-theory
it is enough to consider
unbounded operators~$\HH$ with compact resolvent. However, previous studies of
spectral flow formulae have restrictions on~$\HH$
such as theta summability (requiring the heat operator $e^{-t\HH^2}$ to be trace class for
all $t>0$). In this paper we will relax this condition to capture the more
general situation where only compact resolvent is needed.
By contrast the theory of the spectral shift function was
formulated for Schrodinger operators associated to Euclidean
spaces and in that setting the assumption of
compact resolvent does not generally hold. Thus it is not surprising
that the relationship between the
Krein spectral shift function and spectral flow
should have remained somewhat unexplored until recently.

In~\cite{ACDS} we noticed that under hypotheses that guaranteed that both the Krein
spectral shift function and spectral flow exist then they are essentially the same notion.
This led us to the current investigation where we borrow some ideas from
spectral flow theory as formulated in a spectral triple
to extend the range of situations for which the Krein spectral
shift function may be defined. We are then able to obtain new results on
the spectral shift function and to prove, by a combination
of the methods of~\cite{ACDS} and~\cite{CP98CJM,CP2}
that the spectral shift function gives a wide variety
of analytic formulae for spectral flow. In fact it is
the interaction between these previously disparate theories that makes these advances possible.

An additional very important feature of our analysis is its generality;
we are able to study self adjoint operators that are affiliated
to a general semifinite von Neumann algebra~$\clN$
in contrast to the classical theory which deals with the special
case where~$\clN$ is the algebra of all bounded operators on
a Hilbert space~$\hilb$. (Note that there are examples of quantum mechanical
Hamiltonians whose resolvent lies in a semifinite von Neumann algebra, see e.g. \cite{CHMM})

\subsection*{Summary of results}
In order to keep this paper to a reasonable length we have omitted a
detailed discussion
of spectral flow in von Neumann algebras referring
to the paper~\cite{BCPRSW} for this.
On the other hand we include in Section \ref{S: preliminaries} preliminary
results particularly those from~\cite{ACDS}
which may be less accessible to the reader.
This latter paper is very important for the present investigation because
it provides
a calculus of functions of operators more effective than has been available in
the past. The key technical advance that we exploit is the use of
double operator integrals (DOI)
not in its original form~\cite{BS66DOI}
but in the much more effective form discovered in~\cite{ACDS}.
In particular it is this latter paper
that develops a calculus for functions of operators affiliated to general
semifinite von Neumann
algebras. This calculus replaces the more elementary
perturbation theory techniques in~\cite{CP98CJM,CP2} and
is more user friendly than that of \cite{dPS04FA}.

The main new results on the Krein spectral shift function
are in Section \ref{S: spectral shift measure}. Here we show
that by starting with self adjoint operators $\HH_0$ with compact
resolvent in a semifinite von Neumann algebra we are
able to define the Krein spectral shift function
for general bounded self adjoint perturbations $\VV.$
Section \ref{S: spectral shift measure}
also exhibits some properties of the spectral shift function for
operators satisfying these hypotheses mostly in the context
of preparing the ground for the subsequent discussion of spectral flow.
Specifically we observe that the
spectral projections $E_{(a,b)}^{\HH_0}$ and  $E_{(a,b)}^{\HH_0 +\VV}$ of $\HH_0$
and $\HH_0+\VV=\HH_1$ respectively are finite in $\clN$ for bounded intervals $(a,b).$
Using a fixed faithful normal semifinite trace $\tau$ on $\clN$
we define (following~\cite{BS72SM}) the spectral shift measure for the pair $\HH_0,\HH_1$, by
\begin{gather*}\label{F: BS formula}
  \Xi_{\HH_1,\HH_0}(\Delta) = \int_0^1 \taubrs{\VV E_{\Delta}^{\HH_\rr}}\,d\rr,
\end{gather*}
where $ \Delta$ is a bounded Borel subset of the real line and $\HH_r = \HH_0 + rV, \ r \in [0,1].$
We consider this Birman-Solomyak formula
(formula (\ref{F: def of Xi}) in the text)
for the spectral shift measure as fundamental. It is this formula which
in our opinion must be taken as the definition of the generalized spectral shift
function of a pair of operators, whenever this
expression makes sense for that pair.
We then prove that this measure is absolutely continuous with respect
to Lebesgue measure and the resulting Radon-Nikodym derivative we
define to be the spectral shift function  $\xi_{\HH_1,\HH_0}(\lambda)$.
This function can be related to the original spectral shift function
of Krein which was introduced in a completely different fashion.
Furthermore, under our assumptions, the spectral shift function
$\xi_{\HH_1,\HH_0}(\lambda)$ is in fact a function of a bounded variation and
there is a canonical representative which is an everywhere defined function.
We remark that in \cite{Si98PAMS} B.\,Simon outlined (without proof) some
similar results on the spectral shift function.

In Section \ref{S: spectral flow} we present a series
of new analytic formulae for spectral flow. These require the use
of the spectral shift function and its properties derived
in Section \ref{S: spectral shift measure}. The formulae we obtain
subsume those of~\cite{CP2}.  We are able to provide a variety
of analytic formulae whenever $\HH_0$ has compact resolvent and show
how these formulae can be specialized when there are summability
hypotheses imposed.

For illustrative purposes we describe two of our main theorems here.
Suppose that there is a unitary operator $u$ with $\HH_1=u\HH_0 u^*$ and such that
$\VV=u[\HH_0,u^*]\in \clN.$ That is, for a dense subalgebra
of the $C^*$-algebra generated by $u$
we have a semifinite spectral triple.
Then spectral flow from $\HH_0-\lambda$ to $\HH_1-\lambda$ for any real $\lambda$ is equal
to  $\xi_{\HH_1,\HH_0}(\lambda)$.
We relate our results to previous formulae for spectral flow via the following result.
For any positive integrable function $f$ on the real line
with $f(\HH_0+rV)$ trace class for $r\in [0,1]$ and  $\rr \mapsto \norm{f(\HH_0+rV)}_1$
being 1-summable on $[0,1],$ we find that spectral flow from
$\HH_0-\lambda$ to $\HH_1-\lambda$ for any real $\lambda$ is given by
$$
  \sflow(\lambda; \HH_0,\HH_1) = C^{-1} \int_0^1 \taubrs{\VV f(\HH_\rr - \lambda)}\,d\rr,
$$
where $C = \int_{-\infty}^\infty f(x)\,dx.$ In fact in this restricted situation of
unitarily equivalent endpoints we can prove that the spectral shift function
is constant implying that spectral flow occurs uniformly past any point
on the real line not just zero. We also obtain analogues of these results
when the endpoint operators are not unitarily equivalent. Then the relation between
spectral flow and the spectral shift function is modified by the addition of
endpoint correction terms (in a fashion analogous to~\cite{CP2}).
The key idea we exploit in this part of the paper is the observation~\cite{Ge93Top}
that spectral flow can be written as the integral of an
exact one form on the affine space of bounded perturbations of a fixed
unbounded operator $\HH_0$. It eventuates that there is a way to construct
a large class of such
one forms that will compute spectral flow when the only constraint
on $\HH_0$ is that it have compact resolvent. This advance is made possible
by the understanding of the spectral shift measure afforded by Section \ref{S: spectral shift measure}.
It then follows that the spectral shift function computes spectral flow
with the caveat that there  are simple  correction terms arising from
the zero eigenspace of each endpoint of the path.

One motivation for this investigation arises from~\cite{CPRS2} where
the odd local index theorem in noncommutative geometry was deduced
from the analytic formula for spectral flow in~\cite{CP2} via an
intermediate formula in terms of a cyclic cocycle which was called
the `resolvent cocycle'. The results and viewpoint of this paper may have other
applications to noncommutative geometry besides spectral flow. Specifically we have
in mind situations where we would like to avoid `summability constraints' on $\HH_0$.

\section{Notation and Preliminary results} \label{S: preliminaries}
\subsection{Notation}
We write  $\hilb$ for a fixed separable complex Hilbert space, and
by $\clN$ we denote a semifinite von Neumann algebra acting on $\hilb$
\refsemifinvNa.
We use the usual notation $\HH \affl \clN$ for
operators~$\HH$ affiliated with~$\clN$ \refaffloper.
We let $\tau$ be a \nsf trace on $\clN$ \refnsftrace\
and by $L^1(\clN,\tau)$ we denote the set of $\tau$-trace class
operators affiliated with $\clN$ \reftautraceclassaffl.
Then we use the notation $\LpN{1}=L^1(\clN,\tau)\cap \clN$
for the (unitarily) invariant operator ideal \refinvoperideal\
of all bounded $\tau$-trace class operators
with norm $\clLnorm{\cdot} = \norm{\cdot} + \norm{\cdot}_1,$
where $\norm{\cdot}$ is the usual operator norm and $\norm{\cdot}_1 := \tau(\abs{\cdot})$
\reftautracenorm.
By $\clB(\mbR)$ we denote the $\sigma$-algebra of all Borel subsets of $\mbR.$
For a self-adjoint operator $T$ let $E^T_\Delta$ be
the spectral projection of $T$ corresponding to $\Delta \in \clB(\mbR),$
and let $E^T_\lambda$ be the spectral projection of $T$ corresponding to $(-\infty,\lambda].$
Write $C^\infty_\csupp(\mbR)$ (respectively, $C^\infty_\csupp(\Omega)$) for
the set of all compactly supported $C^\infty$-smooth functions on $\mbR$
(respectively, $\Omega \subseteq \mbR$),
and  $B(\mbR)$ (respectively, $B_\csupp(\mbR)$) for
set of all bounded Borel functions on $\mbR$
(respectively, compactly supported bounded Borel functions on $\mbR$).

If an operator $T$ affiliated with $\clN$ is $\tau$-measurable
\reftaumeasurable, then the $t$-th generalized
$s$-number $\mu_t(T)$ of the operator $T$ is defined by
$$
  \mu_t(T) = \inf \set{ \norm{TE} \colon E \ \text{is a projection in} \ \clN \ \text{with} \ \tau(1-E) \leq t}.
$$
If $E \in \clN$ is a $\tau$-finite projection (i.e. $\tau(E) < \infty$),
then
\begin{gather}\label{F: mu t(E) = chi [0,tau(E)]}
  \mu_t(E) = \chi_{[0,\tau(E))}(t), \ t \geq 0.
\end{gather}
An operator $T \in \clN$ is said to be $\tau$-compact iff $\liminfty {t} \mu_t(T) = 0.$
The set of all $\tau$-compact operators from $\clN$ forms an ideal in $\clN$ which is denoted $\clKN.$
The ideal $\clKN$ coincides with the norm closure of the ideal in $\clN$
generated by $\tau$\tire finite projections.

Let $R_z(\HH) = (z-\HH)^{-1}$ denote the resolvent of $\HH.$
If $\HH=\HH^* \affl \clN$ and $R_z(\HH) \in \clKN$ for some (and hence for all)
$z \in \mbC\setminus \mbR,$ then we say that $\HH$ has $\tau$-compact resolvent.

An operator $T \in \clN$ is called $\tau$-Fredholm \reftauFredholm\ if the projections $[\ker T],$ $[\ker T^*]$
are $\tau$-finite and there exists a $\tau$-finite projection $E \in \clN$ such that $\rng(1-E) \subseteq \rng(T).$
For a $\tau$-Fredholm operator $T$ one can define its $\tau$-index by $\tind(T) = \tau([\ker T]) - \tau([\ker T^*]).$
\begin{defn} \label{D: defn of Frechet der}
Let $\clN$ be a semifinite \vNa\ and let $\clE, \clE'$ be two invariant operator
ideals over $\clN.$ We denote self-adjoint part of $\clE$ by $\clE_{sa}.$
Let $\HH_0$ be a fixed self-adjoint operator affiliated with $\clN.$
A function $f \colon \HH_0 + \clE_{sa} \mapsto f(\HH_0) + \clE'_{sa}$
is called affinely $(\clE,\clE')$-Fr\'echet differentiable at $\HH \in \HH_0 + \clE_{sa},$
if there exists a (necessarily unique)
bounded operator $L \colon \clE_{sa} \to \clE'_{sa}$ such that
the following equality holds
$$
  f(\HH + \VV) - f(\HH) = L(\VV) + o(\norm{\VV}_\clE), \ \ \VV \in \clE_{sa}.
$$
If $\clE=\clE'$ then we write $L = D_{\clE} f(\HH),$ otherwise we write
$L = D_{\clE,\clE'}f(\HH).$
\end{defn}
In our case the ideals $\clE$ and $\clE'$ will be $\clN$ or $\LpN{1}.$

We constantly use some parameters for specific purposes.
The parameter~$\rr$ will always be an operator path parameter, i.\,e. the letter~$\rr$
is used when we consider paths of operators such as  $\HH_\rr = \HH_0 + \rr\VV.$
Very rarely we need another path parameter which we denote by~$s.$
We do not use $t$ as path parameter, since $t$ is used for other purposes
later in the paper.
The letter $\lambda$ is always used as a spectral parameter. If we need another spectral parameter
we will use~$\mu.$

We finish this subsection with the following elementary fact which we will use repeatedly.
\begin{lemma}\label{L: f = f1 - f2, where f1,f2 are nonnegative}
If $\Omega$ is an open interval in $\mbR$ and if $f \in C^k_\csupp(\Omega),$
then there exist functions $f_1, f_2 \in C^k_\csupp(\Omega)$
such that $f_1, f_2$ are non-negative, $f = f_1 - f_2$ and $\sqrt {f_1}, \sqrt {f_2} \in C^k_\csupp(\Omega).$
\end{lemma}
\begin{proof} Let $[a,b]$ be a closed interval, $[a,b] \subseteq \Omega$ and $\supp(f) \subseteq (a,b).$
Take a non-negative $C^\infty$-function $f_1 \geq f$ on $[a,b]$ which vanishes at $a$ and $b$ in such a way that
$\sqrt {f_1}$ is $C^\infty$-smooth at $a$ and $b,$ and take $f_2 = f_1 - f.$
\end{proof}
\subsection{Self-adjoint operators with $\tau$-compact resolvent}
The remaining subsections in this Section are technical and dry.
The reader may wish to browse this Section then move on
to Sections \ref{S: spectral shift measure} and \ref{S: spectral flow}
which contain the main results referring
back to the technicalities of this Section when needed.

In this current subsection we collect some facts about operators with compact resolvent
in a semifinite von Neumann algebra. We do not claim any great originality for these.
\begin{lemma}\label{L: D+V has comp resolvent} If $\HH=\HH^* \affl \clN$ has $\tau$-compact resolvent
   and $\VV=\VV^* \in \clN,$ then $\HH+\VV$ also has $\tau$-compact resolvent.
\end{lemma}
\begin{proof} This follows from the equality
$R_z(\HH+\VV) = \brs{R_z(\HH+\VV) \VV + 1}R_z(\HH),$ $\quad z \in \mbC \setminus \mbR.$
\end{proof}
\begin{lemma}\label{L: If D has comp rslv then E Delta if finite}
   If $\HH=\HH^* \affl \clN$ has $\tau$-compact resolvent, then for all compact sets $\Delta \subseteq \mbR$
   the spectral projection $E^\HH_\Delta$ is $\tau$-finite.
\end{lemma}
\begin{proof} If $\HH$ has $\tau$-compact resolvent then the operator
$(1+\HH^2)^{-1} = (\HH+i)^{-1}(\HH-i)^{-1}$ is $\tau$-compact.
Since for every finite interval $\Delta$ there exists a constant $c>0,$ not depending on $\HH$ such that
$E_\Delta^\HH \leq c (1+\HH^2)^{-1},$ the projection $E_\Delta^\HH$ is also $\tau$-compact, and hence
$\tau$-finite.
\end{proof}
\begin{cor}\label{C: if B has comp rslv then f(B) is trace class}
   If $\HH=\HH^* \affl \clN$ has $\tau$-compact resolvent, then for all $f \in B_\csupp(\mbR)$
   the operator $f(\HH)$ is $\tau$-trace class.
\end{cor}
\begin{proof}
  There exists a finite segment $\Delta \subseteq \mbR$
  such that $\abs{f} \leq \const \chi_\Delta,$ so that $\abs{f(\HH)} \leq \const E^\HH_\Delta.$
\end{proof}
\begin{lemma}\cite[Appendix B, Lemma 6]{CP98CJM} \label{CP, App B, Lemma 6}
If $\HH_0$ is an unbounded self-adjoint operator, $A$ is a bounded self-adjoint operator,
and $\HH = \HH_0 + A$ then
$$
  (1+\HH^2)^{-1} \leq f(\norm{A}) (1+\HH_0^2)^{-1},
$$
where $f(a) = 1 + \frac 12 a^2 + \frac 12 a\sqrt{a^2+4}.$
\end{lemma}
\begin{lemma}\label{L: lemma with def of B R}
   Let $\HH_0=\HH_0^* \affl \clN$ have $\tau$-compact resolvent,
   and let $B_R = \set{\VV = \VV^* \in \clN \colon \norm{\VV} \leq R}.$
Then for any compact subset $\Delta \subseteq \mbR$ the function
$$
  \VV \in B_R \mapsto E_\Delta^{\HH_0+\VV}
$$
is $\LpN{1}$-bounded.
\end{lemma}
\begin{proof} We have
\ $
   E^{\HH_0+\VV}_\Delta \leq c_0 (1+(\HH_0+\VV)^2)^{-1}
$ \
for some constant $c_0 = c_0(\Delta)>0$ and for every $\VV=\VV^* \in \clN.$
Now, by Lemma \ref{CP, App B, Lemma 6}
there exists a constant $c_1=c_1(R) > 0,$ such that for all $\VV \in B_R$
$$
  (1+(\HH_0+\VV)^2)^{-1} \leq c_1 (1+\HH_0^2)^{-1}.
$$
Hence, since $\HH_0$ has $\tau$-compact resolvent,
all projections $E^{\HH_0+\VV}_\Delta,$ $\VV \in B_R,$ are bounded from above
by a single $\tau$-compact operator $T=c_0 c_1 (1+\HH_0^2)^{-1}.$
This means, that for $t>0$
$$
  \mu_t(E^{\HH_0+\VV}_\Delta) \leq \mu_t(T).
$$
Further, by (\ref{F: mu t(E) = chi [0,tau(E)]})
 $\mu_t(E^{\HH_0+\VV}_\Delta) = \chi_{[0,\tau(E^{\HH_0+\VV}_\Delta))}(t)$
and there exists $t_0 > 0$ such that \mbox{$\mu_{t_0}(T) \leq 1.$} This implies that for all $\VV \in B_R$,
\ $
  \tau(E^{\HH_0+\VV}_\Delta) \leq t_0.
$ \
\end{proof}
\begin{cor}\label{C: trace f(D0+V) is bounded}
  If $\HH_0=\HH_0^* \affl \clN$ has $\tau$-compact resolvent,
  then for any function $f \in B_c(\mbR)$ the function
  $\VV \in B_R \mapsto \clLnorm{f(\HH_0+\VV)}$ is bounded.
\end{cor}
\begin{cor}\label{C: trace Br Delta is bounded}
  Let $\HH_0=\HH_0^* \affl \clN$ have $\tau$-compact resolvent, $\rr=(\rr_1,\ldots,\rr_m) \in [0,1]^m,$
  $\VV_1,\ldots,\VV_m \in \clN_{sa}$ and set $\HH_\rr = \HH_0+\rr_1\VV_1+\ldots+\rr_m\VV_m.$
  Then \\ (i) for any compact subset $\Delta \subseteq \mbR$ the function
  $\rr \in [0,1]^m \mapsto \norm{E^{\HH_\rr}_\Delta}_1$ is bounded;
  \\ (ii) for any function $f \in B_c(\mbR)$ the function
  $\rr \in [0,1] \mapsto \norm{f(\HH_\rr)}_1$ is bounded.
\end{cor}
An elementary proof of the following lemma can also be found in \cite{CP98CJM}.
\begin{lemma}\label{L: exp(B+rV) - exp(B) is norm continuous}
  If $\HH_0=\HH_0^* \affl \clN$ and if $\VV=\VV^* \in \clN,$
  then for any $t \in \mbR$, \
  $e^{it(\HH_0+\VV)}$ converges in $\norm{\cdot}$-norm to $e^{it\HH_0}$ when $\norm{\VV} \to 0.$
\end{lemma}
\begin{proof} Follows directly from Duhamel's formula \ \
$$
  e^{it(\HH_0+\VV)} - e^{it\HH_0} = \int_0^t e^{i(t-u)(\HH_0+\VV)} i\VV e^{iu\HH_0}\,du.
$$
\end{proof}

If $(S,\nu)$ is a finite measure space then we say that a function
$f \colon S \to \clBH$ is measurable, if for any $\eta \in \hilb$
the functions $f(\cdot)\eta, f(\cdot)^*\eta \colon S \to \hilb$ are Bochner measurable.
We define the integral of a measurable function by the formula
$$
  \int_S f(\sigma)\,d\nu(\sigma)\eta = \int_S f(\sigma)\eta\,d\nu(\sigma), \ \ \eta \in \hilb,
$$
and call it an $so^*$-integral. If $f \colon S \to \LpN{1}$ is an $\clL^1$\tire bounded function then
it is measurable if and only if there exists a sequence of simple (finitely-valued) functions $f_n \colon S \to \LpN{1}$
such that $f_n(\sigma) \to f(\sigma)$ in the $so^*$-topology \cite[Section 3]{ACDS}.
\begin{lemma}\label{L: lemma 3.10 ACDS} \cite[Lemma 3.10]{ACDS}
  If a measurable function $f \colon S \to \LpN{1}$ is uniformly $\LpN{1}$-bounded
  then $\int_S f(\sigma)\,d\nu(\sigma) \in \LpN{1}$ and
  $$
    \taubrs{\int_S f(\sigma)\,d\nu(\sigma)} = \int_S \taubrs{f(\sigma)}\,d\nu(\sigma).
  $$
\end{lemma}
\subsection{Difference quotients and double operator integrals}
Originally, multiple operator integrals of the form
\begin{multline}\label{D: formal MOI}
  T^{\HH_0,\ldots,\HH_n}_{\phi}(\VV_1,\ldots,\VV_n)
    \\ = \int_{\mbR^{n+1}} \phi(\lambda_0,\ldots,\lambda_n)
    \,dE^{\HH_0}_{\lambda_0} \VV_1 d E^{\HH_1}_{\lambda_1} \VV_2 \ldots \VV_n dE^{\HH_n}_{\lambda_n},
\end{multline}
where $\HH_0,\ldots,\HH_n$ are self-adjoint operators and $\VV_1,\ldots,\VV_n$ are bounded operators,
were defined as repeated operator
integrals~\cite{DalKr} or as spectral integrals~\cite{BS66DOI,BS67DOI,BS73DOI}
in the case of double operator integrals.
It was noted in~\cite{ACDS} (see also \cite{Pel2}) that one can give another definition
of multiple operator integrals (Definition \ref{D: def of MOI}). The idea of the new definition
is that for functions $\phi(\lambda_0,\ldots,\lambda_n)$ of the form
$\phi(\lambda_0,\ldots,\lambda_n) = \alpha_0(\lambda_0)\alpha_1(\lambda_1)\ldots\alpha_n(\lambda_n)$
the expression (\ref{D: formal MOI}) can be interpreted as $\alpha_0(\HH_0)\VV_1\alpha_1(\HH_1)\VV_2\ldots\VV_n\alpha_n(\HH_n).$
Hence, for functions $\phi$ of the form (\ref{F: BS representation})
one can try to define the multiple operator integral by the formula (\ref{F: Def of MOI}).
Having proved the correctness of this definition (Proposition \ref{P: Def of MOI is correct}),
we will be able to work with the usual operator integrals of the form (\ref{F: Def of MOI}), provided
that one has a representation of the form of equation (\ref{F: BS representation}) for the function $\phi.$
Actually, one has considerable freedom of choice of this representation, and one should try
to find that representation of $\phi$ which is most suitable for the purpose at hand.

In our case the function $\phi(\lambda_0,\lambda_1)$ is the difference quotient
\begin{gather}\label{F: divided difference}
  f^{[1]}(\lambda_0,\lambda_1) := \frac{f(\lambda_0)-f(\lambda_1)}{\lambda_0-\lambda_1}
\end{gather}
of some function $f.$
We denote by $\CPlus{n}$ the set of functions
$f\in C^n(\mbR),$ such that the $j$-th derivative $f^{(j)}, \ j=1,\ldots,n,$
belongs to the space $\clF^{-1}(L^1(\mbR)),$ where $\fourier f$ is the Fourier transform of $f.$
We emphasize that the Fourier transform of a function $f \in \CPlus{n}$ need not to belong to $L^1(\mbR)$
(see Lemma \ref{L: C1plus is subset of L1}(ii)).
Let
\begin{gather*}
  \Pi := \set{(s_0,s_1) \in \mbR^{2}\colon \abs{s_1}\leq\abs{s_0}, \sign(s_0) = \sign(s_1)},
\end{gather*}
and
\begin{gather} \label{F: def of nu f}
  d\nu_f(s_0,s_1) := \sgn(s_0) \frac {i}{\sqrt{2\pi}} \fourier f(s_0)\,ds_0\,ds_1.
\end{gather}
If $f\in \CPlus{1},$ then it is not difficult to see that $\brs{\Pi, \nu_f}$
is a finite measure space~\cite[Lemma 2.1]{ACDS}, so that $\brs{\Pi, \nu_f}$
can be used for the construction of the double operator integral.
The following Birman-Solomyak or BS-representation of $\psif{f}(\lambda_0,\lambda_1)$~\cite[Lemma 2.2]{ACDS}
  \begin{gather}\label{F: BS repr-n from ACDS}
    \psif{f}(\lambda_0,\lambda_1)
      = \int_\Pi \alpha(\lambda_0,\sigma)\beta(\lambda_1,\sigma)\,d\nu_f(\sigma),
  \end{gather}
where $\sigma = (s_0,s_1),$ $\alpha(\lambda_0,\sigma) = e^{i(s_0-s_1)\lambda_0}$ and
$\beta(\lambda_1,\sigma) = e^{is_1\lambda_1},$ is not suitable for our present purposes.
Lemma \ref{L: BS representation for f(x) - f(y)} provides a modification
of this BS-representation for $f^{[1]}$ with which we will constantly work.

We include the proof of the following fact for completeness.
\begin{lemma}\label{L: C1plus is subset of L1}
  (i) If $f \in C^1_\csupp(\mbR),$ then $\hat f \in L^1(\mbR).$
  \\ (ii) the function $\phi(x) = \frac x{\sqrt{1+x^2}}$ belongs to $\CPlus 2.$
\end{lemma}
\begin{proof} (i) Following the proof of~\cite[Corollary 3.2.33]{BR1}, we have
  \begin{multline*}
    \int_\mbR |\hat f|(\xi)\,d\xi = \int_\mbR |\xi+i|^{-1} |\xi+i| |\hat f|(\xi)\,d\xi
    \\ \leq \brs{\int_\mbR |\xi+i|^{-2}\,d\xi}^{\frac 12} \brs{\int_\mbR|\xi+i|^2 |\hat f|^2(\xi)\,d\xi}^{\frac 12}
    \\ = \const \brs{\int_\mbR \abs{\frac {df(x)}{dx} + f(x)}^2 \,dx}^{\frac 12} < \infty.
  \end{multline*}

  (ii) The proof  is similar to that of  (i).
\end{proof}
Applying part (i) of this lemma to the first $n$ derivatives
of a function $f$ from $C^{n+1}_\csupp(\mbR),$ we obtain
\begin{cor}\label{C: C(n+1)(comp supp) is a subset of C+(n)} $C^{n+1}_\csupp(\mbR) \subseteq \CPlus n,$ $n = 1,2,\ldots$
\end{cor}
The following lemma provides a BS-representation for $f^{[1]},$ $0\leq f \in C^2_\csupp(\mbR),$
which will be used throughout this paper.
\begin{lemma}\label{L: BS representation for f(x) - f(y)} Let
$f \in C^2_\csupp(\mbR)$ be a non-negative function such that $g:=\sqrt f \in C^2_\csupp(\mbR).$
If $\Omega \supseteq \supp(f),$ then
  \begin{gather*}
    \psif{f}(\lambda_0,\lambda_1)
      = \int_\Pi \brs { \alpha_1(\lambda_0,\sigma)\beta_1(\lambda_1,\sigma) +
               \alpha_2(\lambda_0,\sigma)\beta_2(\lambda_1,\sigma) } \,d\nu_g(\sigma),
  \end{gather*}
   where $\sigma = (s_0,s_1)$ and
   \begin{gather}\label{F: alpha1,alpha2,beta1,beta2}
     \alpha_1(\lambda_0,\sigma) = e^{i(s_0-s_1)\lambda_0} g(\lambda_0),   \qquad
     \beta_1(\lambda_1,\sigma) = e^{is_1\lambda_1},   \\ \notag
     \alpha_2(\lambda_0,\sigma) = e^{i(s_0-s_1)\lambda_0},     \qquad
     \beta_2(\lambda_1,\sigma) = e^{is_1\lambda_1}g(\lambda_1),
   \end{gather}
   so that
   $\alpha_1(\cdot,\sigma), \,\, \beta_2(\cdot,\sigma) \in C_\csupp^2(\Omega)$ for all
   $\sigma \in \Pi, and $ $\abs{\alpha_1(\cdot)}, \abs{\beta_2(\cdot)} \leq \norm{g}_\infty,$
   while $\alpha_2(\cdot,\sigma), \,\, \beta_1(\cdot,\sigma) \in C^\infty(\mbR)$ for all
   $\sigma \in \Pi, and $ $\abs{\alpha_2(\cdot)}, \abs{\beta_1(\cdot)} \leq 1.$
\end{lemma}
\begin{proof}
The assumption $g \in C^2_\csupp(\mbR)$ implies that $g \in \CPlus 1$ (see
Corollary \ref{C: C(n+1)(comp supp) is a subset of C+(n)}).
Now,
  \begin{multline*}
    \psif{f}(\lambda_0,\lambda_1)
     = \frac {g^2(\lambda_0) - g^2(\lambda_1)}{\lambda_0-\lambda_1}
     \\ = \frac {g(\lambda_0) - g(\lambda_1)}{\lambda_0-\lambda_1} \brs{g(\lambda_0) + g(\lambda_1)}
     = \psif{g}(\lambda_0,\lambda_1) \brs{g(\lambda_0) + g(\lambda_1)}.
  \end{multline*}
Hence, using (\ref{F: BS repr-n from ACDS}), we have
  \begin{gather*}
    \psif{f}(\lambda_0,\lambda_1)
      = \int_\Pi \Big( \alpha(\lambda_0,\sigma)g(\lambda_0)\beta(\lambda_1,\sigma)
      +  \alpha(\lambda_0,\sigma)g(\lambda_1)\beta(\lambda_1,\sigma)\Big) \,d\nu_g(\sigma).
  \end{gather*}
If we set $\alpha_1(\lambda,\sigma) = \alpha(\lambda,\sigma)g(\lambda),$ $\beta_1(\lambda,\sigma) = \beta(\lambda,\sigma),$
$\alpha_2(\lambda,\sigma) = \alpha(\lambda,\sigma)$ and
$\beta_2(\lambda,\sigma) = g(\lambda)\beta(\lambda,\sigma),$
then we see that all the conditions of the Lemma are fulfilled.
\end{proof}
The next step is to recall the definition of the multiple operator integral as it was given in~\cite{ACDS}.
Let $\phi \in B(\mbR^{n+1})$ be a bounded Borel function on $\mbR^{n+1}$
which admits a representation of the form
\begin{gather}\label{F: BS representation}
  \phi(\lambda_0,\lambda_1,\ldots,\lambda_n) = \int_{S} \alpha_0(\lambda_0,\sigma)\ldots\alpha_n(\lambda_n,\sigma)\,d\nu(\sigma),
\end{gather}
where $(S,\nu)$ is a finite measure space and $\alpha_0,\ldots,\alpha_n$ are bounded Borel functions on $\mbR\times S.$
\begin{defn}\label{D: def of MOI}
   For arbitrary self-adjoint operators $\HH_0,\ldots,\HH_n$ on
   the separable Hilbert space $\hilb,$ bounded operators $\VV_1,\ldots,\VV_{n}$ on $\hilb$
   and any function $\phi \in B(\mbR^{n+1})$ which admits a representation given by (\ref{F: BS representation}),
   the multiple operator integral $T^{\HH_0,\ldots,\HH_n}_\phi(\VV_1,\ldots,\VV_{n})$ is defined as
\begin{gather}\label{F: Def of MOI}
  T^{\HH_0,\ldots,\HH_n}_\phi(\VV_1,\ldots,\VV_{n}) := \int_{S} \alpha_0(\HH_0,\sigma)\VV_1\ldots \VV_{n}\alpha_n(\HH_n,\sigma)\,d\nu(\sigma),
\end{gather}
where the integral is taken in the $so^*$-topology.
\end{defn}
\begin{prop}\label{P: Def of MOI is correct}~\cite[Lemma 4.3]{ACDS}
The multiple operator integral in Definition \ref{D: def of MOI} is well-defined
in the sense that it does not depend on the representation (\ref{F: BS representation}) of $\phi.$
\end{prop}
The following lemma is a corollary of Lemma \ref{L: BS representation for f(x) - f(y)}
and the definition of the multiple operator integral.
\begin{lemma}\label{L: BS representation for DOI} If $\HH_0=\HH_0^*\affl \clN,$
if $\HH_1 = \HH_0+\VV,$ $\VV = \VV^* \in \clN$ and if $f \in C^2_\csupp(\mbR)$ is a non-negative function,
such that $g := \sqrt f \in C^2_\csupp(\mbR),$ then
$$
  T^{\HH_1,\HH_0}_{f^{[1]}}(\VV) = \int_\Pi \brs { \alpha_1(\HH_1,\sigma)\VV\beta_1(\HH_0,\sigma) +
                 \alpha_2(\HH_1,\sigma)\VV\beta_2(\HH_0,\sigma) }  \,d\nu_g(\sigma),
$$
where $\alpha_1,\beta_1,\alpha_2,\beta_2$ are given by (\ref{F: alpha1,alpha2,beta1,beta2}).
\end{lemma}
We need the following weaker version of~\cite[Theorem 5.3]{ACDS}. See also~\cite{BS73DOI}.
\begin{prop} \label{P: DalKreinII}
  \cite[Theorem 5.3]{ACDS} 
  Let $\clN$ be a \vNa. Suppose that $\HH_0=\HH_0^*$ is affiliated with $\clN,$ that
  $\VV \in \clN$ is self-adjoint and set $\HH_1=\HH_0+\VV.$
\\ (i) If $f \in \CPlus 1,$ then
    \begin{gather*}
      f(\HH_1) - f(\HH_0) = T^{\HH_1,\HH_0}_{\psif{f}}(\VV).
    \end{gather*}
\\ (ii)
  If $f \in \CPlus 2,$ then the function $f \colon \HH_0 + \clN_{sa} \mapsto f(\HH_0) + \clN_{sa}$
  is affinely $(\clN,\clN)$-Fr\'echet differentiable, the equality $D_\clN f(\HH) = T^{\HH,\HH}_{f^{[1]}}$
  holds and $D_\clN f(\HH)$ is $\norm{\cdot}$-continuous, where $D_\clN $ is to be understood in the sense
  of Definition \ref{D: defn of Frechet der}.
\end{prop}
%
%
\subsection{Some continuity and differentiability properties of operator functions}
We are going to consider spectral flow along `continuous' paths of
unbounded  Fredholm operators. We will make precise
what we mean by continuity in this setting later.
However our formulae require more than just continuity.
They require us to be able to take derivatives with
the respect to the path parameter. For this to be feasible we need the full force of
the double operator integral formalism.
We present the results we will need as a sequence of lemmas.
\begin{lemma} \label{L: norm of integral < integral of norm}
Let $(S,\nu)$ be a finite measure space and let $f \colon S \to \LpN{1}$
be a $\LpN{1}$-bounded $so^*$-measurable function. Then
$$
  \norm{\int f(\sigma)\,d\nu(\sigma)}_{\clL^1}
    \leq \int \norm{f(\sigma)}_{\clL^1} \,d\abs{\nu}(\sigma).
$$
\end{lemma}
\begin{proof} By definition for any $\eta \in \hilb$ the function
$\sigma \mapsto f(\sigma)\eta$ is Bochner measurable. Hence,
the function $\sigma \mapsto \norm{f(\sigma)} = \sup_{\eta\in \hilb\colon \norm{\eta} \leq 1}\norm{f(\sigma)\eta}$
is also measurable. Similarly, since the function $\sigma \mapsto \taubrs{f(\sigma)B}$ is measurable,
the function $\sigma \mapsto \norm{f(\sigma)}_1 = \sup_{B \in \clN\colon \norm{B} \leq 1} \abs{\taubrs{f(\sigma)B}}$
is also measurable. Hence, the right hand side of the last equality is well defined.

For $\eta \in \hilb$ with $\norm{\eta} \leq 1$ we have
by definition of $so^*$-integral (see~\cite[(2)]{ACDS})
\begin{multline}\label{F: first inequality}
  \norm{\int f(\sigma)\,d\nu(\sigma) \eta}
   = \norm{\int f(\sigma)\eta \,d\nu(\sigma)}
     \\ \leq \int \norm{f(\sigma)\eta} \,d\abs{\nu}(\sigma)
     \leq \int \norm{f(\sigma)} \,d\abs{\nu}(\sigma).
\end{multline}
Hence, the inequality is true for $\norm{\cdot}$-norm.
Since $\norm{A}_1 = \sup_{B\in \clN\colon \norm{B}\leq 1} \abs{\tau(AB)},$
we have
\begin{multline}\label{F: second inequality}
  \norm{\int f(\sigma)\,d\nu(\sigma)}_1
  = \sup_{B\in \clN\colon \norm{B}\leq 1}  \abs{\taubrs{\int f(\sigma) \,d\nu(\sigma) B}}
  \\ = \sup_{B\in \clN\colon \norm{B}\leq 1}  \abs{\int \taubrs{f(\sigma)B} \,d\nu(\sigma)}
  \\ \leq \sup_{B\in \clN\colon \norm{B}\leq 1}  \int \abs{\taubrs{f(\sigma)B}} \,d\nu(\sigma)
    \leq \int \norm{f(\sigma)}_1 \,d\abs{\nu}(\sigma),
\end{multline}
where the second equality follows from the definition of $so^*$-integral and
Lemma \ref{L: lemma 3.10 ACDS}.
Combining (\ref{F: first inequality}) and (\ref{F: second inequality}) completes the proof.
\end{proof}
%
%
In the sequel we will constantly need to take functions of a path of operators. We thus need the
following continuity result. For the definition of $B_R$ see Lemma \ref{L: lemma with def of B R}.
\begin{prop}\label{P: f(B+rV) is L1 continuous}
   If $\HH_0=\HH_0^* \affl \clN$ has $\tau$-compact resolvent
   and if $f \in C^2_\csupp(\mbR)$ then the
   operator-valued function $A \colon \VV \in B_R \mapsto f(\HH_0+\VV)$ takes values in $\LpN{1}$
   and is $\LpN{1}$-continuous.
\end{prop}
  \begin{proof} 
  That $A(\cdot)$ takes values in $\LpN{1}$ follows from Lemma \ref{L: D+V has comp resolvent} and
  Corollary \ref{C: if B has comp rslv then f(B) is trace class}.
  By Lemma \ref{L: f = f1 - f2, where f1,f2 are nonnegative}
  it is enough to prove continuity for a non-negative function $f$ with $g = \sqrt f \in C^2_\csupp(\mbR).$
  By Proposition \ref{P: DalKreinII}(i) and Lemma \ref{L: BS representation for DOI}  we have
  \begin{multline*}
    f(\HH_0+\VV) - f(\HH_0) = T_{f^{[1]}}^{\HH_0+\VV,\HH_0}(\VV)
    \\ = \int_\Pi \brs { \alpha_1(\HH_0+\VV,\sigma)\VV\beta_1(\HH_0,\sigma) +
                   \alpha_2(\HH_0+\VV,\sigma)\VV\beta_2(\HH_0,\sigma) } \,d\nu_g(\sigma).
  \end{multline*}
Hence, by Lemma \ref{L: norm of integral < integral of norm}, we have
  \begin{multline*}
    \clLnorm{f(\HH_0+\VV) - f(\HH_0)}
    \leq \int_\Pi \big[\clLnorm{\alpha_1(\HH_0+\VV,\sigma)} \norm{\VV}\norm{\beta_1(\HH_0,\sigma)}
              \\ + \norm{\alpha_2(\HH_0+\VV,\sigma)}\norm{\VV}\clLnorm{\beta_2(\HH_0,\sigma)} \big] \,d\abs{\nu_g}(\sigma)
    \\ \leq \int_\Pi \brs { \clLnorm{g(\HH_0+\VV)} \norm{\VV} +
                   \norm{\VV}\clLnorm{g(\HH_0)} } \,d\abs{\nu_g}(\sigma)
    \\ \leq \abs{\nu_g}(\Pi) \norm{\VV} (\clLnorm{g(\HH_0+\VV)} + \clLnorm{g(\HH_0)}).
  \end{multline*}
  Now, Corollary \ref{C: trace f(D0+V) is bounded} applied to $g$ completes the proof.
  \end{proof}
\begin{cor}\label{C: f(B+rV) is L1 continuous}
   If $\HH_0=\HH_0^* \affl \clN$ has $\tau$-compact resolvent,
    $\rr=(\rr_1,\ldots,\rr_m) \in [a,b]^m,$ if
   $\VV_1,\ldots,\VV_m \in \clN_{sa}$
   and if $\HH_\rr = \HH_0 + \rr\VV = \HH_0+\rr_1\VV_1+\ldots+\rr_m\VV_m,$
   then for any function $f \in C^2_\csupp(\mbR)$ the
   operator-valued function $A \colon \rr \in [a,b]^m \mapsto f(\HH_0+\rr\VV)$ takes values in $\LpN{1}$
   and is $\LpN{1}$-continuous.
\end{cor}
%
%
Next we prove the main lemmas of this Section. There are several matters
to establish. First we want to be able to differentiate, with respect to
the path parameter, certain functions of paths of operators. Then we need
to determine formulae for the derivatives and the
continuity properties of the derivatives with respect to the path parameter.
\begin{lemma}\label{L: DOI is continuous wrt D1,D2}
  If $\HH_1$ and $\HH_2$ are two self-adjoint operators with $\tau$-compact
  resolvent affiliated with semifinite \vNa\ $\clN,$ if $X \in \clN_{sa}$ and if $f \in C^3_\csupp(\mbR)$ then
  $T^{\HH_1,\HH_2}_{f^{[1]}}(X)$ depends $\clL^1$-continuously on $\norm{\cdot}$
  perturbations of $\HH_1$ and $\HH_2.$
\end{lemma}
\begin{proof}
As usual, we can assume that $f$ is non-negative and its square root $g = \sqrt f$ is $C^3$-smooth.

Let $Y_1,Y_2 \in \clN_{sa}.$ Then by Lemma \ref{L: BS representation for DOI}
  \begin{multline*}
    T_{f^{[1]}}^{\HH_1+Y_1,\HH_2+Y_2}(X) - T_{f^{[1]}}^{\HH_1,\HH_2}(X)
              \\ = \int_\Pi \big[
                \alpha_1(\HH_1+Y_1,\sigma) X \beta_1(\HH_2+Y_2,\sigma) + \alpha_2(\HH_1+Y_1,\sigma) X \beta_2(\HH_2+Y_2,\sigma)
                \\ - \alpha_1(\HH_1,\sigma) X \beta_1(\HH_2,\sigma) - \alpha_2(\HH_1,\sigma) X \beta_2(\HH_2,\sigma)
              \big]\,d\nu_{g}(\sigma)
              \\ = \int_\Pi \Big (
                \big[\alpha_1(\HH_1+Y_1,\sigma) - \alpha_1(\HH_1,\sigma)\big] X \beta_1(\HH_2+Y_2,\sigma)
                \\ + \alpha_1(\HH_1,\sigma) X \big[ \beta_1(\HH_2+Y_2,\sigma) - \beta_1(\HH_2,\sigma) \big]
                \\ + \big[ \alpha_2(\HH_1+Y_1,\sigma) - \alpha_2(\HH_1,\sigma)\big] X \beta_2(\HH_2+Y_2,\sigma)
                \\ + \alpha_2(\HH_1,\sigma) X \big[ \beta_2(\HH_2+Y_2,\sigma) - \beta_2(\HH_2,\sigma) \big]
              \Big)\,d\nu_{g}(\sigma).
  \end{multline*}
For every fixed \mbox{$\sigma \in \Pi$}
by Lemma \ref{L: exp(B+rV) - exp(B) is norm continuous} the norms
$\norm{\beta_1(\HH_2+Y_2,\sigma) - \beta_1(\HH_2,\sigma)}$ and
$\norm{\alpha_2(\HH_1+Y_1,\sigma) - \alpha_2(\HH_1,\sigma)}$ converge to zero
when $\norm{Y_1},$ \, $\norm{Y_2} \to 0,$
and by Corollary \ref{C: trace f(D0+V) is bounded}
the $\clL^1$-norms of $\alpha_1(\HH_1,\sigma)$ and $\beta_2(\HH_2+Y_2,\sigma)$ are bounded
when $\norm{Y_1}, \norm{Y_2} \to 0.$
Hence, for every fixed $\sigma \in \Pi$ the $\clL^1$-norms of the second
and third summands in the last integral converge to zero
when $\norm{Y_1}, \norm{Y_2} \to 0.$

Now we are going to show that the same is true for the first and fourth summands.
It is enough to prove that for every fixed $\sigma \in \Pi,$ for example,
$\norm{\alpha_1(\HH_1+Y_1,\sigma) - \alpha_1(\HH_1,\sigma)}_{\clL^1}$ tends to zero. We have
\begin{multline*}
  \norm{\alpha_1(\HH_1+Y_1,\sigma) - \alpha_1(\HH_1,\sigma)}_{\clL^1}
  \\ = \norm{e^{i(s_0-s_1)(\HH_1+Y_1)} g(\HH_1+Y_1) - e^{i(s_0-s_1)\HH_1} g(\HH_1)}_{\clL^1}
  \\ \leq \norm{\brs{e^{i(s_0-s_1)(\HH_1+Y_1)} - e^{i(s_0-s_1)\HH_1}}g(\HH_1+Y_1)}_{\clL^1}
  \\  + \norm{e^{i(s_0-s_1)\HH_1} \brs{ g(\HH_1+Y_1) - g(\HH_1)}}_{\clL^1}
  \\ \leq \norm{\brs{e^{i(s_0-s_1)(\HH_1+Y_1)} - e^{i(s_0-s_1)\HH_1}}} \norm{g(\HH_1+Y_1)}_{\clL^1}
  \\  + 
        \norm{ g(\HH_1+Y_1) - g(\HH_1)}_{\clL^1}.
\end{multline*}
It follows from Lemma \ref{L: exp(B+rV) - exp(B) is norm continuous} that the first summand converges to zero
when $s_0,s_1$ are fixed and $\norm{Y_1} \to 0,$
and it follows from Proposition \ref{P: f(B+rV) is L1 continuous} that the second summand also converges to zero.

Since by Corollary \ref{C: trace f(D0+V) is bounded}
the trace norm of the expression under the last integral is uniformly $\LpN{1}$-bounded
with respect to $\sigma \in \Pi,$ it follows from Lemma \ref{L: norm of integral < integral of norm} that
$$
  \norm{T_{f^{[1]}}^{\HH_1+Y_1, \HH_2+Y_2}(X) - T_{f^{[1]}}^{\HH_1, \HH_2}(X)}_{\clL^1} \to 0,
$$
when $\norm{Y_1}, \norm{Y_2} \to 0.$
\end{proof}
The following theorem is a version of a well-known Daletskii-\KreinSG\ formula \cite{DalKr}.
We would like to give a heuristic argument explaining the formula (\ref{F: DN = TDD}).
In the resolvent expansion series
$$
  \frac 1{z-H-V} = \frac 1{z-H} + \frac 1{z-H}V\frac 1{z-H} + \frac 1{z-H}V\frac 1{z-H}V\frac 1{z-H} + \ldots
$$
the second summand is a double operator integral $T^{H,H}_\phi(V)$ with
$$
  \phi(\lambda,\mu) = \frac 1{z-\lambda}\cdot \frac 1{z-\mu} = \frac 1{\lambda-\mu}\brs{\frac 1{z-\lambda} - \frac 1{z-\mu}}.
$$
If $H$ and $H+V$ are bounded operators and $f$ is a function analytic in a neighbourhood
of the union of the spectra of $H$ and $H+V$ then the Cauchy integral implies
$$
  f(H+V) = f(H) + T^{H,H}_{f^{[1]}}(V) + \text{terms of second order}.
$$
\begin{thm}\label{T: D L1 f(D) = T D,D}
  If the \vNa\ $\clN$ is semifinite, $\HH_0=\HH_0^* \affl \clN$
  has $\tau$-compact resolvent
  and $f \in C^3_{\csupp}(\mbR)$  then the function
  $f \colon \HH \in \HH_0 + \clN_{sa} \mapsto f(\HH) \in \clN_{sa}$
  takes values in $\LpN{1}_{sa}.$ Moreover, it is affinely $(\clN,\clL^{1})$-Fr\'echet differentiable,
  the equality
  \begin{gather} \label{F: DN = TDD}
    D_{\clN,\clL^1} f(\HH) = T^{\HH,\HH}_{f^{[1]}}
  \end{gather}
  holds, and $D_{\clN,\clL^1} f(\HH)$ is $(\clN,\clL^1)$-continuous, so that
  \begin{gather}\label{F: Newton-Leibnitz}
    f(\HH_b) - f(\HH_a) = \int_a^b T^{\HH_r,\HH_r}_{f^{[1]}} (\VV)\,dr,
  \end{gather}
  where $\VV \in \clN_{sa},$ $\HH_r = \HH_0 + \rr\VV$
  and the integral converges in $\LpN1$-norm.
\end{thm}
\begin{proof} We have by Proposition \ref{P: DalKreinII}(i) and Lemma \ref{L: BS representation for DOI}
  \begin{multline*}
    f(\HH_1) - f(\HH_0) =  T^{\HH_1,\HH_0}_{f^{[1]}}(\VV)
    \\ = \int_\Pi \brs { \alpha_1(\HH_1,\sigma)\VV\beta_1(\HH_0,\sigma) +
                   \alpha_2(\HH_1,\sigma)\VV\beta_2(\HH_0,\sigma) }  \,d\nu_g(\sigma),
    \\ = \int_\Pi \brs { \alpha_1(\HH_0,\sigma)\VV\beta_1(\HH_0,\sigma) +
                   \alpha_2(\HH_0,\sigma)\VV\beta_2(\HH_0,\sigma) }  \,d\nu_g(\sigma)
       \\ + \int_\Pi [\alpha_1(\HH_1,\sigma) - \alpha_1(\HH_0,\sigma)]\VV\beta_1(\HH_0,\sigma) \,d\nu_g(\sigma)
       \\ + \int_\Pi [\alpha_2(\HH_1,\sigma) - \alpha_2(\HH_0,\sigma)]\VV\beta_2(\HH_0,\sigma) \,d\nu_g(\sigma)
    \\ = T^{\HH_0,\HH_0}_{f^{[1]}}(\VV) + (II) + (III).
  \end{multline*}

  Since $\alpha_2$ is just an exponent and since $g \in \CPlus{2}$ that $\norm{(III)}_{\clL^1} = O(\norm{\VV}^2)$ can be shown
  by Duhamel's formula. The argument is
as in the proof of~\cite[Theorem 5.5]{ACDS}.
  So, it is left to show that $\norm{(II)}_{\clL^1}$ is $o(\norm{\VV}).$
  By Lemma \ref{L: norm of integral < integral of norm} we have
  \begin{gather*}
    \norm{(II)}_{\clL^1}
    = \norm{\int_\Pi [\alpha_1(\HH_1,\sigma) - \alpha_1(\HH_0,\sigma)]\VV\beta_1(\HH_0,\sigma) \,d\nu_g(\sigma)}_{\clL^1}
    \\ \leq \int_\Pi \norm{\alpha_1(\HH_1,\sigma) - \alpha_1(\HH_0,\sigma)}_{\clL^1} \norm{\VV} \norm{\beta_1(\HH_0,\sigma)} \,d\nu_g(\sigma)
    \\ = \norm{\VV} \int_\Pi \norm{\alpha_1(\HH_1,\sigma) - \alpha_1(\HH_0,\sigma)}_{\clL^1} \,d\nu_g(\sigma).
  \end{gather*}
Now, it follows from $\alpha_1(\cdot,\sigma) \in C_\csupp^2(\mbR)$ (see (\ref{F: alpha1,alpha2,beta1,beta2}))
and Proposition \ref{P: f(B+rV) is L1 continuous} that
$\norm{\alpha_1(\HH_1,\sigma) - \alpha_1(\HH_0,\sigma)}_{\clL^1} \to 0, \ \sigma \in \Pi,$ so that by the Lebesgue
dominated convergence theorem we conclude
that the last integral converges to 0, and hence $\norm{(II)}_{\clL^1} = o(\norm{\VV}).$

Finally, that $D_{\clN,\clL^1} f(\HH)$ is $(\clN,\clL^1)$-continuous follows from Lemma \ref{L: DOI is continuous wrt D1,D2}.
\end{proof}
\subsection{A class $\clFNab$ of $\tau$-Fredholm operators}
Our technique for handling spectral flow of paths of unbounded operators is to map
them into the space of bounded operators using a particular function.  We thus need
to discuss some continuity properties of paths of bounded $\tau$-Fredholm operators,
analogous to those we described in the unbounded case.

Let $a < b$ be two real numbers.
Let $\clFNab$ be the set of bounded self-adjoint $\tau$-Fredholm operators $F \in \clN$
such that $(F-a)(F-b) \in \clKN.$
For $F_0 \in \clFNab$ let $\clA_{F_0} = F_0 + \clKN_{sa}$ be the affine
space of $\tau$-compact self-adjoint perturbations of $F_0.$
\begin{lemma}\label{L: If F in clFNpm then F + K in clFNpm}
If $F_0 \in \clFNab$ then
$$\clA_{F_0} \subseteq \clFNab.$$
\end{lemma}
\begin{proof}
If $K \in \clKN_{sa}$ then $(F_0+K-a)(F_0+K-b) = (F_0-a)(F_0-b) + (F_0-a)K + K(F_0+K-b) \in \clKN.$
\end{proof}
\begin{lemma}\label{L: if F in clFNpm then h(F) is trace class}
  If $F \in \clFNab$ and $h \in B_\csupp(a,b)$ then $h(F) \in \LpN{1}.$
\end{lemma}
\begin{proof} The proof is similar to the proof of Lemma \ref{L: If D has comp rslv then E Delta if finite}.
For any compact subset $\Delta$ of $(a,b)$
there exists a constant $c_0>0$ such that
$\chi_\Delta(x) \leq c_0 \chi_{[a,b]}(x) (b-x)(x-a),$
so that
  \begin{gather}\label{F: E Delta < c(1-F2)...}
    E^F_\Delta \leq c_0 E^F_{[a,b]}(b-F)(F-a).
  \end{gather}
Since $(b-F)(F-a) \in \clKN,$ it follows that $E^F_\Delta \in \clKN$
and hence $E^F_\Delta$ is $\tau$-finite. Now, for any $h \in B_\csupp(a,b)$ there exists
a compact subset $\Delta$ of $(a,b)$ and a constant $c_1$ such that $\abs{h} \leq c_1 \chi_\Delta,$
so that $\abs{h(F)} \leq c_1 E_\Delta^F$ and hence $h(F) \in \LpN{1}.$
\end{proof}
\begin{lemma}\label{L: E Delta (F0+rK) is L1 bounded}
  If $F_0 \in \clFNab,$ $K=K^* \in \clKN,$ and if $\Delta$ is a compact subset of $(a,b),$
  then
  \\ (i) the function $\rr \in [0,1] \mapsto E_\Delta^{F_0+\rr K}$ takes values in $\LpN{1}$ and is
$\LpN{1}$-bounded;
  \\ (ii) there exists $R > 0$ such that the function
  $K \in B_R \cap \clKN \mapsto E_\Delta^{F_0+K}$ takes values in $\LpN{1}$ and is
$\LpN{1}$-bounded.
\end{lemma}
\begin{proof} (i) That $E_\Delta^{F_\rr} = E_\Delta^{F_0+\rr K} \in \LpN{1}$ follows from Lemmas
\ref{L: if F in clFNpm then h(F) is trace class} and \ref{L: E Delta (F0+rK) is L1 bounded}.
By (\ref{F: E Delta < c(1-F2)...}) we have
$E_\Delta^{F_r} \leq c_0 E^{F_r}_{[a,b]} (b-F_r)(F_r-a)$ for all $r \in [0,1]$
and hence by~\cite[Lemma 2.5]{FK86PJM}
$$
  \mu_t(E_\Delta^{F_r}) \leq c_0 \mu_t \brs{ E^{F_r}_{[a,b]} (b-F_r)(F_r-a) }
    \leq c_0 \mu_t \sqbrs{(b-F_r)(F_r-a)}.
$$
Since $(b-F_r)(F_r-a) = (b-F_0)(F_0-a) + r L_1 - r^2L_2,$
where $L_1, L_2 \in \clKN,$
we have by~\cite[Lemma 2.5(v)]{FK86PJM}
\begin{gather}\label{F: mu t E Delta Fr < ...}
  \mu_t(E_\Delta^{F_r})
     \leq c_0 \brs {  \mu_{t/3}\sqbrs{(b-F_0)(F_0-a)} + r\mu_{t/3}(L_1) + r^2 \mu_{t/3}(L_2) }
  \\ \leq c_0 \brs {  \mu_{t/3}\sqbrs{(b-F_0)(F_0-a)} + \mu_{t/3}(L_1) + \mu_{t/3}(L_2) }, \notag
\end{gather}
so that $\mu_t(E_\Delta^{F_r}) = \chi_{[0, \taubrs{E_\Delta^{F_r}}]}(t)$ is majorized for all $r \in [0,1]$
by a single function decreasing to 0 when $t \to \infty,$ since all three operators
$(b-F_0)(F_0-a),$ $L_1$ and $L_2$ are $\tau$-compact. The same argument as
in Lemma \ref{L: lemma with def of B R} now completes the proof.

(ii) If $F=F_0+K$ then $(b-F)(F-a) = (b-F_0)(F_0-a) + L,$
where $L = (b-F_0)K-K(F_0-a)-K^2 \in \clKN.$
Choose the number $R>0$ such that $\norm{K}<R$ implies $\norm{L} < 1.$ Then
by (\ref{F: mu t E Delta Fr < ...}) the function
$t \mapsto \mu_t(E_\Delta^{F+K}) = \chi_{[0, \taubrs{E_\Delta^{F+K}}]}(t)$ will be majorized by a single
function decreasing to a number $<1,$ so that the same argument as in
Lemma \ref{L: lemma with def of B R} again completes the proof.
\end{proof}
\begin{prop} \label{P: h(F+rK) is L1 continuous}
Let $F_0 \in \clFNab,$ $K=K^* \in \clKN,$ and let $h \in C^2_\csupp(a,b).$
Then
\\ (i) the function $r \in \mbR \mapsto h(F_0+\rr K)$ takes values in $\LpN{1}$ and
is $\LpN{1}$-continuous;
\\ (ii) there exists $R>0$ such that the function $K \in B_R \cap \clKN \mapsto h(F_0+K)$
takes values in $\LpN{1}$ and is $\LpN{1}$-continuous;
\end{prop}
\begin{proof} The proof of this proposition follows verbatim
the proof of Proposition \ref{P: f(B+rV) is L1 continuous}
with references to Lemmas
\ref{L: If F in clFNpm then F + K in clFNpm}, \ref{L: if F in clFNpm then h(F) is trace class}
and \ref{L: E Delta (F0+rK) is L1 bounded}
instead of Lemmas \ref{L: D+V has comp resolvent}, \ref{C: if B has comp rslv then f(B) is trace class}
and Corollary \ref{C: trace f(D0+V) is bounded}.
\end{proof}
\begin{lemma}\label{L: DOI is L1 continuous wrt F,F}
  If $F_1, F_2 \in \clFNab,$ $X \in \clKN_{sa}$ and
  $h \in C^3_\csupp(a,b),$ then the double operator integral
  $$
    T_{h^{[1]}}^{F_1, F_2}(X)
  $$
  takes values in $\LpN{1}$ and is $\LpN{1}$-continuous with respect
  to norm perturbations of $F_1$ and $F_2$ by $\tau$-compact operators.
\end{lemma}
The proof of this lemma is similar to that of Lemma~\ref{L: DOI is continuous wrt D1,D2}
with references to Lemma~\ref{L: E Delta (F0+rK) is L1 bounded}(ii)
and Proposition~\ref{P: h(F+rK) is L1 continuous}(ii)
instead of Corollary~\ref{C: trace f(D0+V) is bounded}
and Proposition~\ref{P: f(B+rV) is L1 continuous}.
\begin{thm} Let $\clN$ be a semifinite \vNa.
If $F_0 \in \clFNab,$ $h \in C^3_\csupp(a,b),$ then
the function $h \colon F \in F_0 + \saKN \mapsto h(F_0)+\saKN$
takes values in $\LpN{1}_{sa}.$ Moreover, it is affinely $(\clK,\clL^1)$-Fr\'echet
differentiable, the equality
$$
  D_{\clK,\clL^1} h(F) = T^{F,F}_{h^{[1]}}
$$
holds, and $D_{\clK,\clL^1} h(F)$ is $(\clK,\clL^1)$ continuous, so that
\begin{gather}\label{F: Newton-Leibnitz bdd case}
  h(F_{r_1}) - h(F_{r_0}) =  \int_{r_1}^{r_0} T^{F_r,F_r}_{h^{[1]}}(K)\,dr, \ \ r_0, r_1 \in \mbR,
\end{gather}
where $K \in \saKN,$ $F_r = F_0 + rK$ and the integral is in $\LpN{1}$-norm.
\end{thm}
The proof is similar to that of Theorem \ref{T: D L1 f(D) = T D,D}
with use of Proposition \ref{P: h(F+rK) is L1 continuous}(ii) and Lemma \ref{L: DOI is L1 continuous wrt F,F}
instead of Proposition \ref{P: f(B+rV) is L1 continuous} and Lemma \ref{L: DOI is continuous wrt D1,D2},
and therefore it is omitted.
\section{Spectral shift function} \label{S: spectral shift measure}
We will take an approach to the notion of spectral shift function
suggested by Birman-Solomyak formula (\ref{F: def of Xi}).
The key point is that once one appreciates that the spectral shift function
of M.\,G.\,Krein is related to spectral flow in a specific fashion one can
reformulate the whole approach to take advantage of what is known about spectral flow
as expounded for example in~\cite{BCPRSW}.
The theorem in~\cite{ACDS} which connects spectral flow and
the spectral shift function contains the germ of the idea but one needs
the technical machinery of the last Section to exploit this.

We now explain this different way to approach spectral shift theory
which is influenced by ideas from noncommutative geometry. %
\subsection{The unbounded case}
\subsubsection{Spectral shift measure}
In order to make our main definition we need to prove a preliminary result
which complements~\cite[Lemma 6.2]{ACDS}. The latter asserts that
the function $\gamma(\lambda,\rr) = \taubrs{ \VV E^{\HH_\rr}_\lambda}$ is measurable for every
$\VV \in \LpN{1}$ and $\HH=\HH^* \affl \clN.$
\begin{lemma}\label{L: trace of VE(l,s) is measurable II}
 Let $(\clN,\tau)$ be a semifinite von Neumann algebra and let $\HH=\HH^* \affl \clN$ have
 $\tau$-compact resolvent. If $\VV=\VV^* \in \clN$
 then the function $f \colon (a,b,\rr) \in \mbR^3 \mapsto \taubrs{\VV E_{(a,b)}^{\HH_\rr}}$ is measurable.
\end{lemma}
\begin{proof} W.l.o.g. we can assume that $\VV \geq 0.$
It is enough to prove that the function $f$ is measurable with respect to
the second variable $b$ and with respect to $\rr.$
Since $\taubrs{\VV E_{(a,b)}^{\HH_\rr}}=\taubrs{\sqrt \VV E_{(a,b)}^{\HH_\rr}\sqrt \VV},$ we know by~\cite[Lemma 5.9]{dPS04FA} that
it is enough to prove that the operator function $(\rr,b) \mapsto \sqrt \VV E_{(a,b)}^{\HH_\rr}\sqrt \VV$ is $so^*$-measurable.
By~\cite[Proposition 3.2]{ACDS} it is enough to prove that for any $\xi,\eta \in \hilb$ the scalar function
$\scal{\sqrt \VV E_{(a,b)}^{\HH_\rr}\sqrt \VV\xi}{\eta} = \Tr(\theta_{\sqrt \VV\xi, \sqrt \VV\eta}E_{(a,b)}^{\HH_\rr})$
is measurable, where $\theta_{\xi,\eta}(\zeta) := \scal{\xi}{\zeta}\eta.$ Since the operator
$\theta_{\sqrt \VV\xi, \sqrt \VV\eta}$ is trace class, the measurability of this function follows from
\cite[Lemma 6.2]{ACDS}.
\end{proof}
\begin{defn}\label{D: def of Xi} If $\HH_0=\HH_0^* \affl \clN$
has $\tau$-compact resolvent and if $\HH_1 = \HH_0 + \VV,$ $\VV \in \clN_{sa},$
then the \emph{spectral shift measure} for the pair $(\HH_0,\HH_1)$ is
defined to be the following Borel measure on $\mbR$
\begin{gather}\label{F: def of Xi}
  \Xi_{\HH_1,\HH_0}(\Delta) = \int_0^1 \taubrs{\VV E_{\Delta}^{\HH_\rr}}\,d\rr, \quad \Delta \in \clB(\mbR).
\end{gather}
The generalized function
\begin{gather}\label{F: def of xi}
  \xi_{\HH_1,\HH_0}(\lambda) = \frac d{d\lambda} \Xi_{\HH_1,\HH_0}(a,\lambda)
\end{gather}
is called the \emph{spectral shift distribution} for the pair $(\HH_0,\HH_1).$
\end{defn}
Evidently, this definition does not depend on a choice of $a.$
By Lemmas \ref{L: D+V has comp resolvent}, \ref{L: If D has comp rslv then E Delta if finite},
\ref{L: trace of VE(l,s) is measurable II} and Corollary  \ref{C: trace Br Delta is bounded}(i)
the measure $\Xi$ exists and is locally finite.

Our task now is to show that the spectral shift distribution is in fact a function of locally bounded
variation.
\subsubsection{Spectral shift function}
The main result we wish to establish next is that the spectral shift measure
is absolutely continuous with respect to Lebesgue measure. 
Moreover its density, which we previously referred to as the spectral shift distribution,
is in fact a function of locally bounded variation which we will then refer to as
the spectral shift function.
It is our extension of M.\,G.\,Krein's function to the setting of this paper.

Our method of proof is to first establish some trace formulae.
\begin{lemma}\label{L: trace of V alpha(D)}
(i) Let $\HH=\HH^* \affl \clN$
have $\tau$-compact resolvent.
A function $\alpha \in B(\mbR)$ is 1-summable with respect to the measure $\tau(E_\Delta^\HH)$
($\Delta \in \clB(\mbR)$),  if and only if $\alpha(\HH) \in \LpN{1}$ and in this case
$$
  \taubrs{\alpha(\HH)} = \int_\mbR \alpha(\lambda)\,\taubrs{dE_\lambda^{\HH}}.
$$
Furthermore, for any $\VV=\VV^* \in \clN$
the function $\alpha$ is 1-summable with respect to the measure $\tau(\VV E_\Delta^\HH),$
and
$$
  \taubrs{\VV\alpha(\HH)} = \int_\mbR \alpha(\lambda)\,\taubrs{\VV dE_\lambda^{\HH}}.
$$
\\ (ii)
Let $F \in \clFNab.$
A function $\alpha \in B(a,b)$ is 1-summable with respect to the measure $\tau(E_\Delta^F)$
($\Delta \in \clB(a,b)$), if and only if $\alpha(F) \in \LpN{1}$ and in this case
$$
  \taubrs{\alpha(F)} = \int_a^b \alpha(\lambda)\,\taubrs{dE_\lambda^{F}}.
$$
Furthermore, for any $\VV=\VV^* \in \clN$
the function $\alpha$ is 1-summable with respect to the measure $\tau(\VV E_\Delta^F),$
and
$$
  \taubrs{\VV\alpha(F)} = \int_a^b \alpha(\lambda)\,\taubrs{\VV dE_\lambda^{F}}.
$$
\end{lemma}
\begin{proof} We give only the proof of (i).
  W.l.o.g. we can assume that $\alpha$ is a non-negative function.
  If $\alpha$ is a simple function
  then the first part of the claim follows from Lemma \ref{L: If D has comp rslv then E Delta if finite}.
  Let $\alpha_n$ be an increasing sequence of simple non-negative functions,
  converging pointwise to $\alpha.$

  Then for each of the functions $\alpha_n$ the first equality is true.
  The supremum of the increasing sequence of
  non-negative operators $\alpha_n(\HH)$ is $\alpha(\HH)$
  and the supremum of the increasing sequence of numbers
  $\int_\mbR \alpha_n(\lambda)\,\taubrs{dE^{\HH}_\lambda}$
  is $\int_\mbR \alpha(\lambda)\,\taubrs{dE^{\HH}_\lambda}.$
  Hence, both non-negative numbers $\int_\mbR \alpha(\lambda)\,\taubrs{dE^{\HH}_\lambda}$
  and $\taubrs{\alpha(\HH)}$ are finite or infinite simultaneously, which proves the first part
  of the lemma.

  For the second part we can assume w.l.o.g. that $\VV \geq 0.$
  Then again the both parts of the second equality make
 sense and they are equal for simple functions.

  Since the measure $\taubrs{\VV E^\HH_\Delta} = \taubrs{\sqrt \VV E^\HH_\Delta \sqrt \VV},
  \ \Delta \in \clB(\mbR),$ is non-negative and the supremum of $\sqrt \VV \alpha_n(\HH) \sqrt \VV \in \LpN{1}$ is
  $\sqrt \VV \alpha(\HH) \sqrt \VV$ we have that
  \begin{multline*}
    \int_\mbR \alpha(\lambda)\,\taubrs{\VV dE^{\HH}_\lambda}
    = \lim_{n\to\infty} \int_\mbR \alpha_n(\lambda)\,\taubrs{\VV dE^{\HH}_\lambda}
    \\ = \lim_{n\to\infty}  \taubrs{\sqrt \VV \alpha_n(\HH) \sqrt \VV}
    = \taubrs{\sqrt \VV \alpha(\HH) \sqrt \VV},
  \end{multline*}
  so that $\int_\mbR \alpha(\lambda)\,\taubrs{\VV dE^{\HH}_\lambda}$ and
  $\taubrs{\VV \alpha(\HH)} = \taubrs{\sqrt \VV \alpha(\HH) \sqrt \VV}$ are finite or infinite simultaneously.
\end{proof}
We need the following version of Fubini's theorem.
\begin{lemma}\label{L: Fubini for variable measures}
(i)
For any self-adjoint operator $\HH \affl \clN$ with $\tau$-compact resolvent
and $\VV=\VV^* \in \clN,$ let $m_{\HH,\VV}(\Delta) = \taubrs{\VV E^{\HH}_\Delta}.$
Let $\HH_0=\HH_0^* \affl \clN$ have $\tau$-compact resolvent and let $\HH_r = \HH_0 + \rr\VV.$
If $g \in B_\csupp(\mbR),$ then
  \begin{gather}\label{F: Fubini for variable measures}
     \int_0^1 d\rr \int_\mbR g(\lambda) \,m_{\HH_\rr,\VV}(d\lambda) = \int_\mbR g(\lambda) \Xi_{\HH_1,\HH_0}(d\lambda).
  \end{gather}
\\ (ii)
For any $F \in \clFNab$ and $\VV=\VV^* \in \clN,$ let $m_{F,\VV}(\Delta) = \taubrs{\VV E^{F}_\Delta},$
$\Delta \in \clB(a,b).$
Let $F \in \clFNab$ and let $F_r = F_0 + \rr\VV.$
If $g \in B_\csupp(a,b),$ then
  \begin{gather*}\label{F: Fubini for variable measures bdd case}
     \int_0^1 d\rr \int_a^b g(\lambda) \,m_{F_\rr,\VV}(d\lambda) = \int_a^b g(\lambda) \Xi_{F_1,F_0}(d\lambda).
  \end{gather*}
\end{lemma}
\begin{proof} (See also~\cite[VI.2]{Ja}). We give only the proof of (i).
The measurability of the function $r \mapsto \int_\mbR g(\lambda) \,m_{\HH_\rr,\VV}(d\lambda)$
follows from Lemma \ref{L: trace of VE(l,s) is measurable II}.

Note, that both integrals are repeated ones. 
Let $\Omega \supseteq \supp(g)$ be a finite interval.
By Corollary \ref{C: trace Br Delta is bounded}\,(i) there exists $M>0$ such
that for all $\rr \in [0,1]$ we have $\abs{m_{\HH_\rr,\VV}(\Omega)} \leq M.$

If $g(\lambda) = \chi_\Delta(\lambda),$ \ $\Delta \in \clB(\Omega),$ then
\begin{gather*}
  \int_0^1 d\rr \int_\Omega \chi_\Delta(\lambda) \,m_{\HH_\rr,\VV}(d\lambda) = \int_0^1 m_{\HH_\rr,\VV}(\Delta)\,d\rr
  = \Xi(\Delta)
   = \int_\Omega \chi_\Delta(\lambda)\, \Xi(d\lambda).
\end{gather*}
So, (\ref{F: Fubini for variable measures}) is true for simple functions.
Let now $g$ be an arbitrary function from $B_\csupp(\Omega),$ let $\eps > 0$ and let $h$ be
a simple function such that $\norm{g - h}_\infty < \eps.$
Then the LHS of (\ref{F: Fubini for variable measures}) is equal to
$$
  \int_0^1 d\rr \int_\Omega (g-h)(\lambda) \,m_{\HH_\rr,\VV}(d\lambda)
      + \int_0^1 d\rr \int_\Omega h(\lambda) \,m_{\HH_\rr,\VV}(d\lambda)
  = (I) + (II),
$$
and the RHS of (\ref{F: Fubini for variable measures}) is equal to
$$
  \int_\Omega (g-h)(\lambda)\, \Xi(d\lambda) + \int_\Omega h(\lambda)\, \Xi(d\lambda) = (III) + (IV).
$$
We have $(II) = (IV).$
Further,
$\abs{(I)} \leq M \norm{g-h}_\infty  \leq M\eps$
and
$\abs{(III)} \leq M \norm{g-h}_\infty \leq M\eps.$
\end{proof}
The following theorem complements~\cite[Theorem 6.3]{ACDS}.
\begin{thm}\label{T: trace formula} 
  If $\HH=\HH^* \affl \clN$ has $\tau$-compact resolvent,
  if $\VV=\VV^* \in \clN,$ and if $\HH_1=\HH_0+\VV,$
  then the measure $\Xi_{\HH_1,\HH_0}$ is absolutely continuous, its density
  is equal to
  \begin{gather}\label{F: xi = trace E D0 - trace E D1 + const}
    \xi_{\HH_1,\HH_0}(\cdot) = \taubrs{E^{\HH_0}_{(a,\lambda]}} - \taubrs{E^{\HH_1}_{(a,\lambda]}} + \const
  \end{gather}
  for almost all $\lambda \in \mbR.$ Moreover, for all $f \in C_\csupp^3(\mbR)$ \ \
  $f(\HH_1) - f(\HH_0) \in \LpN{1}$ and
  \begin{gather}\label{F: trace formula}
      \taubrs{f(\HH_1) - f(\HH_0)} = \int_\mbR f'(\lambda) \xi_{\HH_1,\HH_0}(\lambda)\,d\lambda.
  \end{gather}
\end{thm}
\begin{proof}
By Lemma \ref{L: D+V has comp resolvent} and Corollary \ref{C: if B has comp rslv then f(B) is trace class}
$f(\HH_1) - f(\HH_0) \in \LpN{1}.$

By Lemma \ref{L: f = f1 - f2, where f1,f2 are nonnegative}
we need only consider the case of a non-negative function $f$ with $g := \sqrt f \in C^2_\csupp(\mbR).$

We have by (\ref{F: Newton-Leibnitz})
$$
  f(\HH_1) - f(\HH_0) = \int_0^1 T^{\HH_\rr,\HH_\rr}_{f^{[1]}}(\VV)\,d\rr,
$$
where the integral converges in $\LpN{1}$-norm.
Hence, it follows from Lemma \ref{L: BS representation for DOI} that
\begin{multline}\label{F: f(B1) - f(B0) = int 01 int Pi ...}
  f(\HH_1) - f(\HH_0)
   \\ = \int_0^1 \int_\Pi \brs {
    \alpha_1(\HH_\rr,\sigma)\VV \beta_1(\HH_\rr,\sigma) + \alpha_2(\HH_\rr,\sigma)\VV \beta_2(\HH_\rr,\sigma)
    }\,d\nu_g(\sigma)\,d\rr.
\end{multline}
Now, for a fixed $\sigma \in \Pi,$ we have
\begin{multline*}
   \taubrs{\alpha_1(\HH_\rr,\sigma)\VV \beta_1(\HH_\rr,\sigma) + \alpha_2(\HH_\rr,\sigma)\VV \beta_2(\HH_\rr,\sigma)}
  \\ = \taubrs{\VV (\alpha_1(\HH_\rr,\sigma)\beta_1(\HH_\rr,\sigma) + \alpha_2(\HH_\rr,\sigma)\beta_2(\HH_\rr,\sigma))}
  \\ = \int_\mbR \brs { \alpha_1(\lambda,\sigma)\beta_1(\lambda,\sigma)
      + \alpha_2(\lambda,\sigma)\beta_2(\lambda,\sigma) }\,\taubrs{\VV dE_\lambda^{\HH_\rr}},
\end{multline*}
where the last equality uses Lemma \ref{L: trace of V alpha(D)}.
(That $\alpha_1(\lambda)\beta_1(\lambda) + \alpha_2(\lambda)\beta_2(\lambda)$
belongs to $B_\csupp(\mbR)$ follows from Lemma \ref{L: BS representation for f(x) - f(y)})

Hence using (\ref{F: f(B1) - f(B0) = int 01 int Pi ...}),
and by Lemma \ref{L: lemma 3.10 ACDS}
applied to the finite measure space $\brs{[0,1] \mytimes \Pi, d\rr \mytimes \nu_g}$
our previous equality implies that we have:
\begin{multline*}
  A:=\taubrs{f(\HH_1) - f(\HH_0)}
  \\  = \int_0^1 \int_\Pi \taubrs{   \alpha_1(\HH_\rr,\sigma)\VV \beta_1(\HH_\rr,\sigma)
     + \alpha_2(\HH_\rr,\sigma)\VV \beta_2(\HH_\rr,\sigma)   }\,d\nu_g(\sigma)\,d\rr
  \\  = \int_0^1 \int_\Pi \int_\mbR ( \alpha_1(\lambda,\sigma)\beta_1(\lambda,\sigma)
            + \alpha_2(\lambda,\sigma)\beta_2(\lambda,\sigma)     )\,\taubrs{\VV dE_\lambda^{\HH_\rr}}\,d\nu_g(\sigma)\,d\rr.
\end{multline*}
Now, by Lemma \ref{L: trace of V alpha(D)}, Fubini's theorem, and Lemma \ref{L: BS representation for f(x) - f(y)} we have
\begin{gather*}
   A = \int_0^1 \int_\mbR \int_\Pi (\alpha_1(\lambda,\sigma)\beta_1(\lambda,\sigma)
            + \alpha_2(\lambda,\sigma)\beta_2(\lambda,\sigma) )
               \,d\nu_g(\sigma)\,\taubrs{\VV dE_\lambda^{\HH_\rr}}\,d\rr
   \\ = \int_0^1 \int_\mbR f'(\lambda) \,\taubrs{\VV dE_\lambda^{\HH_\rr}}\,d\rr.
\end{gather*}
Finally, by Lemma \ref{L: Fubini for variable measures} we have
\begin{gather}\label{F: A = int f' d Xi}
  A  = \int_\mbR f'(\lambda) \int_0^1 \,\taubrs{\VV dE_\lambda^{\HH_\rr}}\,d\rr
     = \int_\mbR f'(\lambda)\,d\Xi_{\HH_1,\HH_0}(\lambda).
\end{gather}

Let $f \in C^1_\csupp(\mbR)$ and take a point $a$ outside of the support of $f.$
Then we have (see~\cite[Proposition 8.5.5]{AB})
  \begin{multline}\label{F: trace f(D1) - f(D0) = int f' trace E ...}
     A = \tau (  f(\HH_1) - f(\HH_0)   ) = \tau( f(\HH_1) ) - \tau ( f(\HH_0) )
      \\ =  \int_\mbR f(\lambda) d \taubrs{E^{\HH_1}_{(a,\lambda]}}
          - \int_\mbR f(\lambda) d \taubrs{E^{\HH_0}_{(a,\lambda]}}
           \qquad ( \text{integrating by parts} )
      \\ = - \int_\mbR f'(\lambda) \brs{  \taubrs{E^{\HH_1}_{(a,\lambda]}} - \taubrs{E^{\HH_0}_{(a,\lambda]}} } \,d\lambda.
  \end{multline}
Comparing (\ref{F: trace f(D1) - f(D0) = int f' trace E ...}) and (\ref{F: A = int f' d Xi})
we see that $\Xi$ is absolutely continuous with density equal to
  \begin{gather}\label{F: trace E1 lamb - trace E0 lamb}
    \xi_{\HH_1,\HH_0}(\lambda) = \taubrs{E^{\HH_0}_{(a,\lambda]}} - \taubrs{E^{\HH_1}_{(a,\lambda]}} + \const.
  \end{gather}
\end{proof}
It is worth noting that the formula (\ref{F: trace formula}) does not
determine the function $\xi$ uniquely, but only up to an additive constant.
\begin{rems*} This theorem is an analogue of~\cite[Theorem 3.1]{ADS},
in which the existence and absolute continuity of the spectral shift measure
were proved for any self-adjoint operator $\HH$ affiliated with $\clN$
and $\tau$-trace class operator $\VV \in \clN.$
\end{rems*}
As a result of what we have proved to this point
we are now in a position to assert that in fact the spectral shift distribution is
an everywhere defined function and hence to change our terminology and refer to $\xi$ as a function.
Moreover this last
lemma enables one to modify  $\xi$ so as to make it
a function defined everywhere in a natural way.
\begin{defn}\label{D: xi at discontinuity point}
If the expression (\ref{F: trace E1 lamb - trace E0 lamb}) is continuous at a point $\lambda \in \mbR,$
then we define $\xi_{\HH_1,\HH_0}(\lambda)$ via formula  (\ref{F: trace E1 lamb - trace E0 lamb}).
Otherwise, we define the value of the spectral shift function $\xi$
at a discontinuity point to be half sum of left and right limits.
\end{defn}
\begin{cor} The spectral shift function $\xi$ is a function of locally bounded variation.
\end{cor}
\begin{proof} This is immediate because $\xi$ is the difference of two increasing functions by the last formula.
\end{proof}
\begin{lemma}\label{L: int ... d lambda = int ... dr} Let $\HH_0 \affl \clN$ be a self-adjoint operator with $\tau$-compact resolvent,
let $\VV \in \clN_{sa}$ and let $\HH_1 = \HH_0+\VV.$ If $f \in B_\csupp(\mbR)$ then
$$
  \int_{-\infty}^\infty f(\lambda) \xi_{\HH_1,\HH_0}(\lambda)\,d\lambda = \int_0^1 \taubrs{\VV f(\HH_r)}\,dr.
$$
\end{lemma}
\begin{proof} It follows from Lemma \ref{L: trace of V alpha(D)} and Lemma \ref{L: Fubini for variable measures}
that
  \begin{multline*}
    \int_0^1 \taubrs{\VV f(\HH_r)}\,dr
      = \int_0^1 \int_{-\infty}^\infty f(\lambda) \taubrs{\VV dE^{\HH_r}_\lambda}\,dr
      \\ = \int_{-\infty}^\infty f(\lambda) \int_0^1 \taubrs{\VV dE^{\HH_r}_\lambda}\,dr.
  \end{multline*}
\end{proof}
\subsubsection{The spectral shift function for unitarily equivalent operators}
  The situation where the operators $\HH$ and $\HH+\VV,$
are unitarily equivalent arises naturally in noncommutative geometry in the context of
spectral triples. One thinks of the unitary implementing the equivalence as a gauge transformation
by analogy with the study of gauge transformations of Dirac type operators.
It thus warrants special consideration especially in view of our
first result below.
\begin{thm}\label{T: SSF for unitary equivalent endpoints}
  Let $\HH$ be a self-adjoint operator affiliated with $\clN$ having $\tau$-compact resolvent
  and let $\VV=\VV^*\in \clN$ be such that the operators $\HH+\VV$ and $\HH$ are unitarily equivalent.
  Then the spectral shift function $\xi_{\HH+\VV,\HH}$ of the pair $(\HH+\VV,\HH)$
  is constant on $\mbR.$
\end{thm}
\begin{proof}
  The operators $f(\HH+\VV)$ and $f(\HH)$ are unitarily equivalent
  and for $f \in C^\infty_\csupp(\mbR)$ they are $\tau$-trace class by
  Corollary \ref{C: if B has comp rslv then f(B) is trace class}. Hence,
$$
  \taubrs{f(\HH+\VV) - f(\HH)} = 0,
$$
so that by Theorem \ref{T: trace formula} the equality
$$
  \int_\mbR f'(\lambda) \xi_{\HH+\VV,\HH}(\lambda)\,d\lambda = 0
$$
holds for any $f \in C^\infty_\csupp(\mbR).$
Now, integration by parts shows that $\xi'(\lambda)$ is zero as generalized
function on $\mbR,$ which by~\cite[Ch. I.2.6]{GSh}
implies that $\xi$ is equal to a constant function.
\end{proof}
Note, function $\xi$ in this theorem is equal to a constant function everywhere, not just almost everywhere.

Our second major result on the spectral shift function in this special context
is the following theorem. We shall show in Section \ref{S: spectral flow}
below that this theorem extends one of the main results of~\cite{CP2}.
\begin{thm}\label{T: tau(V f(D))} Let $\HH_0$ be a self-adjoint operator with $\tau$-compact resolvent,
affiliated with $\clN.$ Let $\VV=\VV^* \in \clN$ be such that the operators
$\HH_1=\HH_0+\VV$ and $\HH_0$ are unitarily equivalent.
If $f \in C^2_\csupp(\mbR)$ then
  \begin{gather}\label{F: tau(V f(D))}
    \xi_{\HH_1, \HH_0}(\mu) = C^{-1} \int_0^1 \tau(\VV f(\HH_\rr-\mu))\,d\rr, \quad \forall \ \mu \in \mbR,
  \end{gather}
where $C = \int_\mbR f(\lambda)\,d\lambda.$
\end{thm}
\begin{proof} For any fixed $\mu$ the operator $\HH_r - \mu$ has $\tau$-compact resolvent
by Lemma \ref{L: D+V has comp resolvent} and the function $r \mapsto \tau(\VV f(\HH_\rr-\mu))$
is continuous by Proposition \ref{P: f(B+rV) is L1 continuous}, so that the integral
on the RHS of (\ref{F: tau(V f(D))}) exists.
By Lemma \ref{L: int ... d lambda = int ... dr} and Theorem \ref{T: SSF for unitary equivalent endpoints} we have
$$
  \int_0^1 \taubrs{\VV f(\HH_\rr - \mu)}\,d\rr = \int_\mbR f(\lambda-\mu) \xi_{\HH_1,\HH_0}(\lambda)\,d\lambda
  = \xi(0) \int_\mbR f(\lambda)\,d\lambda.
$$
\end{proof}

\subsection{The bounded case}
As we remarked previously, our technique in the next Section for discussing spectral flow in
the unbounded case is to map into the space of bounded $\tau$-Fredholm operators.
We thus need to develop the theory described in the previous
subsections {\em ab initio} for the bounded case. Fortunately this is not a difficult task as the proofs
are much the same. As we will see, because we are considering bounded perturbations of our
unbounded operators, it suffices to consider compact perturbations in the bounded case.
\begin{defn}\label{D: def of Xi for F's} If $F_0 \in \clFNab,$
$K \in \clKN_{sa},$ $F_1 = F_0 + K$ and if $F_\rr = F_0 + \rr K,$
then the \emph{spectral shift measure} for the pair $(F_0,F_1)$ is
defined to be the following Borel measure on $(a,b)$
\begin{gather}\label{F: def of Xi for F's}
  \Xi_{F_1,F_0}(\Delta) = \int_0^1 \taubrs{K E_{\Delta}^{F_\rr}}\,d\rr, \quad \Delta \in \clB(a,b).
\end{gather}
The generalized function
\begin{gather}\label{F: def of xi for F's}
  \xi_{F_1,F_0}(\lambda) = \frac d{d\lambda} \Xi_{F_1,F_0}(c,\lambda), \ \ c \in (a,b),
\end{gather}
is called the \emph{spectral shift distribution} for the pair $(F_0,F_1).$
\end{defn}
Evidently, this definition does not depend on a choice of $c \in (a,b).$
The measurability of the function $r \mapsto \taubrs{K E_{\Delta}^{F_\rr}}$
may be established following the argument of Lemma \ref{L: trace of VE(l,s) is measurable II},
using Lemma \ref{L: E Delta (F0+rK) is L1 bounded}.
It follows that the measure $\Xi$ exists and is locally-finite on $(a,b).$
\begin{prop}\label{P: xi function for pair F0 and F1} 
If $F_0 \in \clFNab,$ $K \in \clKN$ and if $F_1 = F_0 + K,$ then
  \\ (i) the measure $\Xi_{F_1,F_0}$ is absolutely continuous
       and its density is equal to
     $$
       \xi_{F_1,F_0}(\lambda) = \taubrs{E^{F_0}_{(c,\lambda]} - E^{F_1}_{(c,\lambda]}} + \const, \ \lambda \in (c,b),
     $$
     where $c$ is an arbitrary number from $(a,b)$;
  \\ (ii) there exists a unique function
      $\xi_{F_1,F_0}(\cdot)$ of locally bounded variation on $(a,b),$
      such that for any $h \in C^2_\csupp(a,b)$
      the following equality holds true
      $$
        \taubrs{\h(F_1) - \h(F_0)} = \int_a^b \h'(\lambda) \xi_{F_1,F_0}(\lambda)\,d\lambda.
      $$
\end{prop}
The proof is identical to the proof of Theorem \ref{T: trace formula},
with references to
\ref{L: If F in clFNpm then F + K in clFNpm}, \ref{L: if F in clFNpm then h(F) is trace class},
(\ref{F: Newton-Leibnitz bdd case})
instead of \ref{L: D+V has comp resolvent}, \ref{C: if B has comp rslv then f(B) is trace class},
(\ref{F: Newton-Leibnitz})
and hence we  omit it.

\begin{cor}\label{C: if F0 and F1 are unit equiv then xi = const}
  In the setting of Proposition \ref{P: xi function for pair F0 and F1},
  if $F_0$ and $F_1$ are unitarily equivalent,
  then $\xi_{F_1,F_0}$ is constant on $(a,b).$
\end{cor}
The proof is similar to the proof of Theorem \ref{T: SSF for unitary equivalent endpoints}.
%
\begin{defn}\label{D: def of xi F F}
  We redefine the function $\xi_{F_1,F_0}$ at discontinuity points
  to be half the sum of the left and the right limits of the RHS of the last equality.
\end{defn}
Thus, the function $\xi_{F_1,F_0}$ is defined everywhere on $(a,b).$
\begin{lemma} \label{L: int ... d lambda = int ... dr for bdd case}
If $F_0 \in \clFNab,$ $K \in \clKN,$ if $F_r = F_0 + rK, \ r \in [0,1]$
and if $h \in B_\csupp(a,b)$ then
$$
  \int_\mbR h(\lambda)\xi_{F_1,F_0}(\lambda)\,d\lambda = \int_0^1 \tau(K h(F_r))\,dr.
$$
\end{lemma}
This Lemma and its proof are bounded variants
of Lemma \ref{L: int ... d lambda = int ... dr},
so we omit the details.
\section{Spectral flow} \label{S: spectral flow}
\subsection{The spectral flow function}
We wish to avoid a long excursion into the analytic theory of spectral flow preferring the reader
to read the early Sections of~\cite{BCPRSW} for the relevant background and definitions.
With those prerequisites it is possible to appreciate that the following notions are well defined.

We now introduce the spectral flow function on the real line.
As usual we have $\HH_0=\HH_0^* \affl \clN$ with $\tau$-compact resolvent
and $\HH_1 =\HH_0 +\VV$ with $\VV$ bounded.
Then we define the spectral flow function to be
$$
  \mu \mapsto \sflow(\mu,\HH_0,\HH_1), \quad \mu \in \mbR,
$$
where $\sflow(\mu,\HH_0,\HH_1)$ is spectral flow from
$\HH_0-\mu$ to $\HH_1-\mu$ along the path $\HH_r-\mu= \HH_0-\mu+rV$.
We remark that the homotopy invariance of spectral flow means that, in the
affine space of bounded perturbations of a fixed operator $\HH_0$,
spectral flow does not depend on the choice of continuous path
but only on the endpoints. However one does need to make precise what one means by
continuity in a path parameter in the unbounded case.
We define the space of
self adjoint, unbounded $\tau$-Fredholm operators
to be those operators
that under the map
$$\HH \mapsto F_{\HH}=\HH(1+\HH)^{-1/2}$$ are sent to bounded $\tau$-Fredholm
operators in $\clN$.
This definition is equivalent to the usual definition for unbounded Fredholm
operators in a semifinite von Neumann algebra (more details on this can be found in \cite{CPRS3}.)

A path of unbounded $\tau$-Fredholm operators is said to be continuous if its image under this
map is continuous in the norm topology on the bounded $\tau$-Fredholm operators.
This topology is usually called the Riesz topology and it is different from the graph norm topology
used in \cite{BLP}. A more detailed discussion of topologies on the set of
unbounded self-adjoint Fredholm operators and the relevance of these for
spectral flow may be found in \cite{Wa2}.

We note that the condition of having compact resolvent implies the $\tau$-Fredholm property
for the unbounded operator.

We now recall, to provide some motivation for the point of view of this Section, some ideas
from~\cite{CP2} which is couched in the language of
 noncommutative geometry~\cite{CoNG}. In~\cite{CP2} the condition of theta summability
($\tau(e^{-t\HH_0^2})<\infty$ for all real $t>0$) is imposed and then the main result of~\cite{CP2}
is an analytic formula for spectral flow from $\HH_0$ to $\HH_1=\HH_0+\VV$.
The ideas behind  this formula go back to~\cite{APS76} and a more complete history may be
found in~\cite{BCPRSW}.
The formula of~\cite{CP2} is as follows
\begin{multline}\label{F: Carey-Phillips}
  \sflow(\HH_0,\tHH) = \epspi \int_0^1 \taubrs{\VV e^{-\eps{\HH_\rr}^2}}\,d\rr
     \\ + \frac12 \brs{\eta_\eps(\HH_1) - \eta_\eps(\HH_0)}
     +\frac12 \taubrs{[\ker(\HH_1)]-[\ker(\HH_0)]},
\end{multline}
where $\eta_\eps$ is a `truncated eta invariant'.
For an unbounded self-adjoint operator $\HH$
for which $e^{-t\HH^2}$ is $\tau$-trace class for all $t>0,$ it is defined
in~\cite[Definition 8.1]{CP2} following~\cite{Ge93Top} by
$$
  \eta_\eps(\HH)
:= \frac 1{\sqrt \pi} \int_\eps^\infty \taubrs{\HH e^{-t\HH^2}} \
t^{-1/2}\,dt.
$$
When the endpoints are unitarily equivalent the two $\eta_\eps$ terms and the kernel correction terms
in the formula for spectral flow cancel.

In~\cite{ACDS} we showed that in the theta summable case and for trace class perturbations
the spectral shift function and the spectral flow function differ only by kernel correction terms.
Our aim in this paper is to show that this is the case more generally and in fact to
go further.  We will demonstrate that the spectral shift
function provides a way to prove more general analytic
formulae for spectral flow than is achieved in~\cite{CP2}.

The plan of this Section is to first establish a geometric framework that
is analogous to that of Getzler~\cite{Ge93Top} and~\cite{CP2}.
Then we derive analytic formulae for spectral flow in terms of the spectral shift function
in both the case of paths of bounded $\tau$-Fredholm operators and for paths of unbounded
$\tau$-Fredholm operators. The starting point is basically the same as that of~\cite{CP2} in that
we need a formula for the relative index of two projections.
This next result is a strengthening of~\cite[Theorem 4.1]{CP2}.
\begin{lemma}\label{L: sf lemma for P and Q}
Let $P$ and $Q$ be two    
projections in the semifinite \vNa\ $\clN$
and let $a < 0 < b$ be two real numbers.
Let $\af$ be a continuous function such that for any
$s \in [0,\frac {(b-a)^2}4]$  \ $\af(s(P-Q)^2)$ is $\tau$-trace class. Then $F_0 = (b-a)P+a$
and $F_1 = (b-a)Q+a$ are self-adjoint $\tau$-Fredholm operators from
$\clFNab$ as is the path $F_\rr = F_0 + \rr(F_1-F_0),$ and
$$
  \sflow(\set{F_\rr}) = C_{a,b}^{-1} \int_0^1 \taubrs{\dot F_\rr \af\sqbrs{(b-F_\rr)(F_\rr-a)}}\,d\rr,
$$
where $C_{a,b} = \int_0^1 (b-a) \af\Brs{(b-a)^2(\rr-\rr^2)} \,d\rr$
is a constant,
and the derivative $\dot F_\rr$ is $\norm{\cdot}$-derivative.
\end{lemma}
\begin{proof} We have
$$
  \dot F_\rr = F_1 - F_0 = (b-a)(Q-P)
$$
and
$$
  (b-F_\rr)(F_\rr-a) = (b-a)^2\rr(1-\rr)(Q-P)^2,
$$
so that by assumption $\af\sqbrs{(b-F_\rr)(F_\rr-a)}$ is $\tau$-trace class for $\rr \in [0,1].$
For each $\rr \in (0,1)$ define
$$
  f_\rr(x) = (b-a)x \af\Brs{(b-a)^2(\rr-\rr^2)x^2}.
$$
Then
\begin{multline}
  \int_0^1 \taubrs{\dot F_\rr \af\sqbrs{(b-F_\rr)(F_\rr-a)}}\,d\rr
  \\ = \int_0^1 \tau\BRS{(b-a)(Q-P) \af\sqbrs{(b-a)^2\rr(1-\rr)(Q-P)^2}} \,d\rr,
\end{multline}
and by~\cite[Theorem 3.1]{CP2} we have
\begin{multline*}
  \int_0^1 \taubrs{\dot F_\rr \af\sqbrs{(b-F_\rr)(F_\rr-a)}}\,d\rr
    = \int_0^1 \tau(f_\rr(Q-P))\,d\rr
  \\  = \int_0^1 f_\rr(1) \ind(QP)\,d\rr
    = C_{a,b} \ind (QP) = C_{a,b} \sflow(\set{F_\rr}),
\end{multline*}
where the last equality follows from $a<0<b$ and the definition of spectral flow.
\end{proof}
\begin{rems*} In the proof of this lemma we use~\cite[Theorem 3.1]{CP2} for functions $f$ without the condition
$f(1) \neq 0.$ But an inspection of the proof of this theorem shows that this condition becomes superfluous, if we
rewrite the statement of this theorem as $f(1) \ind(QP) = \tau [f(P-Q)],$ which for functions $f$ with $f(1) = 0$
becomes just $0=0.$
\end{rems*}

\subsection{Spectral flow one-forms in the unbounded case}
The strategy of~\cite{CP2} is geometric and follows ideas of~\cite{Ge93Top}.
The first step in this strategy is summarized in
Proposition \ref{P: one form tau(Xf(B)) is closed} in preparation for which we
need an explicit formula for the derivative of function of a path of operators.
The method by which this is achieved in~\cite{CP2} does not apparently generalise
sufficiently far to cover the situations considered in this paper.
The double operator integral approach of Section 2 overcomes this problem.
\begin{lemma} \label{L: trace of Y TDD X = trace of X TDD Y}
  Let $\HH=\HH^* \affl \clN,$ let $X,Y \in \clN$ and let $f \in C^2_\csupp(\mbR)$
be a non-negative function such that $g:=\sqrt f \in C^2_\csupp(\mbR).$
If $Y T^{\HH,\HH}_{f^{[1]}}(X)$ and $X T^{\HH,\HH}_{f^{[1]}}(Y)$ are both $\tau$-trace class then
$$
  \taubrs{Y T^{\HH,\HH}_{f^{[1]}}(X)} = \taubrs{X T^{\HH,\HH}_{f^{[1]}}(Y)}.
$$
\end{lemma}
\begin{proof}
  By Lemma \ref{L: BS representation for DOI} we have
  \begin{multline*}
    A = \taubrs{Y \, T^{\HH,\HH}_{f^{[1]}}(X)}
     \\       = \taubrs{Y  \int_\Pi \brs {\alpha_1(\HH,\sigma) X \beta_1(\HH,\sigma)
               + \alpha_2(\HH,\sigma) X \beta_2(\HH,\sigma) } \, d\nu_g(\sigma)}
     \\ = \taubrs{Y  \int_\Pi \brs {
               e^{i(s_1-s_0)\HH} g(\HH) X e^{is_1\HH} + e^{i(s_1-s_0)\HH} X g(\HH) e^{is_1\HH} }
                       \, d\nu_g(s_0,s_1)}.
  \end{multline*}
  Making the change of variables $s_1-s_0 = t_0, \ s_1 = t_1,$ and using (\ref{F: def of nu f}) we have
  \begin{multline*}
    \int_\Pi \brs { \alpha_1(\HH,\sigma) X \beta_1(\HH,\sigma) + \alpha_2(\HH,\sigma) X \beta_2(\HH,\sigma) } \, d\nu_g(\sigma)
    \\ = \sgn(t_0+t_1)\frac i{\sqrt{2\pi}} \int_{\set{(t_0,t_1) \in \mbR^2,\, t_0t_1 \geq 0}}
                 \big[e^{it_0\HH} g(\HH) X e^{it_1\HH}
      \\ + e^{it_0\HH} X g(\HH) e^{it_1\HH} \big] \fourier{g}(t_0+t_1)\,dt_0 dt_1
    \\  = \sgn(t)\frac i{\sqrt{2\pi}} \int_\mbR \Big( \int_{\Sigma_t} \Big( e^{it_0\HH} g(\HH) X e^{it_1\HH}
       \\ + e^{it_0\HH} X g(\HH) e^{it_1\HH} \Big) \,dl_t\Big)\fourier{g}(t)\,dt,
  \end{multline*}
  where $\Sigma_t = \set{(t_0,t_1) \in \mbR^2 \colon t_0t_1 \geq 0, \ t_0 + t_1 = t}$ and $dl_t$ is the Lebesgue measure on $\Sigma_t.$
  Thus, by Fubini's theorem~\cite[Lemma 3.8]{ACDS}
  \begin{multline*}
    A = \sgn(t)\frac i{\sqrt{2\pi}} \tau \Big(Y \, \int_\mbR \Big( \int_{\Sigma_t}
          \big(e^{it_0\HH} g(\HH) X e^{it_1\HH}
            \\ + e^{it_0\HH} X g(\HH) e^{it_1\HH} \big)   \,dl_t\Big) \fourier{g}(t)\,dt\Big)
    \\ = \sgn(t)\frac i{\sqrt{2\pi}}  \int_\mbR \Big( \int_{\Sigma_t}
       \tau \Big(Y e^{it_0\HH} g(\HH) X e^{it_1\HH}
       \\ + Y e^{it_0\HH} X g(\HH) e^{it_1\HH} \Big)  \,dl_t\Big) \fourier{g}(t)\,dt
    \\ = \sgn(t)\frac i{\sqrt{2\pi}}  \int_\mbR \Big( \int_{\Sigma_t}
          \tau \Big(X e^{it_1\HH}Y g(\HH) e^{it_0\HH} \\
          + X e^{it_1\HH} g(\HH) Y e^{it_0\HH} \Big) \,dl_t\Big) \fourier{g}(t)\,dt,
  \end{multline*}
  where the trace and the integral can be interchanged by Lemma \ref{L: lemma 3.10 ACDS}.
  The integral above coincides with $\taubrs{X T^{\HH,\HH}_{f^{[1]}}(Y)}.$
\end{proof}
The key geometric idea is to regard the analytic formula for spectral flow
of~\cite{Ge93Top} and~\cite{CP2} as expressing it as
an integral of a one form. As we are dealing with an affine space the geometry is easy to
invoke as we see in the next result.

\begin{prop}\label{P: one form tau(Xf(B)) is closed}
  Let $\HH$ be a self-adjoint operator affiliated with~$\clN,$ having $\tau$-compact resolvent
  and let $f \in C^3_\csupp(\mbR).$
  Let $\alpha= \alpha^f$ be a 1-form on the affine space $\HH_0 + \clN_{sa}$ defined
  at the point  $\HH \in \HH_0 + \clN_{sa}$  by the formula
  \begin{gather}\label{F: alpha B 1-form}
    \alpha_\HH^f(X) = \taubrs{ X f(\HH)}, \ \ X \in \clN_{sa}, \ \ \HH \in \HH_0 + \clN_{sa}.
  \end{gather}
Then $\alpha$ is a closed 1-form, and, hence, also exact by the Poincar\'e lemma.
\end{prop}
\begin{proof} The proof follows mainly the lines of~\cite{CP98CJM}, with necessary adjustments.
As usual, we can assume that $ f \geq 0$ and $g := \sqrt f \in C^3_\csupp(\mbR).$
We note that the operator $X f(\HH)$ is $\tau$-trace class,
so that the 1-form $\alpha$ is well-defined.
Now, by the definition of the
exterior differential, for $X,Y \in \clN,$ we have
$$
  d \alpha_\HH(X,Y) = \LieDer{X} \alpha_\HH(Y) - \LieDer{Y} \alpha_\HH(X) - \alpha_\HH([X,Y]),
$$
where $\LieDer{X}$ is the Lie derivative along the constant vector field $X.$
Since the space $\HH_0+\clN_{sa}$ is flat, we have $[X,Y]=0.$ So, we have to prove that
$
  \LieDer{X} \alpha_\HH(Y) = \LieDer{Y} \alpha_\HH(X).
$
It follows from Theorem \ref{T: D L1 f(D) = T D,D} that
\begin{gather*}
  A := \LieDer{X} \alpha_\HH(Y) = \dds \alpha_{\HH+sX}(Y) = \dds \taubrs{Y f(\HH+sX)}
  \\ = \taubrs{Y D_{\clN,\clL^1} f(\HH)(X)} = \taubrs{Y \, T^{\HH,\HH}_{f^{[1]}}(X)}.
\end{gather*}
Hence, by Lemma \ref{L: trace of Y TDD X = trace of X TDD Y}
$$
  \LieDer{X} \alpha_\HH(Y)
   = \taubrs{Y \, T^{\HH,\HH}_{f^{[1]}}(X)}
   = \taubrs{X \, T^{\HH,\HH}_{f^{[1]}}(Y)}
   = \LieDer{Y} \alpha_\HH(X),
$$
which implies that $\alpha_\HH$ is a closed 1-form.
\end{proof}
Though closedness of a 1-form already should imply its exactness by
the \Poincare\ lemma and contractibility of the domain we follow~\cite{CP2}
and give an independent proof of exactness.

\begin{defn} Let $\HH_0$ be a fixed self-adjoint operator with $\tau$-compact resolvent
affiliated with $\clN,$ and let $f \in C_\csupp(\mbR).$ We define the function $\theta^f$
on the affine space $\HH_0 + \clN$ by the formula
$$
  \theta_\HH^f = \int_0^1 \taubrs{\VV f(\HH_\rr)}\,d\rr,
$$
where $\HH \in \HH_0+\clN_{sa},$ $\VV = \HH-\HH_0$ and $\HH_\rr = \HH_0 + \rr\VV.$
Measurability of the function $r \mapsto \taubrs{\VV f(\HH_\rr)}$ follows from Lemma \ref{L: trace of VE(l,s) is measurable II}.
\end{defn}
\begin{prop}\label{P: d theta = alpha} Let $f \in C^3_\csupp(\mbR)$ and let $X \in \clN.$ Then
$$
  d\theta_\HH^f(X) = \alpha_\HH^f(X).
$$
\end{prop}
\begin{proof} W.l.o.g. we can assume that $X$ is self-adjoint.
By definitions
  \begin{multline*}
     (A):=d\theta_\HH^f(X)
     \\ = \dds \theta^f_{\HH+sX} = \dds \int_0^1 \taubrs{(\VV+s X) f(\HH_\rr + s \rr X)}\,d\rr
     \\ = \lims{s\to 0} \frac 1s \int_0^1 \taubrs{(\VV+s X) f(\HH_\rr + s \rr X) - \VV f(\HH_\rr)}\,d\rr
     \\ = \lims{s\to 0} \int_0^1 \tau\Brs{X f(\HH_\rr + s \rr X)}\,d\rr
          + \lims{s\to 0} \frac 1s \int_0^1 \tau\BRS{\VV \Brs{f(\HH_\rr+s \rr X) - f(\HH_\rr)} }\,d\rr.
  \end{multline*}
  The first summand of this sum by Proposition \ref{P: f(B+rV) is L1 continuous} is equal to
  $$
    \int_0^1 \tau\Brs{X f(\HH_\rr)}\,d\rr.
  $$
By Proposition \ref{P: DalKreinII}(i)
the second summand is equal to
  \begin{multline*}
     \lims{s\to 0} \frac 1s \int_0^1 \taubrs{\VV T^{\HH_\rr+srX,\HH_\rr}_{f^{[1]}}(s \rr X)}\,d\rr
       = \lims{s\to 0} \int_0^1 \taubrs{\VV T^{\HH_\rr+srX,\HH_\rr}_{f^{[1]}}(\rr X)}\,d\rr
      \\ = \int_0^1 \taubrs{\VV\, T^{\HH_\rr,\HH_\rr}_{f^{[1]}}(\rr X)}\,d\rr
       = \int_0^1 \taubrs{X\, T^{\HH_\rr,\HH_\rr}_{f^{[1]}}(\VV)}\rr\,d\rr,
  \end{multline*}
where the second equality follows from Lemma \ref{L: DOI is continuous wrt D1,D2}
and the last equality follows from Lemma \ref{L: trace of Y TDD X = trace of X TDD Y}.
Hence, by Lemma \ref{L: lemma 3.10 ACDS}
  \begin{multline*}
    (A) = \int_0^1 \taubrs{ X \big[ f(\HH_\rr) + \rr\,T^{\HH_\rr,\HH_\rr}_{f^{[1]}}(\VV)\big]}\,d\rr
    \\ = \taubrs{ X \int_0^1 \big[ f(\HH_\rr) + \rr\,T^{\HH_\rr,\HH_\rr}_{f^{[1]}}(\VV)\big]\,d\rr},
  \end{multline*}
  where the integral on the RHS is a $so^*$-integral.
  By Proposition \ref{P: f(B+rV) is L1 continuous} and Lemma \ref{L: DOI is continuous wrt D1,D2}
  the function $\rr \in [0,1] \mapsto f(\HH_\rr) + \rr\,T^{\HH_\rr,\HH_\rr}_{f^{[1]}}(\VV) \in \LpN{1}$
  is $\LpN{1}$-continuous, so that the last integral
  $$
     (B) := \int_0^1 \big[ f(\HH_\rr) + \rr\,T^{\HH_\rr,\HH_\rr}_{f^{[1]}}(\VV)\big]\,d\rr
  $$
  can be considered as Riemann integral.
  Let $0 = \rr_0 < \rr_1 < \ldots < \rr_n = 1$ be the partition of $[0,1]$ into $n$ segments of equal length $\frac 1n.$
  By the argument used in the proof of~\cite[Theorem 5.8]{ACDS} it can be shown that
  the $\LpN{1}$-norm of $T^{\HH_{\rr_j},\HH_{\rr_j}}_{f^{[1]}}(\VV) - T^{\HH_{\rr},\HH_{\rr}}_{f^{[1]}}(\VV),$
  \ $r \in [r_{j-1},r_j],$ has order $\frac 1n.$
  Hence
  \begin{multline*}
    \clL^1\mbox{-}\lims{n \to \infty} \frac 1n \sums{j=1}^n \brs{ \frac jn T^{\HH_{\rr_j},\HH_{\rr_j}}_{f^{[1]}}(\VV)
      - j \int_{\rr_{j-1}}^{\rr_j} T^{\HH_{\rr},\HH_{\rr}}_{f^{[1]}}(\VV)\,d\rr }
    \\ = \clL^1\mbox{-}\lims{n \to \infty} \frac 1n \sums{j=1}^n
      j \int_{\rr_{j-1}}^{\rr_j} \brs{ T^{\HH_{\rr_j},\HH_{\rr_j}}_{f^{[1]}}(\VV) - T^{\HH_{\rr},\HH_{\rr}}_{f^{[1]}}(\VV)}\,d\rr
      = 0,
  \end{multline*}
  so that by formula (\ref{F: Newton-Leibnitz}) applied to the pair $(\HH_{\rr_{j-1}},\HH_{\rr_j})$ we have
  \begin{multline*}
    (B) = \clL^1\mbox{-}\lims{n \to \infty} \frac 1n \sums{j=1}^n \brs{ f(\HH_{\rr_{j-1}}) + \frac jn T^{\HH_{\rr_j},\HH_{\rr_j}}_{f^{[1]}}(\VV) }
    \\ = \clL^1\mbox{-}\lims{n \to \infty} \frac 1n \sums{j=1}^n \brs{ f(\HH_{\rr_{j-1}}) + j\brs{ f(\HH_{\rr_{j}}) - f(\HH_{\rr_{j-1}})}}
    \\ +  \clL^1\mbox{-}\lims{n \to \infty} \frac 1n \sums{j=1}^n \brs{ \frac jn T^{\HH_{\rr_j},\HH_{\rr_j}}_{f^{[1]}}(\VV)
      - j \int_{\rr_{j-1}}^{\rr_j} T^{\HH_{\rr},\HH_{\rr}}_{f^{[1]}}(\VV)\,d\rr }
    \\ = \clL^1\mbox{-}\lims{n \to \infty} \frac 1n \sums{j=1}^n \brs{ j f(\HH_{\rr_{j}}) - (j-1)f(\HH_{\rr_{j-1}})} = f(\HH_1).
  \end{multline*}
\end{proof}
The argument before~\cite[Proposition 1.5]{CP98CJM} now implies the following corollary.
\begin{cor}\label{C: integral is independent of path}
   The integral of the 1-form $\alpha^f$ along a piecewise continuously differentiable path
   $\Gamma$ in $\HH_0 + \clN$ depends only on the endpoints of the path~$\Gamma.$
\end{cor}
\begin{prop}\label{P: xi DD is additive}
If a self-adjoint operator $\HH_0$ affiliated with $\clN$ has $\tau$-compact resolvent,
$\HH_1, \HH_2 \in \HH_0 + \clN_{sa},$ then for all $\lambda \in \mbR$
  $$
     \xi_{\HH_2, \HH_0}(\lambda) = \xi_{\HH_2, \HH_1}(\lambda) + \xi_{\HH_1, \HH_0}(\lambda).
  $$
\end{prop}
\noindent{\bf Remark}. We emphasize that this additivity property is not
almost everywhere in the spectral variable but in fact holds everywhere.

\begin{proof}
It follows from (\ref{F: xi = trace E D0 - trace E D1 + const}) that
$$
  \xi_{\HH_2, \HH_0}(\lambda) = \xi_{\HH_2, \HH_1}(\lambda) + \xi_{\HH_1, \HH_0}(\lambda) + C,
$$
where~$C$ is a constant. Multiplying both sides of this equality by a positive $C^2_\csupp$-function~$f,$
and integrating it, by Lemma \ref{L: int ... d lambda = int ... dr} we get
\begin{gather*}
  \int_{\Gamma_{\HH_2,\HH_0}} \alpha^f
    = \int_{\Gamma_{\HH_2,\HH_1}} \alpha^f + \int_{\Gamma_{\HH_1,\HH_0}} \alpha^f + C\int_\mbR f(\lambda)\,d\lambda,
\end{gather*}
where $\Gamma_{\HH_i,\HH_j}$ is the straight line path connecting operators $\HH_i$ and $\HH_j.$
The last equality and Corollary \ref{C: integral is independent of path} imply that $C=0.$
\end{proof}
\subsection{Spectral flow one-forms in the bounded case}
Since we obtain our unbounded spectral flow formula from a bounded one
we need to study the map
$\HH \mapsto F_\HH= \HH(1+\HH^2)^{-1/2}$ which takes the space
of unbounded self adjoint operators
with $\tau$-compact resolvent to the space $\clFNpm$
of bounded $\tau$-Fredholm operators
$F$ satisfying $1-F^2 \in \clKN.$

Let $F_0 \in \clFNab,$ let $h \in C^2_\csupp(a,b)$ and let
$K = F-F_0, \ F_\rr:=F_0 + \rr K.$
We define
a 0-form $\theta$ and a 1-form $\alpha^\h$ on the affine space $\clA_{F_0}$ by the formulae
$$
  \theta_F^\h = \int_0^1 \taubrs{K \,\h(F_\rr)}\,dr,
$$
and
$$
  \alpha_F^\h(X) = \taubrs{X \,\h(F)}, \ \ X \in \clKN.
$$
By Lemmas \ref{L: If F in clFNpm then F + K in clFNpm} and \ref{L: if F in clFNpm then h(F) is trace class},
the operators $\h(F_\rr)$ and $h(F)$ are $\tau$-trace class, so that
the forms $\theta^h$ and $\alpha^h$ are well-defined.
\begin{prop}\label{P: one form tau(Xf(F)) is closed}
If $F_0 \in \clFNab$ and if $\h \in C^2_\csupp(a,b),$
then
$$
  d\theta_F^h(X) = \alpha_F^h(X),
$$
where $X \in \clKN,$ so that the 1-form $\alpha_F^h$ is exact.
\end{prop}
\begin{proof}
The proof follows verbatim the proof of Proposition  \ref{P: d theta = alpha},
with references to Proposition \ref{P: h(F+rK) is L1 continuous} and
Lemma \ref{L: DOI is L1 continuous wrt F,F} instead of
Proposition \ref{P: f(B+rV) is L1 continuous} and Lemma \ref{L: DOI is continuous wrt D1,D2}.
\end{proof}
As in the unbounded case we get the following
\begin{cor} \label{C: integral is independent of path (bdd)}
  The integral of the one-form $\alpha^h$ depends only on the end-points.
\end{cor}
\begin{cor}\label{C: xi FF is additive}
   Let $F_j \in \clFNab, \ j = 0,1,2,$ such that $F_2 - F_1, F_1 - F_0  \in \clKN.$
Then for any $\lambda \in (a,b)$ the following equality holds true
$$
  \xi_{F_2,F_0}(\lambda) = \xi_{F_2,F_1}(\lambda) + \xi_{F_1,F_0}(\lambda).
$$
\end{cor}
We omit the proof as it is similar to that
of Proposition \ref{P: xi DD is additive}.

Let
$$
  \phi(\lambda) = \lambda \brs{1+\lambda^2}^{-1/2}.
$$
It is easy to see that if $\HH=\HH^* \affl \clN$ is an operator with $\tau$-compact resolvent, then
the operator $F_\HH := \phi(\HH)$ belongs to $\clFNpm.$
\begin{prop}\label{P: invariance principle}
  If $\HH_0=\HH_0^* \affl \clN$ is an operator with $\tau$-compact resolvent,
and if $\VV = \VV^* \in \clN,$ $\HH_1=\HH_0+\VV,$
then the following equality holds
  $$
    \xi_{\HH_1,\HH_0}(\lambda) = \xi_{F_{\HH_1}, F_{\HH_0}}(\phi(\lambda)).
  $$
\end{prop}
\begin{proof}
Let $h \in C^3_\csupp(-1,1)$ and $f(\lambda) = h(\phi(\lambda)).$
Then by Theorem \ref{T: trace formula}
  \begin{gather*}
    A:=\taubrs{f(\HH_1) - f(\HH_0)} = \int_\mbR f'(\lambda) \xi_{\HH_1,\HH_0}(\lambda)\,d\lambda
  \end{gather*}
  and since $F_{\HH_1} - F_{\HH_0} \in \clKN$ by~\cite[Lemma 2.7]{CP98CJM},
  we can apply Proposition \ref{P: xi function for pair F0 and F1} to get
  \begin{gather*}
    A = \taubrs{\h(F_{\HH_1}) - \h(F_{\HH_0})} = \int_{-1}^1 \h'(t)\xi_{F_{\HH_1},F_{\HH_0}}(t)\,dt
    \\ = \int_{-\infty}^\infty \h'(\phi(\lambda)) \phi'(\lambda) \xi_{F_{\HH_1},F_{\HH_0}}(\phi(\lambda))\,d\lambda
      = \int_{-\infty}^\infty f'(\lambda) \xi_{F_{\HH_1},F_{\HH_0}}(\phi(\lambda))\,d\lambda.
  \end{gather*}
  Since $f$ is an arbitrary $C^2$-function with compact support,
  comparing the last two formulas we get the equality
  \begin{gather}\label{F: xi DD = xi FF + const}
    \xi_{\HH_1,\HH_0}(\lambda) = \xi_{F_{\HH_1},F_{\HH_0}}(\phi(\lambda)) + C.
  \end{gather}
It is left to show that the constant $C=0.$

%
Let $h$ be a non-negative function from $C^\infty_\csupp(-1,1).$
By Lemma \ref{L: int ... d lambda = int ... dr} we have
\begin{gather}\label{F: int h d lambda = int ... dr (2)}
  \int_\mbR h(\phi(\lambda))\xi_{\HH_1,\HH_0}(\lambda)\,d\lambda
   = \int_0^1 \taubrs{\VV h(F_{\HH_r})}\,dr.
\end{gather}

Multiplying the first term of the RHS of (\ref{F: xi DD = xi FF + const}) by $h(\phi(\lambda)),$
integrating it and using Lemma \ref{L: int ... d lambda = int ... dr for bdd case},
we get
\begin{multline*}
  A := \int_\mbR h(\phi(\lambda)) \xi_{F_{\HH_1},F_{\HH_0}}(\phi(\lambda))\,d\lambda
   = \int_{-1}^1 h(\mu) \xi_{F_{\HH_1},F_{\HH_0}}(\mu) (\phi^{-1})'(\mu)\,d\mu
   \\ = \int_0^1 \taubrs{K h(F_r)(\phi^{-1})'(F_r)}\,dr,
\end{multline*}
where $K = F_{\HH_1} - F_{\HH_0}$ and
$F_r$ is the straight line path connecting $F_{\HH_1}$ and $F_{\HH_0}.$
Let $g(\mu) = h(\mu)(\phi^{-1})'(\mu).$
By Corollary \ref{C: integral is independent of path (bdd)} we have
$$
  A = \int_0^1 \taubrs{K g(F_r)}\,dr = \int_0^1 \taubrs{\dot F_{\HH_r}g(F_{\HH_r})}\,dr.
$$
By Proposition \ref{P: DalKreinII}(ii) we have $\dot F_{\HH_r} = T^{\HH_r,\HH_r}_{\phi^{[1]}}(\VV).$
Hence,
$$
  A = \int_0^1 \taubrs{T^{\HH_r,\HH_r}_{\phi^{[1]}}(\VV) g(F_{\HH_r})}\,dr.
$$
Using the BS-representation for $\phi^{[1]}$ given by (\ref{F: BS repr-n from ACDS}),
it follows from the definition of DOI (\ref{F: Def of MOI}),
Lemma \ref{L: C1plus is subset of L1}(ii) and Lemma \ref{L: lemma 3.10 ACDS}, that
\begin{multline}\label{F: RHS = int trace V h FDr phi'}
  A = \int_0^1 \taubrs{ \int_\Pi e^{i(s-t)\HH_r}\VV e^{it\HH_r}\,d\nu_\phi(s,t)
             \cdot g(F_{\HH_r})}\,dr
    \\ = \int_0^1 \int_\Pi \taubrs{ e^{i(s-t)\HH_r}\VV e^{it\HH_r}g(F_{\HH_r})}\,d\nu_\phi(s,t)\,dr
    \\ = \frac 1{\sqrt{2\pi}} \int_0^1 \int_\mbR \taubrs{ \VV e^{is\HH_r}is \hat\phi(s)
           g(F_{\HH_r})}\,ds\,dr
    \\ = \frac 1{\sqrt{2\pi}} \int_0^1 \taubrs{\VV g(F_{\HH_r}) \int_\mbR e^{is\HH_r}is \hat\phi(s)\,ds}\,dr
    \\ = \int_0^1 \taubrs{\VV g(F_{\HH_r}) \phi'(\HH_r)}\,dr
     = \int_0^1 \taubrs{\VV h(F_{\HH_r})}\,dr,
\end{multline}
since $g(\phi(\lambda)) \phi'(\lambda) = h(\phi(\lambda)).$
It follows from (\ref{F: xi DD = xi FF + const}),
(\ref{F: int h d lambda = int ... dr (2)}) and (\ref{F: RHS = int trace V h FDr phi'})
that $C \int_\mbR h(\phi(\lambda))\,d\lambda = 0$ and, hence, $C=0.$
\end{proof}

\subsection{The first formula for spectral flow} \label{SS: the first formula for sf}
We establish first
a spectral flow formula for bounded $\tau$-Fredholm operators. In this way
we avoid a number of difficulties with unbounded operators. Then we make
a `change of variable' to get to the unbounded case.

First we require some additional notation which is important for establishing
a convention for how we handle the situation when the endpoints have a kernel.
Let $a<0,$ $b>0$ and let $\sign_{a,b}$ be the function defined as
$\sign_{a,b}(x) = b$ if $x \geq 0,$ and $\sign_{a,b}(x) = a$ if $x < 0.$

We will write $\tF = \sign_{a,b}(F),$ when it is clear from the context what the numbers $a$ and $b$ are.
\begin{defn} If $F \in \clFNab$ and $\af$ is a $C^2$-function on $[0,\infty)$ vanishing in a neighbourhood of $0$
then for $h(\lambda) = \af\Brs{(b-\lambda)(\lambda-a)}$ we define $\gamma_\h(F)$ as
$$
  \gamma_\h(F) = \int_0^1 \alpha^\h_{F_\rr}(\dot F_\rr)\,d\rr,
$$
where $\alpha^\h$ is the closed one-form defined before Proposition \ref{P: one form tau(Xf(F)) is closed},
and $\set{F_\rr}_{\rr \in [0,1]}$ is the straight line connecting $F$ and $\tF.$
\end{defn}
The following theorem is the  analogue in our setting of~\cite[Theorem 5.7]{CP2}.
It is the fundamental formula that we need as our starting point.
\begin{thm}\label{T: sf = int h xi + gamma1 - gamma0}
  Let $F_0 \in \clFNab,$ let $K \in \clKN$ and let $F_1 = F_0 + K.$
  Let $\af$ be a $C^2$-function on $[0,\infty)$ vanishing in a neighbourhood of $0,$
such that the integral of $h(\lambda)=\af\Brs{(b-\lambda)(\lambda-a)}$ over $(a,b)$ is equal to $1.$
  Then the spectral flow between $F_0$ and $F_1$ is equal to
  $$
    \sflow(F_0,F_1) = \int_a^b \h(\lambda)\xi_{F_1,F_0}(\lambda)\,d\lambda
          + \gamma_\h(F_1) - \gamma_\h(F_0).
  $$
\end{thm}
\begin{proof} The proof follows ideas of~\cite[Theorem 5.7]{CP2}.
First of all by Phillips' definition of spectral flow~\cite{BCPRSW}
we have
  $$
    \sflow(F_0,F_1) = \sflow(\tF_0,\tF_1).
  $$
 (Note that, {\it  a priori} we would expect to write
$$\sflow(F_0,F_1) = \sflow(F_0,\tF_0) + \sflow(\tF_0,\tF_1) + \sflow(\tF_1,F_1),$$
 however there is no spectral flow along the paths joining $F$ and $\tF$ as noted in the proof
  of~\cite[Theorem 1.7]{CP98CJM}.)

  Now, by Lemma \ref{L: sf lemma for P and Q} we have
  $$
    \sflow(\tF_0,\tF_1) = \int_0^1 \alpha^\h_{\tF_\rr}(\dot \tF_\rr)\,d\rr,
  $$
  where $\{\tF_\rr\}_{\rr \in [0,1]}$ is the straight line path, connecting $\tF_0$ and $\tF_1.$
  By Corollary \ref{C: integral is independent of path (bdd)}
  we can replace this path by the (broken) path given on this diagram
  $$
    \xymatrix{
       F_0 \ar@{-->}[r] & F_1 \ar@{-->}[d]^{\gamma_\h(F_1)} \\
       \tF_0 \ar@{>}[r] \ar@{-->}[u]^{-\gamma_\h(F_0)} & \tF_1 \\
    }
  $$
  Then we get
  $$
    \sflow(\tF_0,\tF_1) =  - \gamma_\h(F_0) + \int_0^1 \alpha^\h_{F_\rr}(\dot F_\rr)\,d\rr + \gamma_\h(F_1),
  $$
  where $\set{F_\rr}_{\rr \in [0,1]}$ is the straight line path, connecting $F_0$ and $F_1.$
  But, setting $F_1 - F_0 = K,$ we have by Lemma \ref{L: int ... d lambda = int ... dr for bdd case}
  \begin{gather*}
    \int_0^1 \alpha^\h_{F_\rr}(\dot F_\rr)\,d\rr
         = \int_0^1 \taubrs{K \h(F_\rr)}\,d\rr
         =  \int_\mbR \h(\lambda)\xi_{F_1,F_0}(\lambda)\,d\lambda.
  \end{gather*}
\end{proof}
\begin{thm}\label{T: sf = int h xi + gamma1 - gamma0: 2}
  Let $F_0 \in \clFNpm,$ let $K \in \clKN$ and let $F_1 = F_0 + K.$
  Let $\af$ be a $C^2$-function on $[0,\infty)$ vanishing in a neighbourhood of $0,$
such that the integral of $h(\lambda)=\af(1-\lambda^2)$ over $(-1,1)$ is equal to $1.$
  Then the spectral flow function for the pair $F_0$ and $F_1$ is equal to
  $$
    \sflow(\mu; F_0,F_1) = \int_{-1}^1 \h(\lambda)\xi_{F_1,F_0}(\lambda)\,d\lambda
          + \gamma_{\h_{-\mu}}(F_1-\mu) - \gamma_{\h_{-\mu}}(F_0-\mu),
  $$
  where $\h_{-\mu}(\lambda) = \h(\lambda+\mu).$
\end{thm}
\begin{proof} By definition we have
$$
  \sflow(\mu; F_0,F_1) = \sflow(F_0-\mu,F_1-\mu).
$$
Since $F_j - \mu \in \FNab{-1-\mu}{1-\mu},$ by previous theorem we have
\begin{multline*}
  \sflow(F_0-\mu,F_1-\mu)
   \\  = \int_{-1-\mu}^{1-\mu} \h_{-\mu}(\lambda)\xi_{F_1-\mu,F_0-\mu}(\lambda)\,d\lambda
          + \gamma_{\h_{-\mu}}(F_1-\mu) - \gamma_{\h_{-\mu}}(F_0-\mu).
\end{multline*}
Since $\xi_{F_1-\mu,F_0-\mu}(\lambda) = \xi_{F_1,F_0}(\lambda+\mu),$ we have
\begin{multline*}
  \sflow(F_0-\mu,F_1-\mu)
  \\ = \int_{-1-\mu}^{1-\mu} \h(\lambda+\mu)\xi_{F_1,F_0}(\lambda+\mu)\,d\lambda
          + \gamma_{\h_{-\mu}}(F_1-\mu) - \gamma_{\h_{-\mu}}(F_0-\mu)
  \\ = \int_{-1}^{1} \h(\lambda)\xi_{F_1,F_0}(\lambda)\,d\lambda
          + \gamma_{\h_{-\mu}}(F_1-\mu) - \gamma_{\h_{-\mu}}(F_0-\mu).
\end{multline*}
\end{proof}
\begin{cor}
  If $F_0$ and $F_1$ are unitarily equivalent, then
  $$
    \sflow(\mu; F_0,F_1) = \xi_{F_1, F_0}(\mu) = \const.
  $$
\end{cor}
\begin{proof}
By Corollary \ref{C: if F0 and F1 are unit equiv then xi = const}
the function $\xi_{F_1,F_0}(\cdot)$ is constant on $(-1,1),$
so that
 $\int_\mbR \h(\lambda)\xi_{F_1,F_0}(\lambda)\,d\lambda = \xi_{F_1, F_0}(0).$

If $F_0$ and $F_1$ are unitarily equivalent, then
$\gamma_{\h_{-\mu}}(F_1-\mu) = \gamma_{\h_{-\mu}}(F_0-\mu).$
Hence, for all $\mu \in (-1,1)$
  $$
    \sflow(\mu; F_0,F_1) = \xi_{F_1, F_0}(\mu) = \xi_{F_1, F_0}(0).
  $$
\end{proof}
\begin{lemma}\label{L: gamma = xi}
If $F \in \clFNpm$ and if
$\set{\h_\eps}_{\eps>0}$ is an approximate $\delta$ function  (by compactly supported even functions)
 then for all $\mu \in (-1,1)$ the limit
$$
  \gamma_{\mu}(F) := \lim_{\eps \to 0} \gamma_{\h_\eps}(F-\mu)
$$
exists and is equal to $\xi_{G,\tG}(0),$ where $G = F - \mu.$
\end{lemma}
\begin{proof}
Since $h$ is an even function we have that
  \begin{gather*}
     \int_{-\infty}^0 {\h_\eps}(\lambda) \xi_{G,\tG}(\lambda)\,d\lambda
     \to \frac 12 \xi_{G,\tG}(0-)
  \end{gather*}
and
  \begin{gather*}
     \int_0^\infty {\h_\eps}(\lambda) \xi_{G,\tG}(\lambda)\,d\lambda
     \to \frac 12 \xi_{G,\tG}(0+),
  \end{gather*}
  as $\eps \to 0.$
If $\set{G_\rr}_{\rr \in [0,1]}$ is the straight line path connecting $G$ and $\tG$ then
by Lemma \ref{L: int ... d lambda = int ... dr for bdd case} we have
  \begin{gather*}
     \gamma_{\h_\eps}(G)
        = \int_0^1 \alpha_{\h_\eps}(\dot G_\rr)\,d\rr
        = \int_0^1 \taubrs{\dot G_\rr {\h_\eps}(G_\rr)}\,d\rr
     \\ = \int_\mbR {\h_\eps}(\lambda) \xi_{G,\tG}(\lambda)\,d\lambda
     \to \frac 12 \brs{ \xi_{G,\tG}(0-) + \xi_{G,\tG}(0+)}
     = \xi_{G,\tG}(0)
  \end{gather*}
 as $\eps \to 0,$ by Definition \ref{D: xi at discontinuity point} of $\xi$ at
 discontinuity points.
\end{proof}

Now we need to handle the situation when the endpoints are not unitarily equivalent.
For this we require some additional facts about the `end-point
correction terms'. The interesting fact which we now establish
is that the approach using the spectral shift function differs in a fundamental way from
the previous point of view in~\cite{CP2}. The next few results demonstrate this by showing that
the spectral shift function absorbs the contribution to the formula due to the
spectral asymmetry of the endpoints leaving only kernel correction terms to be handled.

\begin{lemma}\label{L: gamma = trace of ker}
If $F \in \clFNpm$ and if $\mu \in (-1,1),$ then
the following equality holds true
  $$
    \gamma_{\mu}(F) = \frac 12 \tau [\ker (F - \mu)].
  $$
\end{lemma}
\begin{proof} Let $G = F - \mu.$
We have
$$
  \tau [\ker G] = \taubrs{E^G_{(-\infty,0]} - E^\tG_{(-\infty,0]}}.
$$
By Proposition \ref{P: xi function for pair F0 and F1}(ii) and Definition \ref{D: def of xi F F}
the value $\xi_{\tG,G}(0)$ is the half sum of the last expression and
$$
  \taubrs{E^G_{(-\infty,0)} - E^\tG_{(-\infty,0)}} = 0.
$$
Hence, by Lemma \ref{L: gamma = xi}
$$
  \gamma_{\mu}(F) = \xi_{G,\tG}(0) = \frac 12 \taubrs{[\ker G]}.
$$
\end{proof}
\begin{thm}\label{T: sf = xi + gamma1 - gamma0} If $F_0,F_1 \in \clFNpm$ such that $F_1 - F_0 \in\clKN,$
then for all $\mu \in (-1,1)$
  \begin{gather}\label{F: sf = xi + gamma1 - gamma0}
    \sflow(\mu; F_0,F_1) = \xi_{F_1,F_0}(\mu) + \frac 12 (\tau [\ker (F_1 - \mu)] - \tau [\ker (F_0 - \mu)]).
  \end{gather}
\end{thm}
\begin{proof} Replace $h$ in Theorem \ref{T: sf = int h xi + gamma1 - gamma0: 2} by $h_{\eps,\mu}$
(thus translate the approximate $\delta$ function  $h_{\eps}$ by $\mu$) and then
let $\eps \to 0$ using Lemmas
\ref{L: gamma = xi}, \ref{L: gamma = trace of ker}.
\end{proof}

We now see that under hypotheses that guarantee both are defined
the spectral flow function and the spectral shift function differ only by kernel corrections terms for the endpoints.  We should remark
that the occurrence of the correction terms $\gamma_\mu(F_j), \ j=1,2,$ in the last formula can be explained by the fact that
we actually define the spectral flow function and the spectral shift function at discontinuity
points in different ways. The spectral shift function is defined as a half-sum of the left and the right limits,
while the spectral flow is defined to be left-continuous.

\subsection{Spectral flow in the unbounded case}
The formulae for spectral flow in the bounded case may be now used to establish corresponding
results in our original setting of unbounded self adjoint  operators with compact resolvent.

By Proposition \ref{P: invariance principle} $\xi_{\HH_1,\HH_0}(0) = \xi_{F_{\HH_1},F_{\HH_0}}(0)$
and by definition of spectral flow for unbounded operators~\cite{BCPRSW}
$\sflow(\HH_0,\HH_1) = \sflow (F_{\HH_1},F_{\HH_0}).$ Hence, it follows from
(\ref{F: sf = xi + gamma1 - gamma0}) taken at $\mu = 0$ that
  $$
    \sflow(\HH_0,\HH_1) = \xi_{\HH_1,\HH_0}(0) + \gamma_0(F_1) - \gamma_0(F_0).
  $$
  Since
  $$
    \ker(\HH) = \ker(F_\HH)
  $$
  we have the following equality
  \begin{gather*}\label{F: sf = xi - trace of ker + trace of ker at 0}
    \sflow(\HH_0,\HH_1) = \xi_{\HH_1,\HH_0}(0) + \frac 12 \tau[\ker (\HH_1)] - \frac 12 \tau[\ker (\HH_0)].
  \end{gather*}
If we replace here the operators $\HH_0$ and $\HH_1$ by the operators $\HH_0-\lambda$ and $\HH_1-\lambda$
respectively then we get
  \begin{gather*}
    \sflow(\lambda; \HH_0,\HH_1) = \xi_{\HH_1-\lambda,\HH_0-\lambda}(0)
      + \frac 12 \tau[\ker (\HH_1-\lambda)] - \frac 12 \tau[\ker (\HH_0-\lambda)].
  \end{gather*}
  Since $\xi_{\HH_1-\lambda,\HH_0-\lambda}(0) = \xi_{\HH_1,\HH_0}(\lambda)$
  it follows that
  \begin{gather}\label{F: sf = xi - trace of ker + trace of ker}
    \sflow(\lambda; \HH_0,\HH_1) = \xi_{\HH_1,\HH_0}(\lambda)
      + \frac 12 \tau[\ker (\HH_1-\lambda)] - \frac 12 \tau[\ker (\HH_0-\lambda)].
  \end{gather}

\subsection{The spectral flow formula using infinitesimal spectral flow}
The results on the spectral shift function that we established in Section
\ref{S: spectral shift measure} now suggest a new direction for spectral
flow theory.
\begin{defn} Let $\HH_0$ be a self-adjoint operator affiliated with $\clN$ having $\tau$-compact resolvent.
The \emph{infinitesimal spectral flow one-form} is a distribution-valued one-form $\Phi_{\HH}$
on the affine space $\HH_0 + \clN_{sa},$ defined by formula
$$
   \la \Phi_\HH(X), \phi \ra = \taubrs{X\phi(\HH)}, \qquad X \in \clN_{sa}, \quad \phi \in C^\infty_\csupp(\mbR).
$$
Formally,
$$
  \Phi_\HH(X) = \taubrs{X \delta(\HH)},
$$
where $\delta(\HH)$ is the $\delta$-function of $\HH.$
\end{defn}

We believe that the notion of infinitesimal spectral flow will have application to the problem of
studying spectral flow when the endpoints differ by a $\HH_0$ relatively bounded perturbation.
To this end we establish that spectral flow may be reformulated in terms of it.

\begin{thm} Let $\HH_1 \in \HH_0 + \clN_{sa}.$ Spectral flow between $\HH_0$ and $\HH_1$ is equal to the integral
of  the infinitesimal spectral flow one-form along any piecewise $C^1$-path $\set{\HH_\rr}_{\rr \in [0,1]}$
in $\HH_0+\clN$ connecting $\HH_0$ and $\HH_1$
in the sense that for any $\phi \in C^\infty_\csupp(\mbR)$ the following equality holds
true
$$
  \int_\mbR \sflow(\lambda; \HH_0, \HH_1)\phi(\lambda) \,d\lambda = \int_0^1 \la \Phi_{\HH_\rr}(\dot \HH_\rr), \phi \ra\,d\rr.
$$
Formally,
$$
  \sflow(\HH_0, \HH_1) = \int_0^1 \Phi_{\HH_\rr}(\dot \HH_\rr)\,d\rr,
\qquad \text{or} \qquad
  \sflow(\lambda; \HH_0, \HH_1) = \int_0^1 \Phi_{\HH_\rr-\lambda}(\dot \HH_\rr)\,d\rr.
$$
\end{thm}
\begin{proof} By Corollary \ref{C: integral is independent of path}
we can choose the path $\set{\HH_\rr}_{\rr \in [0,1]}$
to be the straight line path $\HH_\rr = \HH_0 + \rr\VV$.
It follows from Lemmas \ref{L: D+V has comp resolvent} and \ref{L: If D has comp rslv then E Delta if finite}
that the functions $\lambda \mapsto \taubrs{[\ker(\HH_0-\lambda)]}$ and $\lambda \mapsto \taubrs{[\ker(\HH_1-\lambda)]}$
can be non-zero only on a countable set.
%
Hence, by (\ref{F: sf = xi - trace of ker + trace of ker})
and Lemma \ref{L: int ... d lambda = int ... dr} we have
  \begin{multline*}
    \int_\mbR \sflow(\lambda; \HH_0, \HH_1)\phi(\lambda) \,d\lambda
     = \int_\mbR \xi_{\HH_1,\HH_0}(\lambda) \phi(\lambda) \,d\lambda
     \\ = \int_0^1 \taubrs{\VV \phi(\HH_\rr)}\,d\rr
     = \int_0^1 \la \Phi_{\HH_\rr}(\dot \HH_\rr), \phi \ra\,d\rr.
  \end{multline*}
\end{proof}
We remark that the infinitesimal spectral flow one-form is
exact in the sense that its value on every test function is exact.
\subsection{The spectral flow formulae in the $\clI$-summable spectral triple case}
The original approach of~\cite{CP2} required summability constraints on
the operator $\HH_0$. We will now see that if indeed
 $\HH_0$ satisfies such conditions then we can weaken conditions on the function $f$
in Theorem \ref{T: tau(V f(D))}.
\begin{lemma}\label{L: 1-norm of f(D+V) for non-comp supp f}
 Let $\HH_0$ be a self-adjoint operator with $\tau$-compact resolvent
affiliated with~$\clN.$ Let~$g$ be an increasing
continuous function on $[0,+\infty),$ such that $g(0) \geq 0$ and
$g\brs{ c (1+\HH^2)^{-1}} \in \LpN{1}$ for all $c > 0.$
Let $f(x) = g\brs{(1+x^2)^{-1}}.$
Then for any $R > 0$ and for any $\VV=\VV^* \in \clN$
the operator $f(\HH+\VV)$ is trace class and the function
$$
  \VV \in B_R \mapsto \norm{f(\HH+\VV)}_1
$$
is bounded.
\end{lemma}
\begin{proof} By~\cite[Lemma 2.5(iv)]{FK86PJM} we have for all $t>0$
  \begin{gather*}
    \mu_t(f(\HH+\VV)) = \mu_t\Big(g\brs{\brs{1 + (\HH+\VV)^2}^{-1}}\Big)
    = g\Big(\mu_t\brs{\brs{1 + (\HH+\VV)^2}^{-1}}\Big).
  \end{gather*}
By Lemma \ref{L: lemma with def of B R} there exists a constant $c=c(R)>0$ such that
for any $\VV \in B_R$
$$
  \brs{1 + (\HH+\VV)^2}^{-1} \leq  c \brs{1 + \HH)^2}^{-1}.
$$
Hence, by~\cite[Lemma 2.5(iii)]{FK86PJM} we have
$$
  \mu_t(f(\HH+\VV)) \leq g\brs{\mu_t\sqbrs{c \brs{1 + \HH^2}^{-1}}} = \mu_t\brs{g\sqbrs{c \brs{1 + \HH^2}^{-1}}}.
$$
Since $g\brs{ c (1+\HH^2)^{-1}} \in \LpN{1},$ the last function belongs to $L^1[0,\infty),$ which implies that
$f(\HH+\VV) \in \LpN{1}.$
\end{proof}
\begin{lemma}\label{L: alpha h eps is exact}
  Let $\HH_0,$ $g$ and $f$ be as in Lemma \ref{L: 1-norm of f(D+V) for non-comp supp f}.
  An integral of the one-form
  $$
    \alpha^f_\HH(X) = \tau(X f(\HH)), \ X \in \clN, \ \HH \in \HH_0+\clN_{sa},
  $$
  along a piecewise smooth path in $\HH_0+\clN_{sa}$ depends only on endpoints of that path.
\end{lemma}
\begin{proof} Let $f_n$ be a increasing sequence of compactly supported smooth functions
  converging pointwise to $f$ and $\Gamma_1, \ \Gamma_2$ be two piecewise smooth paths in
  in $\HH_0+\clN_{sa}$ with the same endpoints. Then by Lemma \ref{L: 1-norm of f(D+V) for non-comp supp f},
  Lebesgue dominated convergence theorem
  and Corollary \ref{C: integral is independent of path} we have
  \begin{multline*}
    \int_{\Gamma_1} \alpha^f
      = \int_{\Gamma_1} \liminfty{n} \alpha^{f_n}
      = \liminfty{n} \int_{\Gamma_1} \alpha^{f_n}
      \\ = \liminfty{n} \int_{\Gamma_2} \alpha^{f_n}
      = \int_{\Gamma_2} \liminfty{n} \alpha^{f_n}
      = \int_{\Gamma_2} \alpha^{f}.
  \end{multline*}
\end{proof}
The condition that $g(c(1+\HH_0^2)^{-1})$ be trace class
is a generalised summability constraint. This notion arises naturally for certain ideals
$\clI$ of compact operators (for example
for the Schatten ideals $\clI_p$, $p\geq 1$,
$g(x)=x^{p/2}$ and we have the notion of $p$-summability). We have also already remarked on the  $\theta$-summable case.

Now if there is a unitary $u\in \clN$ with $\VV=u^*[\HH_0, u]$ bounded then
we have, for a dense subalgebra $\clA$ of the $C^*$-algebra generated by
$u$, a semifinite `$g$-summable' spectral triple $(\clA, \clN, \HH_0)$.
Moreover $\HH_0+\VV=u\HH_0 u^*$ so we have unitarily equivalent endpoints.

\begin{thm}\label{T: I-summable case} Let $f$ be a non-negative $L^1$-function such that $f(\HH_\rr) \in \LpN{1}$ for all $\rr \in [0,1],$
and let $\rr \mapsto \norm{f(\HH_\rr)}_1$ be integrable on $[0,1].$
If $\HH_0$ and $\HH_1$ are unitarily equivalent then
$$
  \sflow(\lambda; \HH_0,\HH_1) = C^{-1} \int_0^1 \taubrs{\VV f(\HH_\rr - \lambda)}\,d\rr,
$$
where $C = \int_{-\infty}^\infty f(\lambda)\,d\lambda.$
\end{thm}
\begin{proof} Unitary equivalence of $\HH_0$ and $\HH_1$ implies
that two last terms in (\ref{F: sf = xi - trace of ker + trace of ker}) vanish.
In case of $f \in B_\csupp(\mbR),$ multiplying
(\ref{F: sf = xi - trace of ker + trace of ker}) by $f(\lambda)$
and integrating it we get the required equality
by Lemma \ref{L: int ... d lambda = int ... dr} and
Theorem \ref{T: SSF for unitary equivalent endpoints}.
For an arbitrary $f \in L^1$ the claim follows
from Lebesgue's dominated convergence theorem
by approximating
$f$ by an increasing sequence of step-functions converging a.e. to $f.$
\end{proof}
The following corollary recovers two of the main results of~\cite{CP98CJM,CP2}.
\begin{cor} (i) If $\HH_0$ is $\theta$-summable with respect to $\clN$
and if $\HH_0$ and $\HH_1$ are unitarily equivalent then
$$
  \sflow(\HH_0,\HH_1) = \sqrt{\frac \eps\pi} \int _0^1 \taubrs{\VV e^{-\eps\HH_\rr^2}}\,d\rr.
$$

(ii) If $\HH_0$ is $p$-summable (i.e. $(1+\HH_0^2)^{-p/2} \in \LpN{1}$)
with respect to $\clN,$ where $p > 1$ and if $\HH_0$
and $\HH_1$ are unitarily equivalent then
$$
  \sflow(\HH_0,\HH_1) = C_p^{-1}  \int _0^1 \taubrs{\VV \brs{1+\HH_\rr^2}^{-\frac p2}}\,d\rr,
$$
where $C_p = \int_{-\infty}^\infty (1 + \lambda^2)^{-\frac p2}\,d\lambda.$
\end{cor}
\begin{proof} Put $f(\lambda) = e^{-\eps \lambda^2}$ and $f(\lambda) = (1 + \lambda^2)^{-\frac p2}$ for (i) and (ii) respectively
in Theorem \ref{T: I-summable case}. The conditions of that theorem are fulfilled by Lemma
\ref{L: 1-norm of f(D+V) for non-comp supp f}.
\end{proof}

\subsection{Recovering $\eta$-invariants}
To demonstrate that we have indeed generalized previous analytic approaches to
spectral flow formulae we still need some refinements. What is missing is the relationship of the
`end-point correction terms'
to the truncated eta invariants of~\cite{Ge93Top}.

In fact Theorem \ref{T: sf = int h xi + gamma1 - gamma0} combined  with some ideas of~\cite{CP2}
will now enable us to give a new proof of the original formula  (\ref{F: Carey-Phillips})
for spectral flow with unitarily inequivalent endpoints.

Introduce the function
  $$
    \af_\eps(\lambda) = \sqrt{\frac \eps\pi} \lambda^{-3/2} e^{\eps(1-\lambda^{-1})}.
  $$
and let $h_\eps(\lambda) = \af_\eps(1-\lambda^2),$
$f_\eps(\lambda) = \af_\eps\brs{(1+\lambda^2)^{-1}}.$
\begin{lemma} \label{L: 3.25} Let $\HH_0=\HH_0^* \affl \clN$ be $\theta$-summable,
let $\VV \in \clN_{sa}$ and let $\HH_1 = \HH_0 + \VV.$ Then
$$
  \int_{-1}^1 h_\eps(\lambda) \xi_{F_{\HH_1},F_{\HH_0}}(\lambda)\,d\lambda = \sqrt{\frac \eps\pi}\int_0^1 \taubrs{\VV e^{-\eps\HH_r^2}}\,dr.
$$
\end{lemma}
\begin{proof} Since $h_\eps(\phi(\mu)) = f_\eps(\mu),$ by
  Proposition \ref{P: invariance principle} we have
  \begin{gather*}
    (A) := \int_{-1}^1 h_\eps(\lambda) \xi_{F_{\HH_1},F_{\HH_0}}(\lambda)\,d\lambda
        = \int_{-\infty}^\infty h_\eps(\phi(\mu)) \xi_{F_{\HH_1},F_{\HH_0}}(\phi(\mu))\phi(\mu)'\,d\mu
       \\ = \int_{-\infty}^\infty f_\eps(\mu) \phi'(\mu) \xi_{\HH_1,\HH_0}(\mu)\,d\mu.
  \end{gather*}
  Further, by Lemmas \ref{L: Fubini for variable measures} and \ref{L: trace of V alpha(D)}
  \begin{multline*}
     (A) = \int_{-\infty}^\infty f_\eps(\mu) \phi'(\mu) \int_0^1 \taubrs{\VV dE_\mu^{\HH_r}}\,dr
       = \int_0^1 \taubrs{\VV f_\eps(\HH_r) \phi'(\HH_r)}\,dr
      \\ = \int_0^1 \taubrs{\VV \sqrt{\frac\eps\pi} (1+\HH_r^2)^{3/2} e^{-\eps\HH_r^2}(1+\HH_r^2)^{-3/2}}\,dr
      \\  =  \sqrt {\frac \eps \pi}\int_0^1 \taubrs{\VV e^{-\eps\HH_r^2}}\,dr.
  \end{multline*}
\end{proof}
As we have emphasized previously, the strategy of our proof follows that of~\cite{CP2}
in that, we deduce the unbounded version of the spectral flow formula for the theta
summable case from a bounded version. To this end introduce $F_s = \HH (s + \HH^2)^{-1/2}.$
\begin{lemma}\cite[Lemma 8.8]{CP2} \label{L: CP2 Lemma 8.8}
We have
$$
  \lim_{\delta \to 0} \int_{\Gamma_\delta} \alpha^{h_\eps} = 0,
$$
where $\Gamma_\delta$ is the straight line connecting $F_0$ and $F_\delta.$
\end{lemma}

\begin{lemma}\label{L: gamma eps = eta eps + trace ker D} If $\HH=\HH^* \affl \clN$ is $\theta$-summable,
 then the following equality holds true
$$
  \gamma_{h_\eps}(F_\HH) = \frac 12 \brs{\eta_\eps(\HH) + \tau[\ker \HH]}.
$$
\end{lemma}
\begin{proof} We note that $1-F_s^2 = s(s+\HH^2)^{-1}$ and that $\dot F_s = -\frac 12 \HH(s+\HH^2)^{-3/2}.$
The path $\Gamma_1 := \set{F_s}_{s \in [0,1]}$ connects $\sgn(F_\HH)$ with $F_\HH.$
If we denote by $\Gamma_2$
the straight line path connecting $\sgn(F_\HH)$ with $\tF=\sign(F_\HH)$ then
the path $-\Gamma_1+\Gamma_2$ connects $F_\HH$ with $\tF,$ so that by
Lemma \ref{L: alpha h eps is exact} applied to $f = h_\eps,$ and by the argument of~\cite{CP2} and
Lemma \ref{L: CP2 Lemma 8.8} dealing with discontinuity of the path $\Gamma_1$ at zero, it follows that
$$
  \gamma_{h_\eps}(F_\HH) = - \int_{\Gamma_1} \alpha^{h_\eps} + \int_{\Gamma_2} \alpha^{h_\eps}.
$$
We have for the first summand
  \begin{multline*}
    \int_{\Gamma_1} \alpha^{h_\eps} = \int_0^1 \alpha^{h_\eps}_{F_s}(\dot F_s)\,ds = \int_0^1 \taubrs{\dot F_s {h_\eps}(F_s)}\,ds
    \\ = - \frac 12 \int_0^1 \taubrs{\HH (s+\HH^2)^{-3/2} \af_\eps\brs{1-F_s^2}}\,ds
    \\ = - \frac 12 \int_0^1 \taubrs{\HH (s+\HH^2)^{-3/2} \af_\eps\brs{s(s + \HH^2)^{-1}}}\,ds
    \\ = - \frac 12 \sqrt{\frac \eps\pi} \int_0^1 \taubrs{\HH (s+\HH^2)^{-3/2} s^{-3/2}(s+\HH^2)^{3/2} e^{-\frac \eps{s} \HH^2}}\,ds
    \\ = - \frac 12 \sqrt{\frac \eps\pi} \int_0^1 \taubrs{\HH e^{-\frac \eps{s} \HH^2}} s^{-3/2} \,ds
     \\ =  - \frac 12 \sqrt{\frac \eps\pi} \int_1^\infty \taubrs{\HH e^{- \eps t \HH^2}} \frac {dt}{\sqrt{t}}
      = -\frac 12 \eta_\eps(\HH).
  \end{multline*}

Let $E = [\ker \HH]$ and let $G_r = \sgn(F_\HH) + rE$ be the path $\Gamma_2.$
Then the second summand is equal to
  \begin{gather*}
    \int_{\Gamma_2} \alpha^{h_\eps} = \int_0^1 \taubrs{\dot G_r h_\eps(G_r)}\,dr
     = \int_0^1 \taubrs{E \af_\eps(1 - G_r^2)}\,dr.
  \end{gather*}
Since $1 - G_r^2 = E (1-r^2),$ it follows that
  \begin{gather*}
    \int_{\Gamma_2} \alpha^{h_\eps}
     = \int_0^1 \taubrs{E \af_\eps(1 - r^2)}\,dr = \taubrs{E} \int_0^1 \af_\eps(1 - r^2)\,dr = \frac 12 \tau(E),
  \end{gather*}
so that $\gamma_{h_\eps}(F_\HH) = \frac 12 \brs{\eta_\eps(\HH) + \tau[\ker \HH]}.$
\end{proof}

As a direct corollary of these Lemmas and Theorem \ref{T: sf = int h xi + gamma1 - gamma0}
we get a new proof of (\ref{F: Carey-Phillips}).
The idea used in the proof of the next theorem, of approximating the
exponential function by functions of compact support,
has been exploited in the context of spectral flow formula in
\cite{Wa}.

\begin{thm} If $\HH_0$ is $\theta$-summable then the formula\, {\rm (\ref{F: Carey-Phillips})} holds true.
\end{thm}
\begin{proof} Let $h_n$ be a sequence of smooth non-negative functions, compactly supported on $(-1,1)$, and
converging pointwise to $\h_\eps.$
Recall that $\HH_r=\HH_0+rV$. Then the sequence $\gamma_{h_n}(F_{\HH_r})$ converges to
$\gamma_{h_\eps}(F_{\HH_r})$ and the sequence $\int_{-1}^1 h_n(\lambda) \xi_{F_{\HH_1},F_{\HH_0}}(\lambda)\,d\lambda$
converges to $\int_{-1}^1 h_\eps(\lambda) \xi_{F_{\HH_1},F_{\HH_0}}(\lambda)\,d\lambda$
by Lebesgue's DCT, since $\theta$-summability of $\HH_0$ implies 1-summability of $h_\eps(F_{D_r}).$
Hence the claim follows from Theorem \ref{T: sf = int h xi + gamma1 - gamma0},
Lemma \ref{L: 3.25}
and Lemma \ref{L: gamma eps = eta eps + trace ker D}.
\end{proof}



\mathsurround 0pt
\ndef{\AndSoOn}{$\dots$}

\end{document}